\documentclass[letterpaper,11pt]{article}
\usepackage[T1]{fontenc}
\usepackage{lmodern}
\usepackage{amsmath,amsthm,amsfonts,amssymb,mathrsfs,bm}
\usepackage[usenames]{color}
\usepackage{hyperref}
\usepackage{amssymb}
\usepackage{dsfont} 
\usepackage{graphicx}
\usepackage{subcaption}
\usepackage{algorithm}
\usepackage{algpseudocode}
\usepackage{float}

\usepackage[a4paper,margin=1in]{geometry}

\newtheorem{theorem}{Theorem}[section]
\newtheorem{lemma}[theorem]{Lemma}
\newtheorem{corollary}[theorem]{Corollary}
\newtheorem{proposition}[theorem]{Proposition}

\newtheorem{remark}{Remark}[section]

\usepackage{amsmath}
\usepackage{amssymb}
\usepackage{amsthm}
\usepackage{enumitem}
\usepackage{nicematrix}
\usepackage{booktabs}
\usepackage{mathrsfs}

\newcommand{\R}{\mathbb{R}}

\newcommand{\conv}{\textrm{conv}}

\newcommand{\PP}{\mathcal{P}}
\newcommand{\muhat}{\hat{\mu}}
\newcommand{\norm}[1]{\lVert #1 \rVert}
\newcommand{\abs}[1]{\lvert #1 \rvert}
\newcommand{\Lhat}{\hat{L}}
\newcommand{\eps}{\varepsilon}

\newcommand{\KL}{\mathrm{KL}}

\newcommand{\supp}{\mathrm{supp}}

\makeatletter
\renewcommand*{\@cite@ofmt}{\hbox}
\makeatother

\begin{document}
\title{Certified Algorithms for Gaussian Location Mixture NPMLEs via Discretized Gradient Flows}
\author{
	\small\begin{tabular}{c}
		{Pranav Gupta \thanks{Corresponding author}}\\
        Division of Applied Mathematics\\
		 Brown University\\ 
         pranav\_gupta@brown.edu
	\end{tabular}
    \and
  \small\begin{tabular}{c}
		{Satyaki Mukherjee}
		 \\  Department of Mathematics\\
		 National University of Singapore\\
         satyaki.mukhyo@gmail.com
	\end{tabular}
}
\date{}
\maketitle

\begin{abstract}
We study the computation of approximate nonparametric maximum likelihood
estimators (NPMLEs) for Gaussian location mixtures in arbitrary dimension
$d\ge 1$ via the Fisher--Rao gradient flow. We provide a certified
algorithm based on an entropic proximal discretization of the flow which,
given a sufficiently small $\varepsilon > 0$, returns an absolutely continuous mixing
measure $\tilde\mu$ admitting an exact finite-dimensional representation
and satisfying $\hat L_n-L_n(f_{\tilde\mu})\le\varepsilon$.
For almost every dataset $X_1,\cdots, X_n \in \mathbb{R}^d$ satisfying $n>2(d+1)(m+1)^2$, where $m$ is the
number of support points of the NPMLE, the algorithm requires
$O(n\varepsilon^{-1}\log^{d+2}(1/\varepsilon))$ time.
En route, we also establish results on the well-posedness and convergence
properties of the continuous-time Fisher--Rao flow. In particular,
initializing the flow at the uniform probability measure on a ball
containing the convex hull of the observations, we show
that $\hat L_n-L_n(f_{\mu_t})=O(\log t/t)$.
We also establish that, for almost every dataset satisfying the same condition on $n$ as above, the set of global maxima of the optimal score
is the same as the support of the NPMLE, and the Hessian of the optimal
score is nondegenerate at these points, extending a result of
\cite{polyanskiy2025nonparametric} to higher dimensions.
\end{abstract}

\tableofcontents
\newpage

\section{Introduction}
\subsection{Background and Overview}

Given any probability measure $\mu \in \mathcal{P}(\mathbb{R}^d)$, we define the corresponding Gaussian location mixture as a probability measure with density $f_{\mu}$ given by
\begin{equation}
\label{eq:intro-mixture}
 f_\mu(x):=\int_{\mathbb R^d}\phi(x-\theta)\,\mu(d\theta),
 \ \  \phi(x):=(2\pi)^{-d/2}\exp\!\left(-\frac{\|x\|^2}{2}\right).
\end{equation}
We call $\mu$ the mixing measure for this Gaussian location mixture. Given a dataset $X_1,\ldots,X_n \in \mathbb{R}^d$, we consider the log-likelihood and its supremum
\begin{equation}\label{eq:intro-likelihood}
 L_n(f_\mu):=\frac1n\sum_{i=1}^n\log f_\mu(X_i),
 \quad \widehat L_n:=\sup_{\mu\in\mathcal P(\mathbb R^d)}L_n(f_\mu).
\end{equation}
The nonparametric maximum likelihood estimator (NPMLE) $\hat{\mu}$ is defined to be a probability measure that attains this supremum. Starting from the works of Robbins~\cite{robbins1950generalization}
and Kiefer and Wolfowitz~\cite{kiefer1956consistency}, its statistical and
geometric properties have been extensively investigated
\cite{lindsay1983geometry,lindsay1995mixture,
zhang2009generalized, chen2017consistency}.

Gaussian location mixtures arise naturally in empirical Bayes
estimation, where the mixing measure represents the unknown
distribution of latent means observed with Gaussian noise
\cite{jiang2009general,efron2011tweedie,kim2026smooth}.
NPMLE-based methods have been studied for Gaussian denoising
\cite{saha2020nonparametric} and high-dimensional classification
\cite{dicker2016classification}.
Extensions to heteroscedastic Gaussian mixtures have also
been applied to astronomical data and hierarchical regression
models in the social and biological sciences
\cite{soloff2025multivariate}.

Another line of work concerns the self-regularizing property of the NPMLE. Polyanskiy and Wu~\cite{polyanskiy2020selfregularizing} showed that, for one-dimensional Gaussian location mixtures with subgaussian mixing distributions, the NPMLE has $O(\log n)$ support points with high probability. More recently, Han, Shen and Wu~\cite{han2026selfregularization} extended self-regularization to arbitrary dimensions for NPMLEs restricted to a fixed hypercube, establishing a high probability bound of $O((\log n/\log\log n)^d)$ on the number of support points. When the true mixing distribution has finite support, they further show that the number of support points is $O_{\mathbb P}(1)$. 

The computation of Gaussian-mixture
NPMLEs has been extensively studied in the literature. Classical and modern approaches include EM-type procedures and
convex optimization schemes on finite grids
\cite{jiang2009general,koenker2014convex,
feng2018approximate,kim2020fast,zhang2024efficient}.
Soloff, Guntuboyina and Sen~\cite{soloff2025multivariate}
also established quantitative bounds on the discretization
error relative to the unrestricted NPMLE. Although such methods can be highly effective
in practice, an optimization guarantee for the discretized problem
must be distinguished from an \textit{a priori} certificate relative to the unrestricted
optimum $\widehat L_n$. In particular, the existing general-purpose
algorithmic analyses do not give an \textit{a priori} certificate, for
arbitrary fixed data in $\mathbb{R}^d$, of the form
\begin{equation}
\label{certificate}
    \widehat L_n-L_n(f_{\mu_{\rm out}})\leq\varepsilon
\end{equation}
together with an explicit arithmetic complexity specified in terms of $\varepsilon$. An efficient algorithm of this form is of fundamental interest. For instance, the recent results in \cite{ghosh2026gaussian} show that for $X_1,\ldots,X_n$ sampled from a Gaussian location mixture $f_\star$
whose mixing distribution is supported on a fixed compact set, if
$\widetilde f_n$ satisfies $L_n(\widetilde f_n) \geq \widehat L_n-\varepsilon_n,$
then, with high probability,
\begin{equation}
\label{eq:intro-ghosh-hellinger}
    H^2(f_\star,\widetilde f_n)
    \leq
    \varepsilon_n
    +
    C\frac{(\log n)^{d+1}}{n}.
\end{equation}

The main contribution of this paper is Algorithm~\ref{alg:fa-adaptive} which, given any dataset $X_1,\ldots,X_n \in \mathbb{R}^d$ and a prescribed $\varepsilon>0$, returns a measure $\mu_{\rm out}$ satisfying \eqref{certificate}. The algorithm has an arithmetic complexity of
\begin{equation}
\label{eq:intro-complexity}
O\left(
n\varepsilon^{-1}\log^{d+2}(1/\varepsilon)
\right)
\end{equation}
for sufficiently small $\varepsilon$. In particular, we achieve a near linear dependence on $\varepsilon^{-1}$. 

A recent development in certified NPMLE computation is due to Polyanskiy and Sellke~\cite{polyanskiy2025nonparametric}, who study Gaussian location mixtures in one dimension. Their algorithm identifies the number of support points and computes an approximation to the NPMLE in Wasserstein distance in $K(X)+O\left(nm^2\log\log(1/\varepsilon)\right)$, where $m$ is the number of support points of the NPMLE, and $K(X)$ denotes an initialization cost. To the best of our knowledge, for $d>1$, ours is the first algorithm that provides an \textit{a priori} certificate relative to the unrestricted NPMLE likelihood, together with an explicit arithmetic complexity bound in terms of $\varepsilon$.

Our analysis starts with considering optimization directly on the space $\mathcal{P}(\mathbb{R}^d)$ instead of restricting to a discrete subset. Chizat~\cite{chizat2022convergence} developed a general theory for convergence rates of Bregman gradient methods on spaces of measures, including
a suboptimality gap of $O(\log k/k)$  under basic regularity
assumptions. Related convergence results for mirror descent on spaces of
measures, including applications to EM, were established by Aubin-Frankowski, Korba
and L\'eger~\cite{aubinfrankowski2022mirror}. For Gaussian mixtures, Yan, Wang and
Rigollet~\cite{yan2024learning} proposed gradient descent in the
Wasserstein--Fisher--Rao (WFR) geometry. The Wasserstein component moves
the locations of the mixing measure, while the Fisher--Rao component
changes its mass distribution. Their analysis establishes that if the
WFR gradient descent converges weakly, then its limiting measure is an
NPMLE under the appropriate support assumptions. For computation,
\cite{yan2024learning} represents the WFR dynamics by an interacting
particle system: starting from an $m$-atomic measure, one tracks the
$m$ locations and their weights through time. However, they do not
give a prescribed finite time $T=T(\varepsilon)$ at which the likelihood is known to
be within $\varepsilon$ of $\widehat L_n$. 

The starting point of the present work is a different use of the
Fisher--Rao geometry. We consider the score
\begin{equation}
\label{eq:intro-score}
    \Psi_\mu(\theta)
    :=
    \frac1n\sum_{i=1}^n
    \frac{\phi(X_i-\theta)}{f_\mu(X_i)}
\end{equation}
and the Fisher--Rao gradient flow
\begin{equation}
\label{eq:intro-fr}
    \partial_t\mu_t
    =
    \bigl(\Psi_{\mu_t}-1\bigr)\mu_t.
\end{equation}
Instead of initializing \eqref{eq:intro-fr} by finitely many particles,
we initialize it by the uniform probability measure $\mu_0$ on a compact
set containing the convex hull of the observations. Thus $\mu_t$ is
absolutely continuous with respect to $\mu_0$ and, at every finite time,
is non-atomic. At first sight, this appears computationally less convenient
than a particle representation, since the state variable in
\eqref{eq:intro-fr} is a continuous density. A central observation of
this paper is that the entire trajectory along the flow lies in an
$n$-parameter exponential family, allowing for a finite representation of a continuous density.

Thus the use of an absolutely continuous initialization does not require
the numerical evolution of an infinite-dimensional object. The
optimization state consists of $n$ scalar coefficients.
This exact finite-dimensional representation is one of the main reasons
that the pure Fisher--Rao dynamics are computationally useful in the
present problem. Furthermore, we establish an explicit likelihood gap is $O(\log t/t)$ for the flow. 

In principle, our algorithm is an entropic proximal discretization of this flow, which retains an exact finite-dimensional and absolutely continuous representation, much like the original continuous-time flow. A major hurdle in implementing this algorithm is solving the proximal equations, which requires numerical evaluation of integrals in high dimension, and a naive tensorized quadrature rule on the entire domain can lead to an $\varepsilon^{-O(d)}$ term in the algorithmic cost. A key observation is that the integrands we need to evaluate are exponentials of potentials which, when normalized, converge to the optimal score function along the flow. Consequently, these integrands become increasingly concentrated near the maxima of the optimal score. Our algorithm seeks to restrict the integration domain to neighborhoods of these maxima after certifying that the contributions from their complement are negligible. Once this is done, we restrict the integration domains around the approximate maxima and establish uniform analytic bounds for the resulting integrands, which allows us to use only $O(\log(1/\varepsilon))$ quadrature nodes in each coordinate.

This adaptive analysis guarantees the complexity bound \eqref{eq:intro-complexity} when the optimal score function has nondegenerate maxima. We show that this condition holds for almost every dataset $X_1,\ldots,X_n \in \mathbb{R}^d$ whenever
$
n>(2m+1)((m+1)(d+1)-1),
$
where $m$ is the cardinality of the support of the NPMLE $\hat{\mu}$. Additionally, we show that $m$ is equal to the number of contact points of the optimal score function. This extends the genericity result of \cite{polyanskiy2025nonparametric} to higher dimensions, and may be of independent interest. We note that by a recent result of Wang\cite{wang2026almost}, the Gaussian location mixture NPMLE is unique for Lebesgue-almost every dataset. 

\subsection{Organization of the Paper}

We state the main results and algorithmic constructions in Section 2, with proofs given in later sections. In particular,
\begin{enumerate}
    \item In Section~\ref{sec:FR}, we state results regarding well-posedness and convergence of the continuous-time Fisher--Rao flow, along with some of their consequences and extensions. The proofs of the results in this section are given in Section~\ref{sec:FR proofs}.

    \item In Section~\ref{sec:Comp}, we describe the basic entropic proximal discretization, which is solved using a modified fixed-point iteration method. The exact scheme is stated in Algorithm~\ref{alg:fc-global}, which is a preliminary version of our accelerated algorithm. The proofs of the results in this section, including the analysis of Algorithm~\ref{alg:fc-global}, are given in Section~\ref{sec:Comp proofs}.

    \item In Section~\ref{sec:adaptive}, we state Algorithm~\ref{alg:fa-adaptive}, which makes use of the adaptive quadrature scheme described above. In particular, we describe how to select appropriate integration neighborhoods and construct a local quadrature rule to accelerate the quadrature operations. The proofs of this section, including the analysis of Algorithm~\ref{alg:fa-adaptive}, are given in Section~\ref{sec:adaptive proofs}.

    \item Finally, in Section~\ref{sec:generic}, we state our genericity result that enables accelerated computation in Algorithm~\ref{alg:fa-adaptive}. The proof, which draws on the techniques used in \cite{polyanskiy2025nonparametric}, is given in Section~\ref{sec:generic proofs}.
\end{enumerate}

\section{Main results}

\subsection{Results for Continuous-Time Fisher-Rao Gradient Flow} \label{sec:FR}
\subsubsection{Finite Time Guarantees for Fisher-Rao Gradient Flow}
Our main theoretical results in this subsection consists of a guarantee for finite-time convergence for the Fisher-Rao flow 
\begin{equation}
  \partial_t\mu_t = (\Psi_{\mu{t}}-1)\mu_t, \quad \mu_0 = \textrm{Unif}(B(X_1,\rho_n+1)),
  \label{eq:fr}
\end{equation}
where $\Psi_t$ is as in (6) and $\rho_n := \max\{\max_{1\le i,j\le n}\|X_i-X_j\|,1\}$ to an $\varepsilon-$approximate NPMLE. In particular, the convergence result holds without any assumptions on the mixing measure. For notational convenience, we define
\[
K:= B(X_1,\rho_n+1)
\]

Before stating our main theorem regarding finite time guarantees, we first state a result regarding the well-posedness of the Fisher-Rao flow. The well-posedness of Fisher--Rao gradient flows has been
studied in several settings. In particular, Liu, Majka and
Szpruch~\cite{liu2023polyak} establish global well-posedness
for the Fisher--Rao gradient flow associated with
$
V^\sigma(m)
=
F(m)+\frac{\sigma^2}{2}\operatorname{KL}(m\|\pi),
\ \ \sigma>0,
$
where $F$ is a general functional satisfying suitable
regularity assumptions and $\pi\propto e^{-U(x)}\,dx$.
Their argument, however, does not directly apply to
the unregularized flow considered here. In our setting,
we establish global well-posedness using a direct
Banach-space ODE argument.

\begin{lemma}[Well-posedness and monotonicity]\label{lem:wellposed}
Let $\mu_0\in\PP(K)$ have a strictly positive density $p_0$ on $K$. Then
\eqref{eq:fr} has a unique global solution
$(\mu_t)_{t\ge0}\subseteq\PP(K)$ with densities $(p_t)_{t\ge0}$ satisfying
\begin{equation}\label{eq:density}
\log\frac{p_t(\theta)}{p_0(\theta)}
=
\int_0^t\bigl(\Psi_{\mu_s}(\theta)-1\bigr)\,ds,
\quad \theta\in K.
\end{equation}
Moreover, $t\mapsto L_n(f_{\mu_t})$ is nondecreasing, with
\[
\frac d{dt}L_n(f_{\mu_t})
=
\operatorname{Var}_{\mu_t}(\Psi_{\mu_t}).
\]
\end{lemma}

We now state our main continuous time result:

\begin{theorem}[Continuous-time finite-horizon guarantee]\label{thm:fr}
Let $n\ge2$, let $X_1,\dots,X_n\in\R^d$ be arbitrary, and let
$\eps\in(0,1]$. Let $\mu_0$ be the uniform probability measure on
$B(X_1,\rho_n+1)$ and let $(\mu_t)_{t\ge0}$ solve \eqref{eq:fr}. Then
$t\mapsto L_n(f_{\mu_t})$ is nondecreasing, and
\[
T\ge\frac{2d}{\eps}\log\frac{10\rho_n^2}{\sqrt\eps}
\quad\Longrightarrow\quad
L_n(f_{\mu_T})\ge\Lhat-\eps.
\]
\end{theorem}

\begin{corollary}[Direct convergence rate]\label{cor:fr-direct-rate}
For every $T\ge2d/3$,
\begin{equation}\label{eq:fr-direct-rate}
\Lhat-L_n(f_{\mu_T})
\le
\frac{d}{T}
\left(
\frac12+
\log\left(
10\rho_n^2\sqrt{\frac{T}{d}}
\right)
\right).
\end{equation}
In particular, for fixed $d$ and $\rho_n$,
\[
\Lhat-L_n(f_{\mu_T})
=
O\left(\frac{\log T}{T}\right).
\]
\end{corollary}





\subsubsection{Further Consequences and Extensions}
\begin{corollary}\label{cor:fr-cluster}
Along the Fisher--Rao flow,
\[
L_n(f_{\mu_t})\uparrow\Lhat
\quad\text{as }t\to\infty.
\]
Moreover, every weak subsequential limit of $(\mu_t)_{t\ge0}$ is an NPMLE.
In particular, if the NPMLE mixing measure is unique, then
$\mu_t\Rightarrow\muhat$ as $t\to\infty$.
\end{corollary}

\textbf{Notation:}
For a mixing measure $\mu\in\mathcal P(\mathbb R^d)$, we shall refer to
\[
f_\mu(X_i)
=\int_{\mathbb R^d}\phi(X_i-\theta)\,\mu(d\theta),
\quad 1\le i\le n,
\]
as its \textit{fitted values}.

\begin{corollary}\label{cor:fr-fitted-values}
Let
\[
u_t
:=
\bigl(f_{\mu_t}(X_1),\dots,f_{\mu_t}(X_n)\bigr),
\]
and let $u^\star$ denote the unique fitted-value vector attained by an NPMLE.
Then
\begin{equation}\label{eq:fitted-value-bound}
\norm{u_t-u^\star}_2
\le
B\sqrt{2n\bigl(\Lhat-L_n(f_{\mu_t})\bigr)},
\end{equation}
where $B = \phi(0)$. Consequently, for $T\ge2d/3$,
\[
\norm{u_T-u^\star}_2
\le
B\sqrt{
\frac{2nd}{T}
\left(
\frac12+
\log\left(
10\rho_n^2\sqrt{\frac{T}{d}}
\right)
\right)
}.
\]
\end{corollary}

\begin{corollary}\label{cor:fr-strict}
For the initialization of Theorem~\ref{thm:fr}, if $d\ge1$, then
\[
\frac d{dt}L_n(f_{\mu_t})>0
\quad\text{for every finite }t\ge0.
\]
\end{corollary}

\begin{remark}[Dependence on the initialization]\label{rem:fr-init}
The uniform initialization is used only to obtain an explicit entropy bound. More
generally, suppose $\mu_0=p_0\,dx$ on $K$ with
$\inf_{\theta\in K}p_0(\theta)\ge c_0>0$. For the mollified comparator
$\pi_\sigma=\muhat*\mathrm{Unif}(B(0,\sigma))$,
Lemma~\ref{lem:mollify} gives
$\sup p_{\pi_\sigma}\le\abs{B(0,\sigma)}^{-1}$, and therefore
\[
\KL(\pi_\sigma\,\|\,\mu_0)
\le
\log\frac{1}{c_0\abs{B(0,\sigma)}}.
\]
Thus the same argument gives an explicit Fisher--Rao rate for every
initialization whose density is bounded below on $K$.
\end{remark}

\begin{proposition}[Fisher--Rao for general mixture likelihoods]
\label{prop:general-fr}
Let $\Theta$ be a compact metric space and let
$k_1,\dots,k_n:\Theta\to(0,\infty)$ be continuous functions satisfying
\[
0<b\le k_i(\theta)\le B_0<\infty
\quad
\text{for all }i\text{ and }\theta\in\Theta.
\]
For $\mu\in\PP(\Theta)$ define
\[
F_i(\mu):=\int_\Theta k_i(\theta)\,\mu(d\theta),
\quad
L(\mu):=\frac1n\sum_{i=1}^n\log F_i(\mu),
\]
and
\[
\Psi_\mu(\theta)
:=
\frac1n\sum_{i=1}^n\frac{k_i(\theta)}{F_i(\mu)}.
\]
Then, for every $\mu_0\in\PP(\Theta)$, the Fisher--Rao equation
\[
\partial_t\mu_t=(\Psi_{\mu_t}-1)\mu_t
\]
has a unique global solution in $\PP(\Theta)$ satisfying
\[
\frac{d\mu_t}{d\mu_0}(\theta)
=
\exp\left(
\int_0^t(\Psi_{\mu_s}(\theta)-1)\,ds
\right)
\quad \mu_0\text{-a.e.}
\]
Moreover,
\[
\frac d{dt}L(\mu_t)
=
\operatorname{Var}_{\mu_t}(\Psi_{\mu_t})
\ge0,
\]
and, for every $\pi\in\PP(\Theta)$ with
$\KL(\pi\,\|\,\mu_0)<\infty$,
\[
L(\pi)-L(\mu_T)
\le
\frac{\KL(\pi\,\|\,\mu_0)}{T}.
\]
\end{proposition}

\begin{corollary}\label{cor:fr-affine}
Let
\[
A:=\operatorname{aff}(X_1,\dots,X_n),
\quad
s:=\dim A.
\]
If $s=0$, the NPMLE is attained at $\delta_{X_1}$. Suppose $s\ge1$, and let
\[
K_A:=A\cap\overline B(X_1,\rho_n+1).
\]
Initialize the Fisher--Rao flow with the uniform probability measure on the
$s$-dimensional ball $K_A$, where volume is taken with respect to Lebesgue measure
on $A$. Then
\[
T
\ge
\frac{2s}{\eps}
\log\frac{10\rho_n^2}{\sqrt\eps}
\quad\Longrightarrow\quad
L_n(f_{\mu_T})\ge\Lhat-\eps.
\]
\end{corollary}

\subsection{A Computational Algorithm for Fisher Rao Gradient Flow}\label{sec:Comp}
We now provide a numerical implementation of the Fisher-Rao dynamics via an entropic proximal point method. Related
KL-proximal methods for optimization over probability measures
have been studied by Yao, Huang and
Yang~\cite{yao2024minimizing}, including applications to
nonparametric maximum likelihood estimation.

Our key observation is that the continuous Fisher--Rao trajectory, as well as its entropic proximal discretization, remains in an $n$-dimensional
exponential family. This gives us a finite-dimensional representation of the density along the flow, and provides an alternative to grid-based optimization schemes. 
The integrals appearing in this system are approximated by positive
Gauss--Legendre quadrature.

Throughout Sections~2.2--2.3, let $B=\phi(0)$ and set
\begin{equation}\label{eq:fc-domain}
 R=1+\max_{1\le i\le n}\|X_i-X_1\|,
 \quad K=X_1+[-R,R]^d,
 \quad \mu_0=\operatorname{Unif}(K).
\end{equation}
Every point in $\operatorname{conv}\{X_1,\ldots,X_n\}$ has a unit
Euclidean ball contained in $K$, and maximization of the likelihood over
$\mathcal P(K)$ is equivalent to computing the unrestricted NPMLE over $\mathcal P(\mathbb{R}^d)$.
\\
Next, for $0<\varepsilon\le1$, letting  $v_d$ denote the volume of the unit
Euclidean ball and $b: = B/4n$, we set
\begin{equation}
 \begin{aligned}
 &\label{eq:fc-tolerances}
 H_\varepsilon =\log\frac{(2R)^d}{v_d}+d\log\frac4{\sqrt\varepsilon},
 \\
 &\eta_\varepsilon =\frac{\varepsilon}{256n(1+H_\varepsilon)},
 \\
 &\delta_\varepsilon=\frac{\eta_\varepsilon b^2}{512B}.
 \end{aligned}
\end{equation}

\subsubsection{Finite-dimensional representation and proximal discretization}

For $a\in\mathbb R^n$, define
\begin{equation}\label{eq:fc-family}
 \begin{aligned}
 V_a(\theta)&=\sum_{i=1}^n a_i\phi(X_i-\theta),
 & Z(a)&=\int_K e^{V_a(\theta)}\,\mu_0(d\theta),\\
 \mu_a(d\theta)&=\frac{e^{V_a(\theta)}}{Z(a)}\,\mu_0(d\theta),
 & F_i(a)&=\int_K\phi(X_i-\theta)\,\mu_a(d\theta).
 \end{aligned}
\end{equation}
The initialization is represented by $a=0$, and
$L_n(f_{\mu_a})=n^{-1}\sum_i\log F_i(a)$. We refer to the quantities \(F_i(a)\) as the fitted values.
\begin{proposition}[Representation of the flow and its proximal steps]
\label{prop:fc-representation}
The Fisher--Rao flow initialized at $\mu_0$ satisfies
\begin{equation}\label{eq:fc-flow}
 \mu_t=\mu_{a(t)},\quad
 \dot a_i(t)=\frac1{nF_i(a(t))},\quad a_i(0)=0,
 \quad Z(a(t))=e^t.
\end{equation}
Moreover, for $a\in\mathbb R_+^n$ and $\tau>0$, the problem
\begin{equation}\label{eq:fc-prox-original}
 \max_{\mu\in\mathcal P(K)}
 \left\{L_n(f_\mu)-\frac1\tau\operatorname{KL}(\mu\|\mu_a)\right\}
\end{equation}
has a unique solution $\mu_{a^+}$, whose coefficients can be chosen to
satisfy
\begin{equation}\label{eq:fc-original-update}
 a_i^+=a_i+\frac{\tau}{nF_i(a^+)},\quad 1\le i\le n.
\end{equation}
\end{proposition}

A direct fixed-point iteration for \eqref{eq:fc-original-update} may
require a choice of the step size $\tau$ which may be small when $F_i(a)$ is close to zero.
To obtain step sizes that do not become small even when $F_i(a)$ is relatively small, we consider a modified version of the log-likelihood, where we replace logarithm below
$b=B/(4n)$ is by a linear replacement:
\begin{equation}\label{eq:fc-continuation}
\mathcal L_b(\mu)=\frac1n\sum_{i=1}^n h_b(f_\mu(X_i)),
\quad 
\text{where }
 h_b(u)=
 \begin{cases}
  \log u,&u\ge b,\\
  \log b+u/b-1,&0\le u<b.
 \end{cases}
\end{equation}
For this problem, the optimal fitted vector $(f_1^\star,\ldots,f_n^\star)$ is unique,
and the first-order optimality condition gives
$f_i^\star\ge B/n=4b$. As we shall show in
Lemma~\ref{lem:fp-equivalence}, this modification preserves the
maximizers:
\[
\operatorname*{arg\,max}_{\mu\in\mathcal P(K)}\mathcal L_b(\mu)
=
\operatorname*{arg\,max}_{\mu\in\mathcal P(K)}L_n(f_\mu),
\atop
\max_{\mu\in\mathcal P(K)}\mathcal L_b(\mu)=\widehat L_n.
\]
The corresponding proximal update for $\mathcal L_b$ is
\begin{equation}\label{eq:fc-capped-update}
 a_i^+=a_i+\frac{\tau}{n\max\{F_i(a^+),b\}},
 \qquad 1\le i\le n,
\end{equation}
which we solve using a scaled fixed-point iteration.

Although the two objectives have the same maximizers, their
values may differ away from optimality. To account for this, we shall perform
an output correction which, by
Lemma~\ref{lem:fp-equivalence}, converts an approximate maximizer
of $\mathcal L_b$ into an approximate maximizer of $L_n$, with
at most a universal constant factor in the likelihood gap.

\subsubsection{The fixed-point iteration and numerical quadrature}

We now describe how to solve the proximal system
\eqref{eq:fc-capped-update} using a scaled fixed-point iteration.
Since the fitted values $F_i(a)$ cannot generally be evaluated
exactly, we replace them by numerical approximations
$\widehat F_i(a)$, whose construction is given later.

At the $k$-th proximal step, let $p_i^k = \max\{\widehat F_i(a^k),b\} > 0,$
$ 1\le i\le n$ and $\tau_k>0$
be a step size whose exact value is specified in \eqref{eq:fc-step}. To solve
\eqref{eq:fc-capped-update}, we initialize $y^0=\mathbf1$
and, for each $r\ge0$, construct a candidate coefficient vector
\begin{equation}\label{eq:fc-inner}
x_i^r=a_i^k+\frac{\tau_k y_i^r}{np_i^k},
\qquad 1\le i\le n.
\end{equation}
We then evaluate $\widehat F(x^r)$. If
\begin{equation}\label{eq:fc-stop}
\frac Bn\sum_{i=1}^n
\left|\frac{y_i^r}{p_i^k}
-\frac1{\max\{\widehat F_i(x^r),b\}}\right|
\le\frac{\eta_\varepsilon}{2},
\end{equation}
we accept $a^{k+1}=x^r$. Otherwise, we compute the next
candidate using
\begin{equation}\label{eq:fc-picard}
y_i^{r+1}=
\frac{p_i^k}{\max\{\widehat F_i(x^r),b\}},
\qquad 1\le i\le n,
\end{equation}
and repeat the procedure. 

We next choose $\tau_k$ to ensure convergence of this adapted
fixed-point iteration. First, we note that for $y^r\in[1/2,2]^n$, the coefficient vector $x^r$
defined above satisfies
\[
0\le V_{x^r}(\theta)-V_{a^k}(\theta)
\le 2\tau_k\frac1n\sum_{i=1}^n
\frac{\phi(X_i-\theta)}{p_i^k}.
\]
We choose $\tau_k$ so that this difference is bounded uniformly over $\theta\in K$. To obtain a
computable upper bound for the Gaussian sum on the right-hand side, we consider
the grid
\begin{equation}\label{eq:fc-score-grid}
\mathcal G=\left\{
X_1-R\mathbf1+\frac{2R}{\lceil2R\sqrt d\rceil}u:
u\in\{0,\ldots,\lceil2R\sqrt d\rceil\}^d
\right\}.
\end{equation}
We shall show in Lemma \ref{lem:fp-gaussian} that the supremum of the sum over $K$ is bounded
by $e^{1/8}$ times its maximum over $\mathcal G$.
We therefore set
\begin{equation}\label{eq:fc-step}
\tau_k=\frac1{64}
\left[\max\left\{1,
e^{1/8}\max_{\theta\in\mathcal G}
\frac1n\sum_{i=1}^n
\frac{\phi(X_i-\theta)}{p_i^k}
+\frac{B\delta_\varepsilon}{b^2}
\right\}\right]^{-1}.
\end{equation}
The additional term involving $\delta_\varepsilon$
accounts for the error in the approximations
$\widehat F_i$. As shown in Lemma~\ref{lem:fp-inner},
this choice makes the fixed-point iteration contractive
provided \eqref{eq:fc-evaluation-error} holds.

It remains to construct the approximations $\widehat F_i(a)$
used in the preceding iteration. Recall that
\[
F_i(a)=
\frac{\int_K\phi(X_i-\theta)e^{V_a(\theta)}\,d\theta}
{\int_K e^{V_a(\theta)}\,d\theta}.
\]
We approximate these integrals using tensor-product
Gauss--Legendre quadrature.

Let $\{(z_r,\omega_r):1\le r\le q\}$ denote the
Gauss--Legendre rule on $[-1,1]$, with $\omega_r>0$
and $\sum_r\omega_r=2$. For
$u\in\{1,\ldots,q\}^d$, the corresponding tensor-product
rule on $K$ has nodes and normalized weights
\[
\theta_u=X_1+R(z_{u_1},\ldots,z_{u_d}),
\qquad
w_u=2^{-d}\prod_{j=1}^d\omega_{u_j}.
\]
Applying this rule to the numerator and denominator
of $F_i(a)$ gives
\begin{equation}\label{eq:fc-quadrature}
\widehat F_i(a)=
\frac{\sum_u w_u\phi(X_i-\theta_u)
e^{V_a(\theta_u)-\max_v V_a(\theta_v)}}
{\sum_u w_u e^{V_a(\theta_u)-\max_v V_a(\theta_v)}}.
\end{equation}
Note that the positivity of the weights gives
$0<\widehat F_i(a)\le B$.

We choose the quadrature order according to the accumulated
proximal time $t_k=\sum_{j<k}\tau_j$. Define
\begin{equation}\label{eq:fc-potential-bound}
\mathcal M(t)=3t+2\log\frac{(2R)^d}{v_d}
+d\log(d+1)+2.
\end{equation}
For a proximal step starting at time $t$, let
$\overline M$ be the smallest power of two satisfying
$\overline M\ge\mathcal M(t+1)$, and choose
\begin{equation}\label{eq:fc-global-order}
q_K=\left\lceil
16d^{3/2}R\sqrt{1+\overline M}
\log\frac{16dR(1+\overline M)B}{\delta_\varepsilon}
\right\rceil.
\end{equation}
With this choice, the quadrature approximation satisfies
\begin{equation}\label{eq:fc-evaluation-error}
\max_i|\widehat F_i(a)-F_i(a)|\le\delta_\varepsilon
\end{equation}
at the current iterate $a^k$ and at every candidate
coefficient vector $x^r$ generated during the fixed-point
iteration.

Note that $\overline M$ takes only dyadic values, and the
quadrature rule is reconstructed only when $\overline M$
increases. The precise construction of the nodes and weights, together with the required accuracy estimates, is given in Section 4.4.

\subsubsection{Algorithm and convergence guarantees}

\begin{algorithm}[H]
\caption{Fisher--Rao proximal iteration with quadrature on $K$}
\label{alg:fc-global}
\begin{algorithmic}[1]
\Require $X_1,\ldots,X_n\in\mathbb R^d$ and $0<\varepsilon\le1$.
\State Form $K$, $\mathcal G$, and the parameters in
       \eqref{eq:fc-tolerances}; set $a^0=0$, $t_0=0$, and $k=0$.
\While{$t_k<32H_\varepsilon/\varepsilon$}
 \State Construct or reuse the rule \eqref{eq:fc-global-order};
        evaluate $\widehat F(a^k)$ and compute $p^k,\tau_k$ by
        \eqref{eq:fc-step}.
 \State Set $y^0=\mathbf1$; apply
        \eqref{eq:fc-inner}--\eqref{eq:fc-picard} until
        \eqref{eq:fc-stop} holds, and set $a^{k+1}=x^r$.
 \State Set $t_{k+1}=t_k+\tau_k$ and record
        $n^{-1}\sum_i h_b(\widehat F_i(a^{k+1}))$,
        using the fitted values from the accepted trial.
 \State Retain the coefficients and fitted values with the largest
        recorded value; set $k\gets k+1$.
\EndWhile
\State Let $a^{\widehat k}$ be the retained coefficients and set
 $I=\{i:\widehat F_i(a^{\widehat k})<b+\delta_\varepsilon\}$.
\State \Return
 $\displaystyle\mu_{\rm out}=
 (1-|I|/n)\mu_{a^{\widehat k}}
       +\frac1n\sum_{i\in I}\operatorname{Unif}(B(X_i,1))$.
\end{algorithmic}
\end{algorithm}

The algorithm retains the iterate with the largest computed value of
$\mathcal L_b$. If some of its fitted values are below $b$, a mass
$1/n$ is added uniformly on the corresponding balls $B(X_i,1)$;
the remaining mass is assigned to the retained exponential-family
measure.

\begin{theorem}[Accuracy and iteration bounds]
\label{thm:fc-iterations}
For arbitrary $X_1,\ldots,X_n\in\mathbb R^d$ and
$0<\varepsilon\le1$, Algorithm~\ref{alg:fc-global}
terminates after $N$ iterations and returns an absolutely
continuous probability measure $\mu_{\rm out}$ satisfying
\begin{equation}\label{eq:fc-accuracy}
L_n(f_{\mu_{\rm out}})\ge\widehat L_n-\varepsilon.
\end{equation}
For every $0\le k<N$, the fixed-point iteration used to
compute $a^{k+1}$ requires evaluating
$(\widehat F_i(a))_{i=1}^n$ at no more than
\[
2+\left\lceil
\frac{\log(8B/(b\eta_\varepsilon))}{\log 8}
\right\rceil
\]
coefficient vectors $a\in\mathbb R^n$.

Moreover, for every $0\le k\le N$,
\begin{equation}\label{eq:fc-counts}
\begin{aligned}
k&\le C\left[
t_k+n\left(
\log\frac{(2R)^d}{v_d}
+\frac12+\frac d2\log(1+t_k)
\right)\right],\\
N&\le C_dH_\varepsilon
\min\left\{n+\varepsilon^{-1},R^d\varepsilon^{-1}\right\},
\end{aligned}
\end{equation}
where $C$ is a numerical constant.

Furthermore, there exists an integer
$k_\star\le C_dn\log(enR)$, independent of $\varepsilon$,
such that $\tau_k>1/128$ for every $k_\star\le k<N$.
Consequently,
\[
N\le k_\star+\frac{4096H_\varepsilon}{\varepsilon}+2.
\]
Here $C_d$ depends only on $d$.
\end{theorem}

\begin{theorem}[Arithmetic complexity]
\label{thm:fc-global-cost}
For arbitrary $X_1,\ldots,X_n\in\mathbb R^d$ and
$0<\varepsilon\le1$, the total arithmetic complexity of
Algorithm~\ref{alg:fc-global} is
\begin{equation}\label{eq:fc-global-cost}
 O_d\left(
nR^{2d}\varepsilon^{-1-d/2}L^{3d/2 + 2}
+R^2\varepsilon^{-1}L^{4}
\right),
\end{equation}
where $L := (1 + \log(nR/\varepsilon)).$
\end{theorem}

\subsection{Adaptive Quadrature and Accelerated Convergence}\label{sec:adaptive}
In Section~2.2, we perform the quadrature over the whole hypercube $K$ for all the iterations.
However, the integrands considered in the definition of $F_i(a)$ become increasingly concentrated near the
maxima of $V_a$ as the proximal iteration proceeds. We can exploit this
concentration to speed up our algorithm, by restricting the integration domain to a union of
cubes whose centers and widths are adapted to the current coefficients.

For $t>0$, consider $W_{a,t}=V_a/t$. If $z$ is a nondegenerate maximum
of $W_{a,t}$, then
\[
 V_a(z+t^{-1/2}u)-V_a(z)
   =\frac12u^{\mathsf T}\nabla^2W_{a,t}(z)u+o(\|u\|^2)
 \quad (u\to0).
\]
This suggests that we should integrate over cubes of width proportional to $t^{-1/2}$ around
the local maxima. Based on this idea, we consider a  modification of
Algorithm~\ref{alg:fc-global}. The overarching idea is outlined below:
\begin{enumerate}[label=(\roman*),leftmargin=2em]
 \item Proximal steps are initially computed using quadrature on $K$.
 \item At prescribed times, $V_a/t$, $\nabla V_a/t$, and
       $\nabla^2V_a/t$ are evaluated on a uniform grid over $K$
       to select disjoint cubes containing its \textit{dominant} maxima.
 \item Once these cubes have been selected, within them, we construct smaller
      cubes around the approximate maximizers, with widths
      proportional to $t^{-1/2}$. When the mass outside these
      smaller cubes is at most $\delta_\varepsilon/(8B)$,
      we restrict quadrature to them for subsequent evaluations.
\end{enumerate}
The proximal equations and the stopping rule remain unchanged.

\subsubsection{Selection of the cubes}

The procedure corresponds to Step 2 in the above scheme. We attempt to select at the first completed proximal step for which
$t_k\ge2^j$, successively for $j=1,2,\ldots$, until a selection succeeds. At such a time, denote $a=a^k$, $t=t_k$, and $W=V_a/t$.
For centers $c_1,\ldots,c_m$ and a common radius $r$, define
\begin{equation}\label{eq:fa-cubes}
 Q_\ell=c_\ell+[-r,r]^d,\quad
 U_\ell=c_\ell+[-2r,2r]^d,\quad
 \mathcal U=\bigcup_{\ell=1}^m U_\ell.
\end{equation}
The larger cubes $U_\ell$ are used to bound the Hessian and the
discarded integrals; the maxima are computed within $Q_\ell$.

Let $f_1^\star,\ldots,f_n^\star$ be the unique optimal fitted values
and define
\begin{equation}\label{eq:fa-optimal-score}
 \Psi^\star(\theta)=\frac1n\sum_{i=1}^n
                 \frac{\phi(X_i-\theta)}{f_i^\star},
 \quad \mathcal C^\star=\{\theta:\Psi^\star(\theta)=1\}.
\end{equation}

\paragraph{Procedure \textnormal{\textsc{SelectCubes}}$(a,t,j)$.}
The following operations are performed in the given order.
\begin{enumerate}[label=\arabic*.,leftmargin=2em,itemsep=4pt]
\item Consider $h=1,\ldots,\max\{1,\lfloor j/4\rfloor\}$ in
      increasing order. For each $h$, perform Steps~2--4 with
      $s=2R\,2^{-h}$, $r=\min\{1/8,s^{1/4}\}$, and 
      \[
       C_u=c_u+[-s/2,s/2]^d,\ \ 
       c_u=X_1-R\mathbf1+s(u+\tfrac12\mathbf1),
      \]
      for $u\in\{0,\ldots,2^h-1\}^d$. At every midpoint $c_u$, evaluate $W(c_u)$, $\nabla W(c_u)$, and $\nabla^2W(c_u)$.
     
\item Sort the midpoints satisfying
      \[
       W(c_u)\ge\max_v W(c_v)-2\mathcal M(1)\sqrt d\,s
      \]
      by decreasing value of $W$, and break ties in a lexicographic order.
      Retain the first midpoint, and subsequently, retain a midpoint if its
      distance from every previously retained midpoint
      exceeds $5r$ in the infinity-norm. Denote the retained points by $c_1,\ldots,c_m$
      and form \eqref{eq:fa-cubes}.
\item Reject this partition if
      $\|c_\ell-X_1\|_\infty+3r\ge R$ for some $\ell$.
      Otherwise, let
      \[
       \mathcal O=\{u:C_u\cap
           (K\setminus\bigcup_\ell\operatorname{int}U_\ell)
                         \ne\varnothing\}
      \]
      and compute
      \begin{equation}\label{eq:fa-gap}
       \gamma=\frac14\left[
       \max_\ell W(c_\ell)-
       \max_{u\in\mathcal O}
                 \{W(c_u)+\mathcal M(1)\sqrt d\,s\}\right].
      \end{equation}
      Reject the partition if $\gamma\le0$.
\item For $p=0,\ldots,j+1$, set $\kappa=2^{-p}$ and check whether
      \begin{equation}\label{eq:fa-accept}
      \begin{aligned}
       -\nabla^2W(c_u)&\succ
          (4\kappa+24\mathcal M(1)\sqrt d\,s)I_d
                     \ \ \text{if }C_u\cap\mathcal U\ne\varnothing,\\
       \|\nabla W(c_\ell)\|&<\kappa r/8
                    \ \ \text{for } 1\le\ell\le m ,\\
       E_\varepsilon(t)&\le\min\{\gamma/4,\kappa r/16\}.&&
      \end{aligned}
      \end{equation}
      Accept the first pair $(h,p)$ for which all inequalities hold.
      Retain $t_{\rm sel}=t$, $c_1,\ldots,c_m$, $r$, $\kappa$,
      and $\gamma$. Return $\varnothing$ only after all listed
      pairs have been examined without acceptance.
\end{enumerate}
We note that the strict matrix
inequalities can be checked by Cholesky factorization. On acceptance, set $\Lambda=32\mathcal M(1)$.

\subsubsection{Local quadrature}

After selecting the cubes, we retain their centers $c_\ell$ and their
radius $r$. We want to replace global quadrature on $K$ with local quadrature on smaller cubes around the maxima. Hence, at each evaluation of \(F(a)\), we first approximate the maximizer of \(W=V_a/t\) on each \(Q_\ell\) using projected gradient ascent. We then place the quadrature nodes on an appropriate neighbourhood around these approximate maximizers. More specifically, starting at
$z_\ell^0=c_\ell$, we perform
\begin{equation}\label{eq:fa-ascent}
 z_\ell^{v+1}=\Pi_{Q_\ell}\left(
        z_\ell^v+\Lambda^{-1}\nabla W(z_\ell^v)\right),\quad
 0\le v<J(t),
 \quad
 J(t)=\left\lceil\frac\Lambda\kappa
             \max\{0,\log(r\sqrt{\Lambda t})\}\right\rceil.
\end{equation}
Next, let $z_\ell=z_\ell^{J(t)}$ and 
\begin{equation}\label{eq:fa-width}
 L^2=\frac8\kappa\max\left\{1,
 \log\frac{16m e^{17/32}(4\pi\Lambda/\kappa)^{d/2}B}
                    {v_d\delta_\varepsilon}\right\},
\end{equation}
where we recall that $v_d$ is the volume of the unit-ball in dimension $d$.
Then, the cubes used for local quadrature are
\[
z_\ell+t^{-1/2}[-L,L]^d.
\]
We restrict quadrature to these cubes at
the first tested time satisfying
\begin{equation}\label{eq:fa-ready}
 \begin{gathered}
 t\ge\max\{1,1/(\Lambda r^2),d/(2\gamma)\},
 \\
\sqrt d\,L/\sqrt t\le r/2,\\
 t\gamma-\frac d2\log t
 \ge\log\frac{16e^{17/32}(2R)^d\Lambda^{d/2}B}
                   {v_d\delta_\varepsilon}.
 \end{gathered}
\end{equation}

Next, we construct a tensor Gauss--Legendre rule on $[-L,L]^d$,
with per-coordinate order given by
\begin{equation}\label{eq:fa-order}
 q_{\rm loc}=\max\left\{1,
 \left\lceil\frac{1}{2\log2}
 \left[4d\Lambda L^2+
 \log\frac{2^{3d+7}dm e^{1/2}(\Lambda L^2)^{d/2}B}
                         {v_d\delta_\varepsilon}\right]
 \right\rceil\right\}.
\end{equation}
\\
We denote the nodes and the weights of this rule by $(y_u,\omega_u)$, and set
$\theta_{\ell u}=z_\ell+t^{-1/2}y_u$. We evaluate
\begin{equation}\label{eq:fa-quadrature}
 \widehat F_i(a)=
 \frac{\displaystyle\sum_{\ell,u}\omega_u\phi(X_i-\theta_{\ell u})
                  e^{V_a(\theta_{\ell u})-\max_jV_a(z_j)}}
      {\displaystyle\sum_{\ell,u}\omega_u
                  e^{V_a(\theta_{\ell u})-\max_jV_a(z_j)}}.
\end{equation}
Note that the Jacobian $t^{-d/2}$ is common to all cubes and cancels from this
ratio. 


\begin{algorithm}[H]
\caption{Fisher--Rao proximal iteration with adaptive quadrature}
\label{alg:fa-adaptive}
\begin{algorithmic}[1]
\Require $X_1,\ldots,X_n\in\mathbb R^d$ and $0<\varepsilon\le1$.
\State Form $K$, $\mathcal G$, and the parameters in
       \eqref{eq:fc-tolerances}; set $a^0=0$, $t_0=0$, $k=0$, and $j=1$.
\State Use quadrature on $K$ and record that no cubes have been accepted.
\While{$t_k<32H_\varepsilon/\varepsilon$}
 \State If quadrature is on $K$, construct or reuse the rule
        \eqref{eq:fc-global-order} for the current time $t_k$.
 \State Evaluate $\widehat F(a^k)$ and compute $p^k,\tau_k$ by
        \eqref{eq:fc-step}, using the current quadrature rule.
 \State Set $y^0=\mathbf1$; apply
        \eqref{eq:fc-inner}--\eqref{eq:fc-picard} with the current quadrature rule
        until \eqref{eq:fc-stop} holds, and set $a^{k+1}=x^r$.
 \State Set $t_{k+1}=t_k+\tau_k$. Retain $a^{k+1}$ and its fitted values if
$
n^{-1}\sum_{i=1}^n h_b(\widehat F_i(a^{k+1}))
\geq
\max_{1\leq j\leq k}
n^{-1}\sum_{i=1}^n h_b(\widehat F_i(a^j)),
$
with the condition taken to be true when $k=0$.
 \If{$2^j\le t_{k+1}<32H_\varepsilon/\varepsilon$ and quadrature is still on $K$}
  \If{no cubes have been accepted}
   \State Apply \textnormal{\textsc{SelectCubes}}$(a^{k+1},t_{k+1},j)$;
          if cubes are accepted, retain them and compute $L$ by
          \eqref{eq:fa-width}.
  \EndIf
  \If{cubes have been accepted and \eqref{eq:fa-ready} holds at $t_{k+1}$}
   \State Construct the rule \eqref{eq:fa-order} on $[-L,L]^d$.
   \State For every subsequent evaluation, compute the centers by
          \eqref{eq:fa-ascent} and use \eqref{eq:fa-quadrature}.
  \EndIf
  \State Set $j\gets j+1$.
 \EndIf
 \State Set $k\gets k+1$.
\EndWhile
\State Let $a^{\widehat k}$ be the retained coefficients and set
 $I=\{i:\widehat F_i(a^{\widehat k})<b+\delta_\varepsilon\}$.
\State \Return
 $\displaystyle(1-|I|/n)\mu_{a^{\widehat k}}
       +\frac1n\sum_{i\in I}\operatorname{Unif}(B(X_i,1))$.
\end{algorithmic}
\end{algorithm}

For each inner trial in the $k$-th proximal step, we use
$t=t_k$ in \eqref{eq:fa-ascent} and \eqref{eq:fa-quadrature}.
When $a^{k+1}$ is subsequently evaluated as a completed iterate,
we use $t=t_{k+1}$. The same convention applies to the
potential $W_{a,t}$ used in the analysis below.

\subsubsection{Accuracy and arithmetic complexity}

We now state the accuracy and arithmetic complexity guarantees for Algorithm \ref{alg:fa-adaptive}. We first note that the first-order condition for the NPMLE gives that $\Psi^\star\le 1$ on $\mathbb R^d$. Moreover, since $\Psi^\star$ is a finite positive sum of isotropic Gaussian kernels, Theorem 1 of \cite{wang2026finite} implies that $\mathcal C^\star$ is finite. The only additional condition required for the accelerated complexity bound is
\begin{equation}
-\nabla^2\Psi^\star(\xi)\succ 0
\qquad\text{for every }\xi\in\mathcal C^\star.
\label{eq:fa-morse}
\end{equation}
In Theorem \ref{thm:generic-nondegeneracy} in the subsequent section, we show that this condition holds for almost every dataset $X_1,\cdots,X_n \in \mathbb{R}^d$, given that $n > 2(d+1)(m+1)^2$, where $m = |\text{supp}(\hat{\mu})|$

\begin{theorem}[Accuracy and eventual use of local quadrature]
\label{thm:fa-validity}
Let $X_1,\ldots,X_n\in\mathbb R^d$ and $0<\varepsilon\le1$.
Then the following statements hold.

\begin{enumerate}
\item[(i)] Algorithm~\ref{alg:fa-adaptive} satisfies the accuracy
and iteration bounds of Theorem~\ref{thm:fc-iterations}.

\item[(ii)] Suppose that the procedure \textnormal{\textsc{SelectCubes}} succeeds.
Let $Q_1,\ldots,Q_m$ and $U_1,\ldots,U_m$ denote the selected
cubes. For every subsequent iterate or inner trial, with
coefficient vector $a$ and corresponding time $t$, we have
\begin{equation}\label{eq:fa-persistent-geometry}
\begin{aligned}
\kappa I_d\preceq-\nabla^2W_{a,t}(\theta)
\preceq\Lambda I_d
\quad \forall \theta\in\mathcal U,
\atop
\sup_{K\setminus\mathcal U}W_{a,t}
\le\max_\ell W_{a,t}(c_\ell)-\gamma.
\end{aligned}
\end{equation}
Moreover, $W_{a,t}$ has a unique maximizer $z_\ell^\star$
on each $U_\ell$, satisfying $\|z_\ell^\star-c_\ell\|\le r/4.$
Whenever \eqref{eq:fa-ready} holds, the local quadrature rule
\eqref{eq:fa-quadrature} satisfies
\eqref{eq:fc-evaluation-error}, with
\begin{equation}\label{eq:fa-local-cost}
\begin{aligned}
q_{\rm loc}
\le C_d\frac{\Lambda}{\kappa}
\log\frac{em(1+\Lambda/\kappa)B}{\delta_\varepsilon},
\atop
\text{cost per evaluation}
=O\bigl(ndm\{q_{\rm loc}^d+J(t)\}\bigr).
\end{aligned}
\end{equation}
The cost includes all projected gradient iterations in
\eqref{eq:fa-ascent}.

\item[(iii)] Suppose additionally that \eqref{eq:fa-morse}
holds. Then, there exist $\varepsilon_0>0$ and
$T_0<\infty$, depending on the data but independent of $\varepsilon$, such that
the procedure \textnormal{\textsc{SelectCubes}} succeeds by time $T_0$
for every $0<\varepsilon\le\varepsilon_0$.

Furthermore, the number of selected cubes and the quantities
$
r^{-1},\kappa^{-1},\gamma^{-1},\Lambda
$
are bounded uniformly over these accuracies, and local
quadrature begins at time
\[
O_{X,d}\bigl(\log(e/\varepsilon)\bigr).
\]
\end{enumerate}
\end{theorem}

\begin{theorem}[Total arithmetic cost]
\label{thm:fa-total-cost}
Let $X_1,\ldots,X_n\in\mathbb R^d$ and $0<\varepsilon\le1$.
Then the following statements hold.

\begin{enumerate}
\item[(i)] Let $N$ be the total number of proximal steps,
and let $G$ be the number of steps computed using quadrature
on $K$. Setting $t_{\rm loc}=t_G$,
the total arithmetic cost of Algorithm~\ref{alg:fa-adaptive}, suppressing logarithmic factors, is
is
\begin{equation}\label{eq:fa-total-cost}
\begin{split}
 \widetilde{O} \biggl(&
nR^dN+nR^dG(1+\log R+t_{\rm loc})^{d/2}\\
&+R^2(1+\log R+t_{\rm loc})
+nm(N-G)(\Lambda/\kappa)^d
+(\Lambda/\kappa)^2
\biggr),
\end{split}
\end{equation}
where the last two terms are omitted if $G=N$. The iteration counts $N$ and $G$ satisfy \eqref{eq:fc-counts},
with $t=t_N$ and $t=t_G$, respectively. Moreover,
\[
N\le C_dR^dH_\varepsilon/\varepsilon,
\qquad
G\le C_dR^d(1+\log R+t_{\rm loc}).
\]

\item[(ii)] Suppose that \eqref{eq:fa-morse}
holds. For fixed data and sufficiently small $\varepsilon$,
the total arithmetic cost is of the order
\begin{equation}\label{eq:fa-fixed-data}
O\left(
n\log^{2+3d/2}(e/\varepsilon)
+nm_\varepsilon\varepsilon^{-1}
\log^{d+2}(e/\varepsilon)
\right),
\end{equation}
where $m_\varepsilon$ is the number of retained cubes.
\end{enumerate}
\end{theorem}

We remark that $m_\varepsilon$ remains uniformly bounded as $\varepsilon \rightarrow0$. Moreover, in practice, we expect $m_\varepsilon$ to remain close to $m$, the number of global maxima of $\Psi^*$.


\subsection{Generic Non-degeneracy of Optimal Score}\label{sec:generic}
Let $f_1^\star,\ldots,f_n^\star$ denote the unique
fitted values of an NPMLE, and define
\begin{equation}\label{eq:generic-score}
\Psi^\star(\theta)
:=
\frac1n\sum_{i=1}^n
\frac{\phi(X_i-\theta)}{f_i^\star},
\qquad
\mathcal C^\star
:=
\{\theta\in\mathbb R^d:\Psi^\star(\theta)=1\}.
\end{equation}
The first-order optimality condition gives
\[
\Psi^\star(\theta)\le1,
\qquad \theta\in\mathbb R^d,
\]
with equality at every point in the support of an NPMLE.
We show that, for almost every dataset, these are the
only points of equality and each is a nondegenerate
maximum, provided the sample size is sufficiently large
relative to the number of support points.

\begin{theorem}[Generic nondegeneracy]
\label{thm:generic-nondegeneracy}
Let $d,k\ge1$ and suppose that
\begin{equation}\label{eq:generic-sample-size}
n>(2k+1)\bigl((k+1)(d+1)-1\bigr).
\end{equation}
There exists a Lebesgue-null set
$\mathcal N_{n,d,k}\subset(\mathbb R^d)^n$ such that,
for every $(X_1,\ldots,X_n)\notin\mathcal N_{n,d,k}$
and every $k$-atomic global NPMLE $\widehat\mu$,
\begin{equation}\label{eq:generic-conclusion}
\mathcal C^\star=\operatorname{supp}(\widehat\mu),
\qquad
-\nabla^2\Psi^\star(\xi)\succ0,
\quad \xi\in\mathcal C^\star.
\end{equation}
In particular, the conclusion holds almost surely
for observations with a jointly absolutely continuous
distribution.
\end{theorem}

The condition \eqref{eq:generic-sample-size} is satisfied
whenever
\[
n>2(d+1)(k+1)^2.
\]
Thus, under the hypotheses of
Theorem~\ref{thm:generic-nondegeneracy}, the contact
set consists of exactly $k$ nondegenerate maxima.
In particular, the nondegeneracy condition used in
the adaptive quadrature algorithm of the preceding section
is satisfied.

\section*{AI Statement}
In this work, we have used ChatGPT 6 Pro for exploratory simulations, literature review, calculations, and assistance with portions of the proofs. The authors developed the principal research directions, algorithmic constructions, and proof strategies, while AI assisted in elaborating certain technical arguments, particularly those involving detailed calculations and the derivation of explicit constants required for the algorithms.

The central ideas of using the Fisher–Rao flow, its exact finite-dimensional exponential-family representation, the entropic proximal discretization, and localized quadrature near the maxima of the optimal score to accelerate computation were proposed by the authors. The modification of the log-likelihood objective introduced in Section~\ref{sec:Comp} was suggested by AI after the authors observed that controlling the Lipschitz constant of the fixed-point map for the original log-likelihood required very small step sizes, substantially slowing the algorithm.

The authors have verified all mathematical arguments and take full responsibility for the final content of this work.


\bibliographystyle{plain}
\phantomsection
\addcontentsline{toc}{section}{References}
\bibliography{bibliography}

\section{Proofs of Results in Section \ref{sec:FR}}\label{sec:FR proofs}
Throughout this section, we fix the data $\{X_1,\dots,X_n\}$ and let
$K:=\overline B(X_1,\rho_n+1)$. We let
$B:=(2\pi)^{-d/2}$ and
$\beta_1:=Be^{-2(\rho_n+1)^2}$. 
\\
Note that for any $\theta\in K$ we have
$\norm{X_i-\theta}\le2(\rho_n+1)$. Hence, for every
$\mu\in\PP(K)$,
$f_\mu(X_i)\ge\beta_1$ and consequently
$0<\Psi_\mu\le B/\beta_1$ on $K$. 
\\
Finally, we let
\[
D_n:=\max_{1\le i,j\le n}\norm{X_i-X_j},
\quad
\rho_n:=\max(1,D_n).
\]
\subsection{Well-posedness of the Fisher--Rao flow}
\begin{proof}[\textbf{Proof of Lemma \ref{lem:wellposed}}]
We work on the Banach space $\mathcal M(K)$ of finite signed measures on $K$
with the total-variation norm
$\norm{\nu}_{\mathrm{TV}}:=|\nu|(K)$. Since $\PP(K)$ is not open in this
space, define
\begin{equation}\label{eq:Udelta}
U_\delta
:=
\left\{
\nu\in\mathcal M(K):
\operatorname{dist}_{\mathrm{TV}}(\nu,\PP(K))<\delta
\right\},
\quad
\delta:=\frac{\beta_1}{2B}.
\end{equation}
Then $U_\delta$ is an open neighborhood of $\PP(K)$.

For $\nu\in U_\delta$, choose $\mu\in\PP(K)$ such that
$\norm{\nu-\mu}_{\mathrm{TV}}<\delta$. Since $0<\varphi\le B$,
\[
\begin{aligned}
|f_\nu(X_i)-f_\mu(X_i)|
&=
\left|
\int_K\varphi(X_i-\theta)\,(\nu-\mu)(d\theta)
\right|\\
&\le
B\norm{\nu-\mu}_{\mathrm{TV}}
<
B\delta
=
\frac{\beta_1}{2}.
\end{aligned}
\]
Consequently, $f_\nu(X_i)\ge\beta_1/2>0$. Thus $\Psi_\nu$ is well defined
throughout $U_\delta$, and
\[
\norm{\Psi_\nu}_\infty
\le
\frac{2B}{\beta_1}.
\]
Moreover,
\[
\norm{\nu}_{\mathrm{TV}}
\le
\norm{\mu}_{\mathrm{TV}}+\norm{\nu-\mu}_{\mathrm{TV}}
<
1+\delta.
\]

Now let $\nu,\nu'\in U_\delta$. Since
$|f_\nu(X_i)-f_{\nu'}(X_i)|\le B\norm{\nu-\nu'}_{\mathrm{TV}}$ and both
$f_\nu(X_i)$ and $f_{\nu'}(X_i)$ are at least $\beta_1/2$,
\[
\begin{aligned}
\left|
\frac1{f_\nu(X_i)}-
\frac1{f_{\nu'}(X_i)}
\right|
&=
\frac{|f_\nu(X_i)-f_{\nu'}(X_i)|}
{f_\nu(X_i)f_{\nu'}(X_i)}\\
&\le
\frac{4B}{\beta_1^2}
\norm{\nu-\nu'}_{\mathrm{TV}},
\end{aligned}
\]
and consequently
\[
\norm{\Psi_\nu-\Psi_{\nu'}}_\infty
\le
\frac{4B^2}{\beta_1^2}
\norm{\nu-\nu'}_{\mathrm{TV}}.
\]

Let
$V(\nu):=(\Psi_\nu-1)\nu$. Then
\[
V(\nu)-V(\nu')
=
(\Psi_\nu-1)(\nu-\nu')
+
(\Psi_\nu-\Psi_{\nu'})\nu',
\]
and therefore
\[
\begin{aligned}
\norm{V(\nu)-V(\nu')}_{\mathrm{TV}}
&\le
\norm{\Psi_\nu-1}_\infty
\norm{\nu-\nu'}_{\mathrm{TV}}
+
\norm{\Psi_\nu-\Psi_{\nu'}}_\infty
\norm{\nu'}_{\mathrm{TV}}\\
&\le
C\norm{\nu-\nu'}_{\mathrm{TV}},
\end{aligned}
\]
for a finite constant $C$. Also,
\[
\norm{V(\nu)}_{\mathrm{TV}}
\le
\left(1+\frac{2B}{\beta_1}\right)(1+\delta)
=:
M.
\]
Hence $V$ is bounded and Lipschitz on $U_\delta$.

Since $\mu_0\in\PP(K)$, the closed ball
$\overline B_{\mathrm{TV}}(\mu_0,\delta/2)$ is contained in $U_\delta$.
The Banach-space Picard--Lindel\"of theorem, applied on this ball, gives
$T_0>0$ and a unique local solution of
\[
\dot\mu_t=V(\mu_t),
\quad
\mu_{t=0}=\mu_0,
\]
where $T_0$ may be chosen small enough that the solution remains in the ball
(for instance, with $MT_0<\delta/2$). Let $[0,T_*)$ be the maximal interval
of existence of the resulting $U_\delta$-valued solution.

We next show that the solution is in fact a probability measure throughout its
maximal interval. Fix $T<T_*$. For $0\le t\le T$, write
\[
a_t(\theta):=\Psi_{\mu_t}(\theta)-1,
\quad
A_t(\theta):=\int_0^t a_s(\theta)\,ds.
\]
Since the solution takes values in $U_\delta$, $a_t$ is uniformly bounded on
$[0,T]\times K$. For every bounded continuous function $g$,
\[
\begin{aligned}
\frac d{dt}
\int_K g(\theta)e^{-A_t(\theta)}\,\mu_t(d\theta)
&=
-\int_K g a_t e^{-A_t}\,d\mu_t
+
\int_K g e^{-A_t}\,d\dot\mu_t\\
&=
-\int_K g a_t e^{-A_t}\,d\mu_t
+
\int_K g e^{-A_t}a_t\,d\mu_t\\
&=0.
\end{aligned}
\]
Thus
\[
e^{-A_t}\mu_t=\mu_0,
\quad\text{or equivalently}\quad
\mu_t(d\theta)=e^{A_t(\theta)}\mu_0(d\theta).
\]
In particular, $\mu_t$ is a positive measure.

It remains to check that the total mass remains $1$. For every
$\nu\in U_\delta$,
\[
\begin{aligned}
\int_K\Psi_\nu(\theta)\,\nu(d\theta)
&=
\frac1n\sum_{i=1}^n
\frac{\int_K\varphi(X_i-\theta)\,\nu(d\theta)}{f_\nu(X_i)}\\
&=1.
\end{aligned}
\]
Therefore, with $m_t:=\mu_t(K)$,
\[
\dot m_t
=
\int_K(\Psi_{\mu_t}-1)\,d\mu_t
=
1-m_t.
\]
Since $m_0=1$, uniqueness for this scalar ODE gives $m_t\equiv1$. Hence
$\mu_t\in\PP(K)$ for every $t<T_*$. Since $T<T_*$ was arbitrary, this
holds throughout the maximal interval.

We now prove that the solution is global. Since $\mu_t\in\PP(K)$ for every
$t<T_*$, we have $\norm{\Psi_{\mu_t}}_\infty\le B/\beta_1$, and therefore
\[
\norm{\dot\mu_t}_{\mathrm{TV}}
\le
1+\frac{B}{\beta_1}
=:
M_0,
\quad t<T_*.
\]
Suppose for contradiction that $T_*<\infty$. For $0\le s<t<T_*$,
\[
\norm{\mu_t-\mu_s}_{\mathrm{TV}}
\le
\int_s^t\norm{\dot\mu_r}_{\mathrm{TV}}\,dr
\le
M_0(t-s).
\]
Thus $(\mu_t)$ is Cauchy in total variation as $t\uparrow T_*$. Since
$\mathcal M(K)$ is complete, there exists $\mu_*\in\mathcal M(K)$ such that
\[
\mu_t\to\mu_*
\quad\text{in total variation as }t\uparrow T_*.
\]
Because $\PP(K)$ is closed in total variation, $\mu_*\in\PP(K)\subset U_\delta$.
Moreover, letting $t\uparrow T_*$ in the integral equation for the old solution
gives
\[
\mu_*
=
\mu_0+
\int_0^{T_*}V(\mu_s)\,ds.
\]
Picard--Lindel\"of applied at $\mu_*$ gives a local solution beginning at time
$T_*$. The preceding integral identity and uniqueness show that this solution
patches onto the original one, extending it beyond $T_*$ and contradicting
maximality. Hence $T_*=\infty$.

We now record the density formula. Since $\mu_0=p_0\,dx$, the identity
$\mu_t=e^{A_t}\mu_0$ gives
\[
p_t(\theta)
=
p_0(\theta)
\exp\left(
\int_0^t
(\Psi_{\mu_s}(\theta)-1)\,ds
\right)>0.
\]
Therefore
\[
\log\frac{p_t(\theta)}{p_0(\theta)}
=
\int_0^t(\Psi_{\mu_s}(\theta)-1)\,ds,
\]
which is \eqref{eq:density}.

Finally, differentiating the likelihood along the flow,
\[
\begin{aligned}
\frac d{dt}L_n(f_{\mu_t})
=
\frac1n\sum_{i=1}^n
\frac1{f_{\mu_t}(X_i)}
\int_K\varphi(X_i-\theta)\,\dot\mu_t(d\theta)
=
\int_K\Psi_{\mu_t}(\theta)\,\dot\mu_t(d\theta)\\
=
\int_K\Psi_{\mu_t}(\Psi_{\mu_t}-1)\,d\mu_t
=
\int_K\Psi_{\mu_t}^2\,d\mu_t-1,
\end{aligned}
\]
where we used $\int_K\Psi_{\mu_t}\,d\mu_t=1$. Since $\mu_t\in\PP(K)$,
\[
\int_K\Psi_{\mu_t}^2\,d\mu_t-1
=
\int_K\Psi_{\mu_t}^2\,d\mu_t
-
\left(\int_K\Psi_{\mu_t}\,d\mu_t\right)^2
=
\operatorname{Var}_{\mu_t}(\Psi_{\mu_t}).
\]
Thus
\[
\frac d{dt}L_n(f_{\mu_t})
=
\operatorname{Var}_{\mu_t}(\Psi_{\mu_t})
\ge0,
\]
so $t\mapsto L_n(f_{\mu_t})$ is nondecreasing.
\end{proof}

\subsection{Finite-time convergence to approximate maximizers and consequences}
\begin{proof}[\textbf{Proof of Theorem \ref{thm:fr}}]
For any $\eta > 0$, we define the entropy approximation profile $\mathcal{H}_{\mu_0}(\eta)$ of a measure $\mu_0$ as
\begin{equation}\label{eq:entropy-profile}
\mathcal H_{\mu_0}(\eta)
:=
\inf\left\{
\KL(\pi\,\|\,\mu_0):
\pi\in\PP(K),\
\Lhat-L_n(f_\pi)\le\eta
\right\},
\end{equation}
with the convention that the infimum is $+\infty$ if the admissible set is empty.

We shall use the following two lemmas regarding the entropy approximation profile, 

\begin{lemma}[Entropy approximation principle]\label{cor:entropy-profile}
For every $T>0$,
\begin{equation}\label{eq:entropy-profile-rate}
\Lhat-L_n(f_{\mu_T})
\le
\inf_{\eta>0}
\left\{
\eta+\frac{\mathcal H_{\mu_0}(\eta)}{T}
\right\}.
\end{equation}
\end{lemma}

\begin{lemma}[Gaussian entropy approximation profile]
\label{cor:gaussian-entropy-profile}
Let $\mu_0=\mathrm{Unif}(B(X_1,\rho_n+1))$. Then, for
$0<\eta\le3/4$,
\begin{equation}\label{eq:gaussian-entropy-profile}
\mathcal H_{\mu_0}(\eta)
\le
 d\log\frac{7\rho_n^2}{\sqrt\eta}.
\end{equation}
\end{lemma}

Now, take $\eta:=\eps/2$. Since $\eps\le1$, we have $\eta\le1/2$, and
lemma~\ref{cor:gaussian-entropy-profile} gives
\[
\mathcal H_{\mu_0}(\eps/2)
\le
 d\log\frac{7\sqrt2\,\rho_n^2}{\sqrt\eps}
\le
 d\log\frac{10\rho_n^2}{\sqrt\eps}.
\]
Therefore lemma~\ref{cor:entropy-profile} yields
\[
\Lhat-L_n(f_{\mu_T})
\le
\frac{\eps}{2}
+
\frac{d}{T}\log\frac{10\rho_n^2}{\sqrt\eps}.
\]
If
$T\ge2d\eps^{-1}\log(10\rho_n^2/\sqrt\eps)$, the second term is at most
$\eps/2$, proving the result.
\end{proof}

We now give a proof of the two lemmas that we had used in the proof of this result. Each lemma is, in fact, a corollary of more general propositions, which may be of independent interest. 

\begin{proposition}[Fisher--Rao comparator inequality]\label{prop:frrate}
Let $(\mu_t)_{t\ge0}$ be as in Lemma~\ref{lem:wellposed}, and let
$\pi\in\PP(K)$ satisfy
\[
\KL(\pi\,\|\,\mu_0)<\infty.
\]
Then, for every $T>0$,
\begin{equation}\label{eq:fr-integrated}
\int_0^T
\bigl(L_n(f_\pi)-L_n(f_{\mu_t})\bigr)\,dt
\le
\KL(\pi\,\|\,\mu_0).
\end{equation}
Consequently,
\begin{equation}\label{eq:fr-comparator}
L_n(f_\pi)-L_n(f_{\mu_T})
\le
\frac{\KL(\pi\,\|\,\mu_0)}{T}.
\end{equation}
\end{proposition}

\begin{proof}
We first bound the likelihood gap by the Fisher--Rao first variation. For every
$\mu\in\PP(K)$ and every $i$, concavity of the logarithm gives
\[
\log f_\pi(X_i)-\log f_\mu(X_i)
\le
\frac{f_\pi(X_i)-f_\mu(X_i)}{f_\mu(X_i)}.
\]
Averaging over $i$ yields

\begin{equation}
\begin{aligned}
L_n(f_\pi)-L_n(f_\mu)
&\le
\frac1n\sum_{i=1}^n
\left(
\frac{f_\pi(X_i)}{f_\mu(X_i)}-1
\right)\notag\\
&=
\int_K\Psi_\mu(\theta)\,\pi(d\theta)-1\notag\\
&=
\int_K\bigl(\Psi_\mu(\theta)-1\bigr)\,\pi(d\theta).
\label{eq:concavity}
\end{aligned}
\end{equation}

We next compute the evolution of the relative entropy with respect to the flow.
Since $\KL(\pi\,\|\,\mu_0)<\infty$, we have $\pi\ll\mu_0$. By
Lemma~\ref{lem:wellposed},
\[
\frac{d\mu_t}{d\mu_0}(\theta)
=
\exp\left(
\int_0^t
\bigl(\Psi_{\mu_s}(\theta)-1\bigr)\,ds
\right)>0.
\]
Thus $\mu_t$ and $\mu_0$ are mutually absolutely continuous, and in particular
$\pi\ll\mu_t$. Therefore
$\frac{d\pi}{d\mu_t}=\frac{d\pi}{d\mu_0}\frac{d\mu_0}{d\mu_t}$, and hence

\begin{equation*}
\begin{aligned}
 \KL(\pi\,\|\,\mu_t)
&=
\int_K
\log\left(\frac{d\pi}{d\mu_t}\right)\,d\pi\\
&=
\int_K
\log\left(\frac{d\pi}{d\mu_0}\right)\,d\pi
-
\int_K
\log\left(\frac{d\mu_t}{d\mu_0}\right)\,d\pi\\
&=
\KL(\pi\,\|\,\mu_0)
-
\int_K
\int_0^t
\bigl(\Psi_{\mu_s}(\theta)-1\bigr)\,ds\,\pi(d\theta).
\end{aligned}   
\end{equation*}

Since $\Psi_{\mu_s}$ is uniformly bounded on $K$ for $s$ in bounded time
intervals, Fubini's theorem gives
\begin{equation}\label{eq:KL-flow}
\KL(\pi\,\|\,\mu_t)
=
\KL(\pi\,\|\,\mu_0)
-
\int_0^t
\int_K
\bigl(\Psi_{\mu_s}-1\bigr)\,d\pi\,ds.
\end{equation}

Applying \eqref{eq:concavity} with $\mu=\mu_s$ gives
\[
\int_K(\Psi_{\mu_s}-1)\,d\pi
\ge
L_n(f_\pi)-L_n(f_{\mu_s}).
\]
Substituting this into \eqref{eq:KL-flow} at time $T$,
\[
\KL(\pi\,\|\,\mu_T)
\le
\KL(\pi\,\|\,\mu_0)
-
\int_0^T
\bigl(L_n(f_\pi)-L_n(f_{\mu_s})\bigr)\,ds.
\]
Since relative entropy is nonnegative, this proves \eqref{eq:fr-integrated}.
Finally, by Lemma~\ref{lem:wellposed}, $s\mapsto L_n(f_{\mu_s})$ is
nondecreasing. Hence, for every $0\le s\le T$,
\[
L_n(f_\pi)-L_n(f_{\mu_s})
\ge
L_n(f_\pi)-L_n(f_{\mu_T}).
\]
Substituting this into \eqref{eq:fr-integrated} yields
\eqref{eq:fr-comparator}.
\end{proof}

\begin{proof}[\textbf{Proof of Lemma 3.1}]
If $\pi\in\PP(K)$ satisfies
$\Lhat-L_n(f_\pi)\le\eta$ and
$\KL(\pi\,\|\,\mu_0)<\infty$, then Proposition~\ref{prop:frrate} gives
\[
\begin{aligned}
\Lhat-L_n(f_{\mu_T})
&=
\bigl(\Lhat-L_n(f_\pi)\bigr)
+
\bigl(L_n(f_\pi)-L_n(f_{\mu_T})\bigr)\\
&\le
\eta+
\frac{\KL(\pi\,\|\,\mu_0)}{T}.
\end{aligned}
\]
Taking the infimum first over such $\pi$ and then over $\eta>0$ proves the
claim.
\end{proof}

\begin{proposition}[Mollification]\label{lem:mollify}
Let $\mu\in\PP\bigl(B(X_1,\rho_n)\bigr)$, let
$0<\sigma\le1/(4\rho_n)$, and define
\[
\pi:=\mu*\mathrm{Unif}\bigl(B(0,\sigma)\bigr).
\]
Then $\supp(\pi)\subseteq B(X_1,\rho_n+1)$, $\pi$ has a density bounded by
$\abs{B(0,\sigma)}^{-1}$, and
\[
\bigl\lvert L_n(f_\pi)-L_n(f_\mu)\bigr\rvert
\le
12\sigma^2\rho_n^2.
\] 
\end{proposition}

\begin{proof}
Let $U_\sigma:=\mathrm{Unif}(B(0,\sigma))$. Since
$\supp(\mu)\subseteq B(X_1,\rho_n)$ and $\supp(U_\sigma)\subseteq B(0,\sigma)$,
\[
\supp(\pi)
\subseteq
B(X_1,\rho_n+\sigma)
\subseteq
B(X_1,\rho_n+1),
\]
where we used $\sigma\le1/(4\rho_n)\le1$. Moreover, since $U_\sigma$ has
density $\abs{B(0,\sigma)}^{-1}\mathbf 1_{B(0,\sigma)}$, the convolution
$\pi=\mu*U_\sigma$ has density
\[
p_\pi(x)
=
\frac{1}{\abs{B(0,\sigma)}}
\int\mathbf 1_{\{\norm{x-\theta}\le\sigma\}}\,\mu(d\theta)
\le
\frac{1}{\abs{B(0,\sigma)}}.
\]

It remains to control the likelihood. Define
\[
\bar\varphi_\sigma(z)
:=
\frac{1}{\abs{B(0,\sigma)}}
\int_{B(0,\sigma)}\varphi(z-v)\,dv.
\]
By Fubini,
\[
f_\pi(X_i)
=
\int\bar\varphi_\sigma(X_i-\theta)\,\mu(d\theta).
\]
We first compare $\bar\varphi_\sigma(z)$ with $\varphi(z)$. For fixed
$z,v\in\R^d$, Taylor's theorem with integral remainder gives
\[
\varphi(z-v)-\varphi(z)
=
-\langle\nabla\varphi(z),v\rangle
+
\int_0^1(1-s)
\left\langle
\nabla^2\varphi(z-sv)v,v
\right\rangle\,ds.
\]
Averaging over $v\sim U_\sigma$, the linear term vanishes because $U_\sigma$
has mean zero. Hence
\[
\bar\varphi_\sigma(z)-\varphi(z)
=
\frac{1}{\abs{B(0,\sigma)}}
\int_{B(0,\sigma)}
\int_0^1(1-s)
\left\langle
\nabla^2\varphi(z-sv)v,v
\right\rangle\,ds\,dv.
\]
Since $\norm v\le\sigma$ and $\int_0^1(1-s)\,ds=1/2$,
\[
\abs{\bar\varphi_\sigma(z)-\varphi(z)}
\le
\frac{\sigma^2}{2}
\sup_{\norm v\le\sigma}
\norm{\nabla^2\varphi(z-v)}.
\]

For the standard Gaussian,
$\nabla^2\varphi(y)=\varphi(y)(yy^\top-I)$, and therefore
$\norm{\nabla^2\varphi(y)}\le\varphi(y)(\norm y^2+1)$. Now take
$z=X_i-\theta$ with $X_i,\theta\in B(X_1,\rho_n)$. Then
$\norm z\le2\rho_n$. Since
$\sigma\le1/(4\rho_n)\le1/4\le\rho_n/4$, for every $\norm v\le\sigma$,
\[
\norm{z-v}^2+1
\le
\left(2\rho_n+\frac{\rho_n}{4}\right)^2+1
\le
7\rho_n^2.
\]
Also,
\[
\frac{\varphi(z-v)}{\varphi(z)}
=
\exp\left(
\langle z,v\rangle-\frac{\norm v^2}{2}
\right)
\le
e^{\norm z\,\norm v}
\le
e^{2\rho_n\sigma}
\le
e^{1/2}.
\]
It follows that
\[
\sup_{\norm v\le\sigma}
\norm{\nabla^2\varphi(z-v)}
\le
7e^{1/2}\rho_n^2\varphi(z).
\]
Consequently,
\[
\abs{\bar\varphi_\sigma(z)-\varphi(z)}
\le
\frac{7e^{1/2}}{2}\sigma^2\rho_n^2\varphi(z)
\le
6\sigma^2\rho_n^2\varphi(z).
\]

Applying this with $z=X_i-\theta$ and integrating against $\mu$ gives
\[
\begin{aligned}
\abs{f_\pi(X_i)-f_\mu(X_i)}
&\le
\int
\abs{
\bar\varphi_\sigma(X_i-\theta)
-
\varphi(X_i-\theta)
}\,\mu(d\theta)\\
&\le
6\sigma^2\rho_n^2 f_\mu(X_i).
\end{aligned}
\]
Thus, with $\gamma':=6\sigma^2\rho_n^2$, we have
\[
\abs{f_\pi(X_i)-f_\mu(X_i)}
\le
\gamma' f_\mu(X_i).
\]
Since $\sigma\le1/(4\rho_n)$,
$\gamma'\le6/16=3/8<1/2$. Lemma~\ref{lem:log} therefore yields
\[
\abs{\log f_\pi(X_i)-\log f_\mu(X_i)}
\le
2\gamma'
=
12\sigma^2\rho_n^2.
\]
Averaging over $i=1,\dots,n$ proves
\[
\bigl\lvert L_n(f_\pi)-L_n(f_\mu)\bigr\rvert
\le
12\sigma^2\rho_n^2.
\]
\end{proof}

\begin{lemma}\label{lem:log}
For all $u,v>0$, $\abs{\log u-\log v}\le\abs{u-v}/\min(u,v)$. Consequently, if
$\abs{u-v}\le\gamma v$ with $\gamma\le\frac12$, then
$\abs{\log u-\log v}\le 2\gamma$.
\end{lemma}

\begin{proof}
Assume $v\le u$; then $\log u-\log v=\int_v^{u}s^{-1}\,ds\le(u-v)/v$. For the
consequence, $u\ge(1-\gamma)v$, so the bound is
$\gamma v/((1-\gamma)v)=\gamma/(1-\gamma)\le2\gamma$.
\end{proof}

\begin{proof}[\textbf{Proof of Lemma 3.2}]
Choose an NPMLE $\muhat$, which is supported on
$B(X_1,\rho_n)$, and put
$\sigma:=\rho_n^{-1}\sqrt{\eta/12}$. Since $\eta\le3/4$, we have
$\sigma\le1/(4\rho_n)$. For
\[
\pi_\eta
:=
\muhat*\mathrm{Unif}(B(0,\sigma)),
\]
Lemma~\ref{lem:mollify} gives
$\Lhat-L_n(f_{\pi_\eta})\le12\sigma^2\rho_n^2=\eta$.
Moreover, $\pi_\eta$ is supported on $K$ and has density at most
$\abs{B(0,\sigma)}^{-1}$. Since $\mu_0$ has constant density
$\abs{B(0,\rho_n+1)}^{-1}$ on $K$,
\[
\begin{aligned}
\KL(\pi_\eta\,\|\,\mu_0)
\le
\log\frac{\abs{B(0,\rho_n+1)}}{\abs{B(0,\sigma)}}
=
 d\log\frac{\rho_n+1}{\sigma}\\
\le
 d\log\frac{2\sqrt{12}\,\rho_n^2}{\sqrt\eta}
\le
 d\log\frac{7\rho_n^2}{\sqrt\eta}.
\end{aligned}
\]
The definition \eqref{eq:entropy-profile} now gives
\eqref{eq:gaussian-entropy-profile}.
\end{proof}

\begin{proof}[\textbf{Proof of Corollary \ref{cor:fr-direct-rate}}]
In \eqref{eq:entropy-profile-rate}, take $\eta=d/(2T)$. The assumption
$T\ge2d/3$ ensures $\eta\le3/4$, so
Lemma~\ref{cor:gaussian-entropy-profile} gives
\[
\begin{aligned}
\Lhat-L_n(f_{\mu_T})
&\le
\frac{d}{2T}
+
\frac{d}{T}
\log\left(
7\rho_n^2\sqrt{\frac{2T}{d}}
\right)\\
&\le
\frac{d}{T}
\left(
\frac12+
\log\left(
10\rho_n^2\sqrt{\frac{T}{d}}
\right)
\right),
\end{aligned}
\]
since $7\sqrt2<10$.
\end{proof}

\begin{proof}[\textbf{Proof of Corollary \ref{cor:fr-cluster}}]
The likelihood convergence follows from Corollary~\ref{cor:fr-direct-rate}
together with the monotonicity in Lemma~\ref{lem:wellposed}. Since $K$ is
compact, $\PP(K)$ is weakly compact. Thus every sequence $t_j\to\infty$ admits
a subsequence, still denoted $t_j$, such that
$\mu_{t_j}\Rightarrow\bar\mu$ for some $\bar\mu\in\PP(K)$.

For each $i$, the map
$\theta\mapsto\varphi(X_i-\theta)$ is bounded and continuous on $K$, so
$f_{\mu_{t_j}}(X_i)\to f_{\bar\mu}(X_i)$. Hence
\[
L_n(f_{\bar\mu})
=
\lim_{j\to\infty}L_n(f_{\mu_{t_j}})
=
\Lhat.
\]
Thus $\bar\mu$ is an NPMLE. If the NPMLE is unique, every weakly convergent
subsequence has the same limit $\muhat$, which implies
$\mu_t\Rightarrow\muhat$.
\end{proof}

\begin{proof}[\textbf{Proof of Corollary \ref{cor:fr-fitted-values}}]
The set
\[
\mathcal U
:=
\left\{
\bigl(f_\mu(X_1),\dots,f_\mu(X_n)\bigr):
\mu\in\PP(K)
\right\}
\]
is convex, and
$G(u):=n^{-1}\sum_{i=1}^n\log u_i$ is strictly concave on $\mathcal U$.
Hence its maximizing fitted-value vector $u^\star$ is unique. Since every
coordinate in $\mathcal U$ is at most $B$,
\[
\nabla^2G(u)
=
-\frac1n\operatorname{diag}(u_1^{-2},\dots,u_n^{-2})
\preceq
-\frac{1}{nB^2}I.
\]
Thus $G$ is $1/(nB^2)$-strongly concave on $\mathcal U$. The first-order
optimality condition at $u^\star$ therefore gives
\[
\Lhat-L_n(f_{\mu_t})
=
G(u^\star)-G(u_t)
\ge
\frac{1}{2nB^2}\norm{u_t-u^\star}_2^2,
\]
which proves \eqref{eq:fitted-value-bound}. The stated rate follows from
Corollary~\ref{cor:fr-direct-rate}.
\end{proof}

\begin{proof}[\textbf{Proof of Corollary \ref{cor:fr-strict}}]
If the derivative vanished at some finite time $t$, then
Lemma~\ref{lem:wellposed} would give
$\operatorname{Var}_{\mu_t}(\Psi_{\mu_t})=0$. Since $\mu_t$ has a strictly
positive density on $K$ and $\Psi_{\mu_t}$ is continuous,
$\Psi_{\mu_t}\equiv1$ on $K$. For fixed $\mu_t$, however,
$\Psi_{\mu_t}$ is a real-analytic function of $\theta$. Since $K$ has nonempty
interior, the identity theorem implies $\Psi_{\mu_t}\equiv1$ on $\R^d$. This is
impossible because $\Psi_{\mu_t}(\theta)\to0$ as $\norm\theta\to\infty$.
\end{proof}

\begin{proof}[\textbf{Proof of Proposition \ref{prop:general-fr}}]
The proof is the same as that of Lemma~\ref{lem:wellposed} and
Proposition~\ref{prop:frrate}. The lower bound $F_i(\mu)\ge b$ uniformly on
$\PP(\Theta)$ makes the vector field bounded and locally Lipschitz in total
variation. The integrating-factor argument gives positivity and the given
Radon-Nikodym formula, while $\int\Psi_\mu\,d\mu=1$ gives mass preservation and
the dissipation identity. Finally, concavity of the logarithm and the
relative-entropy identity give the comparator bound.
\end{proof}

\begin{proof}[\textbf{Proof of Corollary \ref{cor:fr-affine}}]
By classical theory, the NPMLE is supported on
$\conv(X_1,\dots,X_n)\subseteq A$. We may therefore regard the mixing problem as
a problem over $\PP(A)$.

Let $V_A:=A-X_1$, an $s$-dimensional linear subspace, and mollify an NPMLE
within $V_A$:
\[
\pi_\sigma
:=
\muhat*\mathrm{Unif}\bigl(B_{V_A}(0,\sigma)\bigr).
\]
The proof of Proposition~\ref{lem:mollify} is unchanged, since it uses only the
Euclidean norm and the Hessian bound for the Gaussian kernel. Thus
\[
\Lhat-L_n(f_{\pi_\sigma})
\le
12\sigma^2\rho_n^2.
\]
The entropy calculation now involves $s$-dimensional volumes, and therefore
\[
\KL(\pi_\sigma\,\|\,\mu_0)
\le
s\log\frac{\rho_n+1}{\sigma}.
\]
The proof of Theorem~\ref{thm:fr} now applies verbatim with $d$ replaced by $s$.
\end{proof}

\section{Proofs of Results in Section \ref{sec:Comp}}\label{sec:Comp proofs}
The notation of Section~2.2 is retained. In particular,
$b=B/(4n)$, and $\widehat L_n$ denotes the unrestricted maximum
log likelihood. All norms of vectors are Euclidean unless another
norm is indicated.

\subsection{The exponential family and the modified likelihood}

\begin{proof}[Proof of Proposition~\ref{prop:fc-representation}]
Let $p_t=d\mu_t/d\mu_0$. By the Fisher--Rao equation,
\[
\frac{\partial}{\partial t}\log p_t(\theta)
=\Psi_{\mu_t}(\theta)-1
=\frac1n\sum_{i=1}^n
\frac{\phi(X_i-\theta)}{f_{\mu_t}(X_i)}-1.
\]
Integrating from $0$ to $t$, and using $p_0=1$, gives
\[
\log p_t(\theta)
=-t+\sum_{i=1}^n\phi(X_i-\theta)
\frac1n\int_0^t\frac{ds}{f_{\mu_s}(X_i)}.
\]
Setting
$
a_i(t)=\frac1n\int_0^t\frac{ds}{f_{\mu_s}(X_i)},
\quad 1\le i\le n,
$
we obtain
$
\frac{d\mu_t}{d\mu_0}(\theta)
=\exp\{V_{a(t)}(\theta)-t\}.
$
Since $\mu_t$ is a probability measure,
\[
1=\int_K e^{V_{a(t)}(\theta)-t}\,\mu_0(d\theta)
=e^{-t}Z(a(t)),
\]
and therefore $Z(a(t))=e^t$ and
$\mu_t=\mu_{a(t)}$. Differentiating $a_i(t)$ now gives
\[
\dot a_i(t)=\frac1{nF_i(a(t))},
\quad a_i(0)=0,
\]
which proves \eqref{eq:fc-flow}. We next consider the proximal problem
\eqref{eq:fc-prox-original}. For $a\in\mathbb R_+^n$
and $\tau>0$, define
\[
H(u)=\log Z(a+u)-\frac{\tau}{n}
\sum_{i=1}^n\log u_i,
\quad u\in(0,\infty)^n.
\]
The function $H$ is strictly convex, since $\log Z$
is convex and $-\log u_i$ is strictly convex.
Moreover, $H(u)\to\infty$ as any coordinate $u_i$
tends to zero while the remaining coordinates are bounded.

To establish coercivity, let
$$
c=\min_{\substack{1\le i\le n\\\theta\in K}}
\phi(X_i-\theta)>0.
$$
Since $u_i\ge0$, we have
\[
\begin{aligned}
Z(a+u)
=\int_K
e^{V_a(\theta)+\sum_{i=1}^n u_i\phi(X_i-\theta)}
\,\mu_0(d\theta)
\ge e^{c\|u\|_{TV}}Z(a).
\end{aligned}
\]
Consequently,
\[
\begin{aligned}
H(u)
\ge\log Z(a)+c\|u\|_{TV}
-\frac{\tau}{n}\sum_{i=1}^n\log u_i
\ge\log Z(a)+c\|u\|_{TV}
-\tau\log\|u\|_{TV}.
\end{aligned}
\]
Thus $H(u)\to\infty$ as $\|u\|_{TV}\to\infty$,
and $H$ admits a unique minimizer in $(0,\infty)^n$.

Moreover, 
\[
\frac{\partial}{\partial u_i}\log Z(a+u)
=\frac{\int_K\phi(X_i-\theta)
e^{V_{a+u}(\theta)}\,\mu_0(d\theta)}
{Z(a+u)}
=F_i(a+u).
\]
The first-order condition for the minimizer is therefore
\[
0=\frac{\partial H}{\partial u_i}
=F_i(a+u)-\frac{\tau}{nu_i},
\quad 1\le i\le n.
\]
Setting $a^+=a+u$, we obtain
$
a_i^+=a_i+\frac{\tau}{nF_i(a^+)},
\quad 1\le i\le n,
$
which is \eqref{eq:fc-original-update}.

It remains to show that $\mu_{a^+}$ is the unique
maximizer of \eqref{eq:fc-prox-original}.
Let $\nu\in\mathcal P(K)$ satisfy
$\operatorname{KL}(\nu\|\mu_a)<\infty$.
By concavity of the logarithm,
\[
\begin{aligned}
L_n(f_\nu)-L_n(f_{\mu_{a^+}})
\le\frac1n\sum_{i=1}^n
\frac{f_\nu(X_i)-f_{\mu_{a^+}}(X_i)}
{f_{\mu_{a^+}}(X_i)}
=\int_K\Psi_{\mu_{a^+}}(\theta)
\,(\nu-\mu_{a^+})(d\theta).
\end{aligned}
\]
On the other hand, the definition of the exponential family
and \eqref{eq:fc-original-update} give
\[
\begin{aligned}
\log\frac{d\mu_{a^+}}{d\mu_a}(\theta)
=\sum_{i=1}^n(a_i^+-a_i)\phi(X_i-\theta)
-\log\frac{Z(a^+)}{Z(a)}
=\tau\Psi_{\mu_{a^+}}(\theta)
-\log\frac{Z(a^+)}{Z(a)}.
\end{aligned}
\]
Since $\nu$ and $\mu_{a^+}$ are probability measures,
the constant term vanishes upon integration against
$\nu-\mu_{a^+}$. Hence
\[
\begin{aligned}
\tau\int_K\Psi_{\mu_{a^+}}(\theta)
\,(\nu-\mu_{a^+})(d\theta)
=\int_K
\log\frac{d\mu_{a^+}}{d\mu_a}
\,d(\nu-\mu_{a^+})
\\
\quad=\operatorname{KL}(\nu\|\mu_a)
-\operatorname{KL}(\nu\|\mu_{a^+})
-\operatorname{KL}(\mu_{a^+}\|\mu_a).
\end{aligned}
\]
Combining the preceding inequalities yields
\[
\begin{aligned}
L_n(f_\nu)-L_n(f_{\mu_{a^+}})
-\frac1\tau\left\{
\operatorname{KL}(\nu\|\mu_a)
-\operatorname{KL}(\mu_{a^+}\|\mu_a)
\right\}
\le
-\frac1\tau\operatorname{KL}(\nu\|\mu_{a^+})
\le0.
\end{aligned}
\]
For measures satisfying
$\operatorname{KL}(\nu\|\mu_a)=\infty$, the proximal
objective is $-\infty$. Since
$\operatorname{KL}(\nu\|\mu_{a^+})=0$ if and only if
$\nu=\mu_{a^+}$, the preceding inequality proves that
$\mu_{a^+}$ is the unique maximizer of
\eqref{eq:fc-prox-original}.
\end{proof}

We next consider the proximal problem with $L_n$ replaced by
$\mathcal L_b$. Since
\[
h_b'(u)=\frac1{\max\{u,b\}},
\]
the first variation of $\mathcal L_b$ is
\begin{equation}\label{eq:fp-capped-score}
\Psi_b[\mu](\theta)
=\frac1n\sum_{i=1}^n
\frac{\phi(X_i-\theta)}{\max\{f_\mu(X_i),b\}}.
\end{equation}
Consider the minimization of the function $H$ introduced
in the preceding proof, now subject to
\[
0<u_i\le\frac{\tau}{nb},
\quad 1\le i\le n.
\]
Strict convexity and the fact that $H$ diverges as any
coordinate approaches zero imply the existence of a
unique minimizer. The coordinate optimality conditions are given by
\[
\begin{cases}
F_i(a+u)=\dfrac{\tau}{nu_i},
&u_i<\dfrac{\tau}{nb},\\[6pt]
F_i(a+u)\le b,
&u_i=\dfrac{\tau}{nb}.
\end{cases}
\]
Combining the two cases, we get
\[
u_i=\frac{\tau}{n\max\{F_i(a+u),b\}},
\quad 1\le i\le n.
\]
Setting $a^+=a+u$ proves \eqref{eq:fc-capped-update}.
Moreover,
\[
\log\frac{d\mu_{a^+}}{d\mu_a}
=\tau\Psi_b[\mu_{a^+}]
-\log\frac{Z(a^+)}{Z(a)}.
\]
By concavity of $\mathcal L_b$, the relative-entropy
argument in the proof of Proposition~\ref{prop:fc-representation}
gives
\[
\begin{aligned}
\mathcal L_b(\nu)-\mathcal L_b(\mu_{a^+})
-\frac1\tau\left\{
\operatorname{KL}(\nu\|\mu_a)
-\operatorname{KL}(\mu_{a^+}\|\mu_a)
\right\}
\le
-\frac1\tau\operatorname{KL}(\nu\|\mu_{a^+})\le0
\end{aligned}
\]
for every $\nu$ with finite $\operatorname{KL}(\nu\|\mu_a)$.
Thus $\mu_{a^+}$ is the unique maximizer of the modified
proximal problem for every $\tau>0$.

\begin{lemma}[Equivalence of maximizers and output correction]
\label{lem:fp-equivalence}
Let $\mu^\star$ be an NPMLE and set
$f_i^\star=f_{\mu^\star}(X_i)$, $1\le i\le n$.
The vector $(f_1^\star,\ldots,f_n^\star)$ is independent
of the choice of $\mu^\star$ and satisfies
$f_i^\star\ge B/n=4b$. Moreover,
\[
\operatorname*{arg\,max}_{\mu\in\mathcal P(K)}
\mathcal L_b(\mu)
=
\operatorname*{arg\,max}_{\mu\in\mathcal P(K)}
L_n(f_\mu),
\quad
\max_{\mu\in\mathcal P(K)}\mathcal L_b(\mu)
=\widehat L_n.
\]
\\
Now let $\mu\in\mathcal P(K)$ and suppose that
$\widehat f_1,\ldots,\widehat f_n$ satisfy
$
\max_{1\le i\le n}
|\widehat f_i-f_\mu(X_i)|\le\delta,
$ for $0\le\delta\le b/2.$

Define
\[
I=\{i:\widehat f_i<b+\delta\},
\quad
\widetilde\mu=
\left(1-\frac{|I|}{n}\right)\mu
+\frac1n\sum_{i\in I}
\operatorname{Unif}(B(X_i,1)).
\]
Then
\begin{equation}\label{eq:fp-correction}
\begin{aligned}
\frac{|I|}{n}
&\le
\frac{\widehat L_n-\mathcal L_b(\mu)}
{\log 2-1/2},\\
\widehat L_n-L_n(f_{\widetilde\mu})
&\le
\left(1+\frac1{\log 2-1/2}\right)
\bigl(\widehat L_n-\mathcal L_b(\mu)\bigr).
\end{aligned}
\end{equation}
\end{lemma}

\begin{proof}
Let $\mu^\star$ be an NPMLE. Projection onto
$\operatorname{conv}\{X_1,\ldots,X_n\}$ cannot increase the
distance to any observation. Consequently, pushing a
probability measure forward under this projection cannot
decrease its likelihood. Since the convex hull is compact,
an unrestricted maximizer exists. Moreover, strict concavity
of the function
\[
(f_1,\ldots,f_n)\longmapsto
\frac1n\sum_{i=1}^n\log f_i
\]
implies that the optimal fitted vector is unique.

The first-order optimality condition, obtained by considering
the perturbation $(1-t)\mu^\star+t\delta_\theta$, gives
\begin{equation}\label{eq:fp-optimality}
\frac1n\sum_{i=1}^n
\frac{\phi(X_i-\theta)}{f_i^\star}
\le1,
\quad \theta\in\mathbb R^d.
\end{equation}
Taking $\theta=X_i$, and retaining the $i$-th term in
the sum, we obtain $f_i^\star\ge B/n=4b$.

We next show that $\mathcal L_b$ has the same maximizers
as $L_n$. For $0\le p\le B$ and $q\ge4b$, define
\[
D_b(p,q)
=h_b(q)-h_b(p)-\frac{q-p}{q}.
\]
Since $h_b(q)=\log q$, we have
\[
D_b(p,q)=
\begin{cases}
p/q-1-\log(p/q),&p\ge b,\\
\log(q/b)-p/b+p/q,&p<b.
\end{cases}
\]
For $p\ge b$, the inequality
$r-1-\log r\ge0$, with equality only at $r=1$,
gives $D_b(p,q)\ge0$. For $p<b$, the second
expression is decreasing in $p$, and hence
$D_b(p,q)\ge D_b(b,q)>0.$ Furthermore, if $p\le q/2$ and $p\ge b$, then
\[
D_b(p,q)
\ge\min_{0<r\le1/2}\{r-1-\log r\}
=\log2-\frac12.
\]
The same bound holds for $p<b$, since
$q\ge4b$ and
\[
D_b(p,q)\ge D_b(b,q)
\ge\log4-\frac34>\log2-\frac12.
\]
We have therefore established
\begin{equation}\label{eq:fp-scalar-divergence}
\begin{gathered}
D_b(p,q)\ge0,\quad
D_b(p,q)=0\ \Longleftrightarrow\ p=q,\\
p\le q/2\ \Longrightarrow\
D_b(p,q)\ge\log2-\frac12.
\end{gathered}
\end{equation}

For any $\mu\in\mathcal P(K)$, the definition of $D_b$
and \eqref{eq:fp-optimality} give
\begin{equation}\label{eq:fp-gap-decomposition}
\begin{aligned}
\widehat L_n-\mathcal L_b(\mu)
=\frac1n\sum_{i=1}^n
D_b(f_\mu(X_i),f_i^\star)
+1-\frac1n\sum_{i=1}^n
\frac{f_\mu(X_i)}{f_i^\star}
\ge\frac1n\sum_{i=1}^n
D_b(f_\mu(X_i),f_i^\star).
\end{aligned}
\end{equation}
At an NPMLE, $f_i^\star\ge4b$, so
$\mathcal L_b(\mu^\star)=\widehat L_n$.
Conversely, equality in \eqref{eq:fp-gap-decomposition}
requires $D_b(f_\mu(X_i),f_i^\star)=0$ for every $i$,
and hence $f_\mu(X_i)=f_i^\star$.
This proves the equivalence of the maximizers.

We now consider the correction of an approximate maximizer.
For every $i\in I$, the approximation error gives
$
f_\mu(X_i)<b+2\delta\le2b\le f_i^\star/2.
$
Applying \eqref{eq:fp-scalar-divergence} and
\eqref{eq:fp-gap-decomposition}, we obtain
\[
\widehat L_n-\mathcal L_b(\mu)
\ge\frac{|I|}{n}\left(\log2-\frac12\right),
\]
which proves the first inequality in \eqref{eq:fp-correction}.

To establish the second inequality, observe that
\[
\begin{aligned}
f_{\operatorname{Unif}(B(X_i,1))}(X_i)
=\frac B{v_d}
\int_{\|u\|\le1}e^{-\|u\|^2/2}\,du
\ge Be^{-1/2}>\frac B2.
\end{aligned}
\]
Thus, for every $i\in I$,
$
f_{\widetilde\mu}(X_i)\ge B/2n =2b.
$
Since $f_\mu(X_i)\le2b$ and $h_b$ is increasing,
$
\log f_{\widetilde\mu}(X_i)
\ge\log(2b)\ge h_b(f_\mu(X_i)).
$

For $i\notin I$, we have $f_\mu(X_i)\ge b$.
Moreover,
\[
f_{\widetilde\mu}(X_i)
\ge\left(1-\frac{|I|}{n}\right)f_\mu(X_i),
\]
and consequently,
$
\log f_{\widetilde\mu}(X_i)
\ge h_b(f_\mu(X_i))
+\log\left(1-\frac{|I|}{n}\right).
$
Summing these inequalities gives
\[
\begin{aligned}
L_n(f_{\widetilde\mu})
\ge\mathcal L_b(\mu)
+\left(1-\frac{|I|}{n}\right)
\log\left(1-\frac{|I|}{n}\right)
\ge\mathcal L_b(\mu)-\frac{|I|}{n},
\end{aligned}
\]
where we use $0\log0=0$ if $|I|=n$.
Combining this with the first inequality in
\eqref{eq:fp-correction} yields
\[
\widehat L_n-L_n(f_{\widetilde\mu})
\le
\left(1+\frac1{\log2-1/2}\right)
\bigl(\widehat L_n-\mathcal L_b(\mu)\bigr),
\]
as required.
\end{proof}

\begin{lemma}[Bounds for the score]
\label{lem:fp-score}
For every $\mu\in\mathcal P(K)$, the modified score
$\Psi_b[\mu]$ for $\mathcal{L}_b$ satisfies
\begin{equation}\label{eq:fp-score-gap}
1\le\|\Psi_b[\mu]\|_\infty
\le 2+\frac{4n}{\log2-1/2}
\bigl(\widehat L_n-\mathcal L_b(\mu)\bigr),
\end{equation}
where the supremum is taken over $\mathbb R^d$.
Moreover,
\begin{equation}\label{eq:fp-absolute-score}
\begin{aligned}
&\sup_{\theta\in\mathbb R^d}
\frac1n\sum_{i=1}^n\phi(X_i-\theta)
\left|
\frac1{\max\{f_\mu(X_i),b\}}
-\frac1{f_i^\star}
\right|
\\
&\le
3\sqrt{n\bigl(\widehat L_n-\mathcal L_b(\mu)\bigr)}
+\frac{4n}{\log2-1/2}
\bigl(\widehat L_n-\mathcal L_b(\mu)\bigr).
\end{aligned}
\end{equation}
\end{lemma}

\begin{proof}
We first prove \eqref{eq:fp-score-gap}. For $u\ge0$ and
$0<a\le1/4$, we have
$
1/{\max\{u,a\}}\ge2-u.
$
Since $f_i^\star\ge4b$, applying this inequality with
$u=f_\mu(X_i)/f_i^\star$ and $a=b/f_i^\star$ gives
\[
\begin{aligned}
\int\Psi_b[\mu](\theta)\,\mu^\star(d\theta)
=\frac1n\sum_{i=1}^n
\frac{f_i^\star}{\max\{f_\mu(X_i),b\}}
\ge2-\frac1n\sum_{i=1}^n
\frac{f_\mu(X_i)}{f_i^\star}
\ge1,
\end{aligned}
\]
where the last inequality follows by integrating
\eqref{eq:fp-optimality} against $\mu$.
Consequently, $\|\Psi_b[\mu]\|_\infty\ge1$.

For the upper bound, let
\[
J=\{i:f_\mu(X_i)<f_i^\star/2\}.
\]
By \eqref{eq:fp-scalar-divergence} and
\eqref{eq:fp-gap-decomposition},
\[
\begin{aligned}
\widehat L_n-\mathcal L_b(\mu)
\ge\frac1n\sum_{i\in J}
D_b(f_\mu(X_i),f_i^\star)
\ge\frac{|J|}{n}\left(\log2-\frac12\right).
\end{aligned}
\]
For $i\notin J$, we have
$f_\mu(X_i)\ge f_i^\star/2\ge2b$, and hence
\[
\frac{\phi(X_i-\theta)}
{n\max\{f_\mu(X_i),b\}}
\le\frac{2\phi(X_i-\theta)}{nf_i^\star}.
\]
For $i\in J$, the corresponding term is bounded by
$B/(nb)=4$. Therefore, using
\eqref{eq:fp-optimality},
\[
\begin{aligned}
\Psi_b[\mu](\theta)
\le\frac2n\sum_{i\notin J}
\frac{\phi(X_i-\theta)}{f_i^\star}
+4|J|
\le2+4|J|
\le2+\frac{4n}{\log2-1/2}
\bigl(\widehat L_n-\mathcal L_b(\mu)\bigr).
\end{aligned}
\]
Taking the supremum over $\theta\in\mathbb R^d$
proves \eqref{eq:fp-score-gap}.

We next prove \eqref{eq:fp-absolute-score}.
For $u\ge1/2$, Taylor's formula gives
\[
u-1-\log u
\ge\frac{(u-1)^2}{2\max\{1,u\}^2},
\quad
|u^{-1}-1|\le3\sqrt{u-1-\log u}.
\]
For $i\notin J$, we have
$f_\mu(X_i)\ge2b$, so the denominator is not
affected by the continuation. Thus,
\[
\begin{aligned}
&\frac1n\sum_{i\notin J}\phi(X_i-\theta)
\left|
\frac1{\max\{f_\mu(X_i),b\}}
-\frac1{f_i^\star}
\right|
=\frac1n\sum_{i\notin J}
\frac{\phi(X_i-\theta)}{f_i^\star}
\left|\frac{f_i^\star}{f_\mu(X_i)}-1\right|\\
&\le\frac3n\sum_{i\notin J}
\frac{\phi(X_i-\theta)}{f_i^\star}
\sqrt{D_b(f_\mu(X_i),f_i^\star)}
\le3
\left(\frac1n\sum_{i\notin J}
\frac{\phi(X_i-\theta)}{f_i^\star}\right)^{1/2}
\left(\frac1n\sum_{i\notin J}
\frac{\phi(X_i-\theta)}{f_i^\star}
D_b(f_\mu(X_i),f_i^\star)\right)^{1/2}\\
&\quad\le3\sqrt{
n\bigl(\widehat L_n-\mathcal L_b(\mu)\bigr)},
\end{aligned}
\]
where the last inequality follows from
\eqref{eq:fp-optimality},
$\phi(X_i-\theta)/(nf_i^\star)\le1$,
and \eqref{eq:fp-gap-decomposition}.

For $i\in J$, we have
$\max\{f_\mu(X_i),b\}\le f_i^\star/2$.
Consequently,
\[
\begin{aligned}
\frac1n\sum_{i\in J}\phi(X_i-\theta)
\left|
\frac1{\max\{f_\mu(X_i),b\}}
-\frac1{f_i^\star}
\right| \le\frac1n\sum_{i\in J}
\frac{\phi(X_i-\theta)}{\max\{f_\mu(X_i),b\}}
\le4|J|\\
\le\frac{4n}{\log2-1/2}
\bigl(\widehat L_n-\mathcal L_b(\mu)\bigr).
\end{aligned}
\]
Combining the two estimates and taking the supremum
over $\theta\in\mathbb R^d$ proves
\eqref{eq:fp-absolute-score}.
\end{proof}

\subsection{Positive Gaussian sums and the fixed-point iteration}

\begin{lemma}[Maxima and derivative bounds for Gaussian sums]
\label{lem:fp-gaussian}
Let
\[
P(\theta)=\sum_{i=1}^n v_i\phi(X_i-\theta),
\quad v_i\ge0,
\]
and suppose that $P\not\equiv0$. Set
$M=\sup_{\theta\in\mathbb R^d}P(\theta)$.
Then every maximizer $z$ of $P$ belongs to
$\operatorname{conv}\{X_1,\ldots,X_n\}$. Moreover,
\begin{equation}\label{eq:fp-gaussian-max}
\begin{aligned}
P(z+h)&\ge Me^{-\|h\|^2/2},
&&h\in\mathbb R^d,\\
|P(x+iy)|&\le Me^{\|y\|^2/2},
&&x,y\in\mathbb R^d.
\end{aligned}
\end{equation}
The derivatives of $P$ satisfy
\begin{equation}\label{eq:fp-derivatives}
\begin{aligned}
\sup_{\theta\in\mathbb R^d}\|\nabla P(\theta)\|
\le2M, \ \ 
\sup_{\theta\in\mathbb R^d}\|\nabla^2P(\theta)\|_{\rm op}
\le4M, \ \ 
\sup_{\theta\in\mathbb R^d}\|D^3P(\theta)\|_{\rm op}
\le48M.
\end{aligned}
\end{equation}
In addition,
\[
\sup_{\substack{z\in\mathbb C^d\\
\|\operatorname{Im}z\|\le1}}
\|\nabla^2P(z)\|_{\rm op}\le32M.
\]
If the coefficients $v_i$ take values in $\mathbb{R}$,
the same derivative bounds hold with $M$ replaced by
\[
\sup_{\theta\in\mathbb R^d}
\sum_{i=1}^n|v_i|\phi(X_i-\theta).
\]
\end{lemma}

\begin{proof}
Since $P$ is positive and tends to zero at infinity,
it attains its maximum at some $z\in\mathbb R^d$.
Differentiating $P$ gives
\[
\nabla P(z)=\sum_{i=1}^n
v_i\phi(X_i-z)(X_i-z)=0.
\]
Since $P(z)=M$, we obtain
\[
z=\sum_{i=1}^n\frac{v_i\phi(X_i-z)}{M}X_i,
\quad
\sum_{i=1}^n\frac{v_i\phi(X_i-z)}{M}=1.
\]
Thus $z\in\operatorname{conv}\{X_1,\ldots,X_n\}$.
Moreover, for every $h\in\mathbb R^d$,
\[
\begin{aligned}
&P(z+h)
=Me^{-\|h\|^2/2}
\sum_{i=1}^n\frac{v_i\phi(X_i-z)}{M}
e^{(X_i-z)\cdot h}\\
&\ge Me^{-\|h\|^2/2}
\exp\left(
\sum_{i=1}^n
\frac{v_i\phi(X_i-z)}{M}(X_i-z)\cdot h
\right)
=Me^{-\|h\|^2/2},
\end{aligned}
\]
where the inequality follows from Jensen's inequality.

For $x,y\in\mathbb R^d$, we have
$
\phi(X_i-x-iy)
=\phi(X_i-x)
e^{\|y\|^2/2+i(X_i-x)\cdot y}.
$
Consequently,
$
|P(x+iy)|
\le e^{\|y\|^2/2}P(x)
\le Me^{\|y\|^2/2},
$
which proves \eqref{eq:fp-gaussian-max}.

We next establish the derivative bounds. Fix
$x\in\mathbb R^d$ and a real unit vector $u$, and consider
the entire function
\[
g(\zeta)=P(x+\zeta u),
\quad \zeta\in\mathbb C.
\]
For $|\zeta|\le1$, the preceding bound gives
$
|g(\zeta)|\le Me^{1/2}.
$
By Cauchy's integral formula,
\[
|g^{(k)}(0)|\le k!e^{1/2}M,
\quad k=1,2,3.
\]
Since
$
g^{(k)}(0)=D^kP(x)[u,\ldots,u],
$
we obtain
\[
\|\nabla P(x)\|\le e^{1/2}M<2M,
\quad
\|\nabla^2P(x)\|_{\rm op}
\le2e^{1/2}M<4M.
\]
For the third derivative, the polarization identity
for symmetric trilinear forms gives
\[
\begin{aligned}
\|D^3P(x)\|_{\rm op}
\le\frac{3^3}{3!}
\sup_{\|u\|=1}|D^3P(x)[u,u,u]|
\le27e^{1/2}M<48M.
\end{aligned}
\]
This proves \eqref{eq:fp-derivatives}.

Finally, let $z\in\mathbb C^d$ satisfy
$\|\operatorname{Im}z\|\le1$. For every complex unit
vector $u$ and $|\zeta|\le1$,
$
\|\operatorname{Im}(z+\zeta u)\|\le2,
$
and therefore
$
|P(z+\zeta u)|\le e^2M.
$
Applying Cauchy's formula and the polarization identity
for symmetric bilinear forms yields
\[
\|\nabla^2P(z)\|_{\rm op}
\le4e^2M<32M.
\]

If the coefficients are signed, the complex bound becomes
\[
|P(x+iy)|
\le e^{\|y\|^2/2}
\sum_{i=1}^n|v_i|\phi(X_i-x).
\]
The same applications of Cauchy's formula and polarization
give the stated derivative bounds with $M$ replaced by
$\sup_x\sum_i|v_i|\phi(X_i-x)$.
\end{proof}

\begin{lemma}[Step size and contraction]
\label{lem:fp-inner}
Suppose that all fitted-value approximations used in the
$k$-th proximal step satisfy \eqref{eq:fc-evaluation-error}.
Then the step size $\tau_k$ satisfies
\begin{equation}\label{eq:fp-step-bounds}
\begin{aligned}
\frac1{320n}\le\tau_k\le\frac1{64},
\quad
\|\Psi_b[\mu_{a^k}]\|_\infty
\le\frac1{64\tau_k}
\le3\|\Psi_b[\mu_{a^k}]\|_\infty.
\end{aligned}
\end{equation}

For $a=a^k$, $p=p^k$, and $\tau=\tau_k$, define
\[
x_i(y)=a_i+\frac{\tau y_i}{np_i},
\quad
T_i(y)=\frac{p_i}{\max\{F_i(x(y)),b\}},
\quad 1\le i\le n.
\]
Then $T$ is a contraction on $[1/2,2]^n$ in the
Euclidean norm, with contraction factor at most $1/8$.
Moreover, both $T$ and its numerical approximation
\[
\widehat T_i(y)
=\frac{p_i}{\max\{\widehat F_i(x(y)),b\}}
\]
map $[1/2,2]^n$ into $[3/4,4/3]^n$.

When the stopping condition \eqref{eq:fc-stop} is satisfied,
the accepted iterate $a^{k+1}$ obeys
\begin{equation}\label{eq:fp-residual}
\frac Bn\sum_{i=1}^n
\left|
\frac{n(a_i^{k+1}-a_i^k)}{\tau_k}
-\frac1{\max\{F_i(a^{k+1}),b\}}
\right|
\le\eta_\varepsilon.
\end{equation}
The number of fitted-value evaluations required in each
proximal step is bounded as in
Theorem~\ref{thm:fc-iterations}.
\end{lemma}

\begin{proof}
Let
$\widehat\Psi(\theta)=n^{-1}\sum_i\phi(X_i-\theta)/p_i$.
By Lemma~\ref{lem:fp-gaussian}, its maximum is attained
at some $z\in\operatorname{conv}\{X_1,\ldots,X_n\}$.
Recall that the grid $\mathcal G$ in
\eqref{eq:fc-score-grid} has mesh width
$2R/\lceil2R\sqrt d\rceil$. Since
$\operatorname{conv}\{X_1,\ldots,X_n\}\subset K$,
there exists $\theta\in\mathcal G$ such that
\[
\|z-\theta\|
\le\frac{R\sqrt d}{\lceil2R\sqrt d\rceil}
\le\frac12.
\]
Applying \eqref{eq:fp-gaussian-max}, we obtain
$
\widehat\Psi(\theta)
\ge\widehat\Psi(z)e^{-\|z-\theta\|^2/2}
\ge e^{-1/8}\|\widehat\Psi\|_\infty,
$
and hence
$\|\widehat\Psi\|_\infty
\le e^{1/8}\max_{\mathcal G}\widehat\Psi$.

On the other hand, \eqref{eq:fc-evaluation-error}
and $p_i=\max\{\widehat F_i(a),b\}$ give
\[
\left|
\frac1{p_i}
-\frac1{\max\{F_i(a),b\}}
\right|
\le\frac{\delta_\varepsilon}{b^2}.
\]
Since $\phi\le B$, it follows that
$\|\widehat\Psi-\Psi_b[\mu_a]\|_\infty
\le B\delta_\varepsilon/b^2$.
Combining these estimates with \eqref{eq:fc-step}
gives
\[
\begin{aligned}
\|\Psi_b[\mu_a]\|_\infty
\le\|\widehat\Psi\|_\infty
+\frac{B\delta_\varepsilon}{b^2}
\le e^{1/8}\max_{\mathcal G}\widehat\Psi
+\frac{B\delta_\varepsilon}{b^2}
\le\frac1{64\tau}.
\end{aligned}
\]
Furthermore, $\widehat\Psi\le B/b=4n$ and
$\|\Psi_b[\mu_a]\|_\infty\ge1$ by
\eqref{eq:fp-score-gap}. Using the definition of
$\delta_\varepsilon$ in \eqref{eq:fc-tolerances},
we obtain
\[
\frac1{64\tau}\le5n,
\  \text{and }
\frac1{64\tau}
\le3\|\Psi_b[\mu_a]\|_\infty.
\]
Together with $\tau\le1/64$, these inequalities
prove \eqref{eq:fp-step-bounds}.

We next prove the contraction estimate. For
$y\in[1/2,2]^n$, we have
\begin{equation}\label{eq:fp-trial-bound}
0\le V_{x(y)}-V_a
=\frac{\tau}{n}\sum_{i=1}^n
\frac{y_i\phi(X_i-\cdot)}{p_i}
\le2\tau\|\widehat\Psi\|_\infty
\le\frac1{32}.
\end{equation}
Also, 
\begin{equation}\label{eq:fp-density-ratio}
\frac{d\mu_{x(y)}}{d\mu_a}(\theta)
=\frac{e^{V_{x(y)}(\theta)-V_a(\theta)}}
{\int_K e^{V_{x(y)}-V_a}\,d\mu_a},
 \implies e^{-1/32}\le
\frac{d\mu_{x(y)}}{d\mu_a}
\le e^{1/32}.
\end{equation}
Integrating $\phi(X_i-\theta)$ gives
\[
e^{-1/32}\le
\frac{F_i(x(y))}{F_i(a)}
\le e^{1/32},
\]
and the same inequalities hold for
$(\max\{F_i(x(y)),b\}/\max\{F_i(a),b\})$.

For any fitted-value approximation satisfying
$|\widehat F_i(s)-F_i(s)|\le\delta_\varepsilon$,
we have
\[
\frac1{1+\delta_\varepsilon/b}
\le
\frac{\max\{\widehat F_i(s),b\}}
{\max\{F_i(s),b\}}
\le1+\frac{\delta_\varepsilon}{b},
\]
since both denominators are bounded below by $b$.
Applying this inequality at $s=a$ and $s=x(y)$,
and using \eqref{eq:fp-density-ratio}, gives
\[
(1+\delta_\varepsilon/b)^{-2}e^{-1/32}
\le
\frac{p_i}{\max\{\widehat F_i(x(y)),b\}}
\le
(1+\delta_\varepsilon/b)^2e^{1/32}.
\]
Since $\delta_\varepsilon/b<1/8$, both endpoints
belong to $[3/4,4/3]$. The corresponding bounds
for the exact fitted values follow by the same argument.
Thus both $T$ and $\widehat T$ map $[1/2,2]^n$
into $[3/4,4/3]^n$.

We now estimate the derivative of $T$.
Differentiating \eqref{eq:fc-family} gives
\begin{equation}\label{eq:fp-covariance}
\frac{\partial F_i}{\partial a_j}(a)
=\operatorname{Cov}_{\mu_a}
\bigl(\phi(X_i-\Theta),\phi(X_j-\Theta)\bigr).
\end{equation}
For $v\in\mathbb R^n$, Cauchy--Schwarz gives
\[
\left(\sum_i v_i
\frac{\phi(X_i-\theta)}{p_i}\right)^2
\le
\left(\sum_i v_i^2
\frac{\phi(X_i-\theta)}{p_i}\right)
\left(\sum_i\frac{\phi(X_i-\theta)}{p_i}\right).
\]
Integrating against $\mu_{x(y)}$, and using
$F_i(x(y))/p_i\le2$, we obtain
\[
\begin{aligned}
\frac1n\mathbb E_{\mu_{x(y)}}
\left(\sum_i v_i
\frac{\phi(X_i-\Theta)}{p_i}\right)^2
\le
\|\widehat\Psi\|_\infty
\sum_i v_i^2\frac{F_i(x(y))}{p_i}
\le\frac{2}{64\tau}\|v\|^2.
\end{aligned}
\]
Consequently, the covariance matrix of
$(\phi(X_i-\Theta)/p_i)_i$, divided by $n$,
has operator norm at most $2/(64\tau)$.

At points where $F_i(x(y))\ne b$, differentiation
of $T$ gives
\[
DT(y)=-\frac{\tau}{n}
\operatorname{diag}\left(
\frac{p_i^2\mathbf1_{\{F_i(x(y))>b\}}}
{\max\{F_i(x(y)),b\}^2}\right)
\operatorname{Cov}_{\mu_{x(y)}}
\left(\left(\frac{\phi(X_i-\Theta)}{p_i}\right)_i\right).
\]
The diagonal matrix has operator norm at most $4$,
since $T_i(y)\le4/3$. Hence
$\|DT(y)\|_{\rm op}\le1/8$.
Since $T$ is Lipschitz, integration along line
segments gives the same bound across points
where $F_i(x(y))=b$. This proves that $T$ is
a contraction on $[1/2,2]^n$.

Let $y^\star$ denote its unique fixed point.
By the preceding bounds,
$\|T(y)-\widehat T(y)\|_\infty
\le2\delta_\varepsilon/b$.
Starting from $y^0=\mathbf1$, the contraction
estimate gives
\[
\begin{aligned}
\frac{\|y^r-y^\star\|}{\sqrt n}
&\le8^{-r}+\frac{16\delta_\varepsilon}{7b},\\
\frac{\|y^r-\widehat T(y^r)\|}{\sqrt n}
&\le\frac98\,8^{-r}
+\frac{32\delta_\varepsilon}{7b}.
\end{aligned}
\]
Since $p_i\ge b$, the left-hand side of
\eqref{eq:fc-stop} is bounded by
$(B/b)\|y^r-\widehat T(y^r)\|/\sqrt n$.
For
\[
r\ge\left\lceil
\frac{\log(8B/(b\eta_\varepsilon))}{\log8}
\right\rceil,
\]
this is at most
$9\eta_\varepsilon/64+\eta_\varepsilon/112
<\eta_\varepsilon/2$.
Thus the stopping condition is satisfied.

At acceptance,
$n(a_i^{k+1}-a_i^k)/\tau_k=y_i^r/p_i$.
By \eqref{eq:fc-evaluation-error},
\[
\left|
\frac1{\max\{\widehat F_i(a^{k+1}),b\}}
-\frac1{\max\{F_i(a^{k+1}),b\}}
\right|
\le\frac{\delta_\varepsilon}{b^2}.
\]
The resulting additional error in
\eqref{eq:fc-stop} is at most
$B\delta_\varepsilon/b^2=\eta_\varepsilon/512$,
which proves \eqref{eq:fp-residual}.
Including the evaluation at $a^k$ gives the
stated bound on the number of fitted-value evaluations.
\end{proof}

\subsection{Entropy estimates and iteration counts}

\begin{lemma}[Inexact proximal inequalities]
\label{lem:fp-entropy}
Suppose that $a^+$ is obtained from $a$ by an accepted
proximal step of size $\tau>0$ satisfying
\eqref{eq:fp-residual}. Let $\mu=\mu_a$ and
$\mu^+=\mu_{a^+}$. Then, for every $\nu\in\mathcal P(K)$
with $\operatorname{KL}(\nu\|\mu)<\infty$, we have
\begin{align}
\tau\bigl(\mathcal L_b(\nu)-\mathcal L_b(\mu^+)\bigr)
&\le
\operatorname{KL}(\nu\|\mu)
-\operatorname{KL}(\nu\|\mu^+)
-\operatorname{KL}(\mu^+\|\mu)
+2\tau\eta_\varepsilon.
\label{eq:fp-entropy}
\end{align}
Moreover, the accepted step satisfies
\begin{align}
&\mathcal L_b(\mu^+)
\ge\mathcal L_b(\mu)-\frac{\tau\eta_\varepsilon^2}{4},
\label{eq:fp-monotonicity} \\
&\log Z(a^+)-\log Z(a)
\le(1+\eta_\varepsilon)\tau.
\label{eq:fp-partition-step}
\end{align}
\end{lemma}

\begin{proof}
By \eqref{eq:fp-residual}, we may write
\[
\frac{a_i^+-a_i}{\tau}
=\frac1{n\max\{F_i(a^+),b\}}+r_i,
\quad
B\sum_{i=1}^n|r_i|\le\eta_\varepsilon.
\]
Setting $s(\theta)=\sum_i r_i\phi(X_i-\theta)$,
we have $\|s\|_\infty\le\eta_\varepsilon$ and
\[
\begin{aligned}
\log\frac{d\mu^+}{d\mu}
=V_{a^+}-V_a-\log\frac{Z(a^+)}{Z(a)}
=\tau\{\Psi_b[\mu^+]+s\}
-\log\frac{Z(a^+)}{Z(a)}.
\end{aligned}
\]
For any $\nu$ with finite
$\operatorname{KL}(\nu\|\mu)$, concavity of
$\mathcal L_b$ gives
\[
\mathcal L_b(\nu)-\mathcal L_b(\mu^+)
\le\int_K\Psi_b[\mu^+]\,d(\nu-\mu^+).
\]
Moreover,
\[
\begin{aligned}
\int_K\log\frac{d\mu^+}{d\mu}\,d(\nu-\mu^+)
=\operatorname{KL}(\nu\|\mu)
-\operatorname{KL}(\nu\|\mu^+)
-\operatorname{KL}(\mu^+\|\mu).
\end{aligned}
\]
Combining these identities, and using
$\|s\|_\infty\le\eta_\varepsilon$, yields
\[
\begin{aligned}
&\tau\{\mathcal L_b(\nu)-\mathcal L_b(\mu^+)\}
\quad\le\operatorname{KL}(\nu\|\mu)
-\operatorname{KL}(\nu\|\mu^+)
-\operatorname{KL}(\mu^+\|\mu)
-\tau\int_K s\,d(\nu-\mu^+)\\
&\quad\le\operatorname{KL}(\nu\|\mu)
-\operatorname{KL}(\nu\|\mu^+)
-\operatorname{KL}(\mu^+\|\mu)
+2\tau\eta_\varepsilon.
\end{aligned}
\]
This proves \eqref{eq:fp-entropy}.

Taking $\nu=\mu$ in the same argument, we obtain
\[
\begin{aligned}
\mathcal L_b(\mu)-\mathcal L_b(\mu^+)
\le-\frac1\tau
\left\{\operatorname{KL}(\mu\|\mu^+)
+\operatorname{KL}(\mu^+\|\mu)\right\}
+\eta_\varepsilon\|\mu-\mu^+\|_{TV} 
\\
\le-\frac{\|\mu-\mu^+\|_{TV}^2}{\tau}
+\eta_\varepsilon\|\mu-\mu^+\|_{TV}
\le\frac{\tau\eta_\varepsilon^2}{4},
\end{aligned}
\]
where the second inequality follows from Pinsker's
inequality. This proves \eqref{eq:fp-monotonicity}.

Finally, integrating the identity for $\log(d\mu^{+}/d\mu) $against
$\mu^+$ gives
\[
\begin{aligned}
\log Z(a^+)-\log Z(a)
&=\tau\int_K\{\Psi_b[\mu^+]+s\}\,d\mu^+
-\operatorname{KL}(\mu^+\|\mu)\\
&\le\tau(1+\eta_\varepsilon),
\end{aligned}
\]
since
\[
\int_K\Psi_b[\mu^+]\,d\mu^+
=\frac1n\sum_{i=1}^n
\frac{f_{\mu^+}(X_i)}
{\max\{f_{\mu^+}(X_i),b\}}
\le1.
\]
This establishes \eqref{eq:fp-partition-step}.
\end{proof}

We next construct a comparison measure to which
\eqref{eq:fp-entropy} can be applied. Let $\mu^\star$
be an NPMLE and, for $0<r\le1$, define
$\mu_r^\star=\mu^\star*
\operatorname{Unif}(B(0,r))$.
Recall that $\mu^\star$ is supported on
$\operatorname{conv}\{X_1,\ldots,X_n\}$, which lies at a distance of at least one from the boundary of $K$. Consequently,
$\operatorname{supp}(\mu_r^\star)\subset K$.

The density of $\mu_r^\star$ is bounded above by
$(v_dr^d)^{-1}$. Since $\mu_0$ is uniform on $K$, we have 
\[
\begin{aligned}
\operatorname{KL}(\mu_r^\star\|\mu_0)
=\int_K\log\left(
(2R)^d\frac{d\mu_r^\star}{d\theta}
\right)\,d\mu_r^\star
\le\log\frac{(2R)^d}{v_dr^d}.
\end{aligned}
\]
Moreover, for $U\sim\operatorname{Unif}(B(0,r))$,
Jensen's inequality gives
\[
\begin{aligned}
\mathbb E\,\phi(X_i-\theta-U)
&=\phi(X_i-\theta)
\mathbb E\exp\left\{
(X_i-\theta)\cdot U-\frac{\|U\|^2}{2}
\right\}\\
&\ge\phi(X_i-\theta)
\exp\left\{-\frac{\mathbb E\|U\|^2}{2}\right\}
\ge e^{-r^2/2}\phi(X_i-\theta).
\end{aligned}
\]
Integrating against $\mu^\star$ gives
$f_{\mu_r^\star}(X_i)
\ge e^{-r^2/2}f_i^\star$.
Since $h_b(u)\ge\log u$, we obtain
\begin{equation}\label{eq:fp-smoothing}
\begin{gathered}
\operatorname{supp}(\mu_r^\star)\subset K,\quad
\operatorname{KL}(\mu_r^\star\|\mu_0)
\le\log\frac{(2R)^d}{v_d}+d\log(1/r),\\
\widehat L_n-\mathcal L_b(\mu_r^\star)
\le\widehat L_n-L_n(f_{\mu_r^\star})
\le\frac{r^2}{2}.
\end{gathered}
\end{equation}

Similarly, since $B(X_i,1)\subset K$, we have 
\[
\begin{aligned}
f_{\mu_0}(X_i)
=\frac1{(2R)^d}
\int_K\phi(X_i-\theta)\,d\theta
\ge\frac{Be^{-1/2}v_d}{(2R)^d}.
\end{aligned}
\]
Using $\widehat L_n\le\log B$ and
$\mathcal L_b(\mu_0)\ge L_n(f_{\mu_0})$, we obtain
\begin{equation}\label{eq:fp-initial-gap}
\widehat L_n-\mathcal L_b(\mu_0)
\le\log\frac{(2R)^d}{v_d}+\frac12.
\end{equation}

\subsection{Proof of Theorem~\ref{thm:fc-iterations}}
\begin{proof}[\textbf{Proof of Theorem~\ref{thm:fc-iterations}}]
We first assume that all fitted-value evaluations
satisfy \eqref{eq:fc-evaluation-error}. By
\eqref{eq:fp-step-bounds}, the step sizes satisfy
$1/(320n)\le\tau_k\le1/64$. The algorithm therefore
terminates after finitely many proximal steps, with
\[
\frac{32H_\varepsilon}{\varepsilon}
\le t_N
<\frac{32H_\varepsilon}{\varepsilon}+\frac1{64}.
\]

We first establish the accuracy of the output.
Applying \eqref{eq:fp-entropy} at each step with
$\nu=\mu_{\sqrt\varepsilon/4}^\star$, and summing
over $0\le k<N$, gives
\[
\begin{aligned}
&\sum_{k<N}\tau_k
\{\mathcal L_b(\mu_{\sqrt\varepsilon/4}^\star)
-\mathcal L_b(\mu_{a^{k+1}})\}
\le
\operatorname{KL}(\mu_{\sqrt\varepsilon/4}^\star\|\mu_0)
+2\eta_\varepsilon t_N.
\end{aligned}
\]
By \eqref{eq:fp-smoothing} and the definition of
$H_\varepsilon$, we have
\[
\widehat L_n-
\mathcal L_b(\mu_{\sqrt\varepsilon/4}^\star)
\le\frac{\varepsilon}{32},
\ \text{and } 
\operatorname{KL}(\mu_{\sqrt\varepsilon/4}^\star\|\mu_0)
\le H_\varepsilon.
\]
Consequently,
\[
\begin{aligned}
\frac1{t_N}\sum_{k<N}\tau_k
\{\widehat L_n-\mathcal L_b(\mu_{a^{k+1}})\}
\le\frac{\varepsilon}{32}
+\frac{H_\varepsilon}{t_N}
+2\eta_\varepsilon
\le\frac{\varepsilon}{16}+2\eta_\varepsilon.
\end{aligned}
\]
By \eqref{eq:fc-evaluation-error} and the
$b^{-1}$-Lipschitz continuity of $h_b$, the error
in evaluating $\mathcal L_b(\mu_{a^k})$ is at most
$\delta_\varepsilon/b$. Since $a^{\widehat k}$
maximizes the computed likelihood over the
completed iterates, we obtain
\[
\begin{aligned}
\widehat L_n-\mathcal L_b(\mu_{a^{\widehat k}})
\le\frac1{t_N}\sum_{k<N}\tau_k
\{\widehat L_n-\mathcal L_b(\mu_{a^{k+1}})\}
+\frac{2\delta_\varepsilon}{b}
\le\frac{\varepsilon}{16}
+2\eta_\varepsilon+\frac{2\delta_\varepsilon}{b}
<\frac{\varepsilon}{8}.
\end{aligned}
\]
Applying Lemma~\(\ref{lem:fp-equivalence}\) with \(\mu=\mu_{a^{\widehat k}}\), we obtain
\[
\begin{aligned}
\widehat L_n-L_n(f_{\mu_{\rm out}})
\le\left(1+\frac1{\log2-1/2}\right)
\{\widehat L_n-\mathcal L_b(\mu_{a^{\widehat k}})\}
<\frac{\varepsilon}{8}
\left(1+\frac1{\log2-1/2}\right)
<\varepsilon.
\end{aligned}
\]
This proves \eqref{eq:fc-accuracy}.
The output is absolutely continuous, since both
$\mu_{a^{\widehat k}}$ and the uniform measures
used in the correction are absolutely continuous.

We next bound the number of proximal steps.
For any \(1\le k\le N\), let \(t=t_k=\sum_{j=0}^{k-1}\tau_j\). Applying \eqref{eq:fp-entropy}
with $\nu=\mu_{(1+t)^{-1/2}}^\star$, and using
\eqref{eq:fp-smoothing}, gives
\begin{equation}\label{eq:fp-weighted-gap}
\begin{aligned}
\sum_{j<k}\tau_j
\{\widehat L_n-\mathcal L_b(\mu_{a^{j+1}})\}
\le
\log\frac{(2R)^d}{v_d}+\frac12
+\frac d2\log(1+t)+2\eta_\varepsilon t.
\end{aligned}
\end{equation}
We next derive a bound for the sum involving \(\mathcal L_b(\mu_{a^j})\), rather than \(\mathcal L_b(\mu_{a^{j+1}})\).
By \eqref{eq:fp-monotonicity}, we have
\[
\mathcal L_b(\mu_{a^{j+1}})
-\mathcal L_b(\mu_{a^j})
\ge-\frac{\tau_j\eta_\varepsilon^2}{4}.
\].
Summing over the indices for which
$\mathcal L_b(\mu_{a^{j+1}})<\mathcal L_b(\mu_{a^j})$,
we obtain
\[
\sum_{j<k}
\max\left\{
\mathcal L_b(\mu_{a^j})
-\mathcal L_b(\mu_{a^{j+1}}),0
\right\}
\le\frac{\eta_\varepsilon^2}{4}
\sum_{j<k}\tau_j
=\frac{\eta_\varepsilon^2t}{4}.
\] 
Since
$
\sum_{j<k}
\bigl(\mathcal L_b(\mu_{a^{j+1}})
-\mathcal L_b(\mu_{a^j})\bigr)
=\mathcal L_b(\mu_{a^k})-\mathcal L_b(\mu_0)
\le\widehat L_n-\mathcal L_b(\mu_0),
$
the preceding estimate gives
\[
\sum_{j<k}\max\left\{
\mathcal L_b(\mu_{a^{j+1}})
-\mathcal L_b(\mu_{a^j}),0
\right\}
\le\widehat L_n-\mathcal L_b(\mu_0)
+\frac{\eta_\varepsilon^2t}{4}.
\]
Using $\tau_j\le1/64$, we obtain
\[
\begin{aligned}
\sum_{j<k}\tau_j
\{\mathcal L_b(\mu_{a^{j+1}})
-\mathcal L_b(\mu_{a^j})\}
\le\frac1{64}
\left\{\widehat L_n-\mathcal L_b(\mu_0)
+\frac{\eta_\varepsilon^2t}{4}\right\}.
\end{aligned}
\]
Combining this with \eqref{eq:fp-weighted-gap}
and \eqref{eq:fp-initial-gap} yields
\begin{equation}\label{eq:fp-old-gap}
\begin{aligned}
&\sum_{j<k}\tau_j
\{\widehat L_n-\mathcal L_b(\mu_{a^j})\}\\
&\le
\log\frac{(2R)^d}{v_d}+\frac12
+\frac d2\log(1+t)+2\eta_\varepsilon t
+\frac1{64}
\left\{
\log\frac{(2R)^d}{v_d}+\frac12
+\frac{\eta_\varepsilon^2t}{4}
\right\}.
\end{aligned}
\end{equation}

By \eqref{eq:fp-step-bounds} and
\eqref{eq:fp-score-gap},
\[
\begin{aligned}
k
=\sum_{j<k}\tau_j\frac1{\tau_j}
\le192\sum_{j<k}\tau_j
\|\Psi_b[\mu_{a^j}]\|_\infty
\le384t+
\frac{768n}{\log2-1/2}
\sum_{j<k}\tau_j
\{\widehat L_n-\mathcal L_b(\mu_{a^j})\}.
\end{aligned}
\]
Substituting \eqref{eq:fp-old-gap}, and using
$n\eta_\varepsilon\le1/256$, gives
\[
k\le C\left[
t+n\left\{
\log\frac{(2R)^d}{v_d}
+\frac12+\frac d2\log(1+t)
\right\}\right].
\]
This proves the first bound in \eqref{eq:fc-counts}.
At termination, using the comparison measure
$\mu_{\sqrt\varepsilon/4}^\star$ constructed above gives
\[
N\le C H_\varepsilon(n+\varepsilon^{-1}).
\]

We next derive an alternative bound on the number of proximal steps using the increments of the coefficients $a_i^k$. Since $b\le p_i^k\le B$, we have
$B n^{-1}\sum_i(p_i^k)^{-1}\ge1$. Moreover,
\eqref{eq:fc-step} and $\phi\le B$ give
\[
\frac1{64\tau_k}
\le\frac{2B}{n}\sum_{i=1}^n\frac1{p_i^k}.
\]
By Lemma~\ref{lem:fp-inner}, every accepted fixed-point
iterate satisfies $y_i\ge1/2$. Therefore,
\begin{equation}\label{eq:fp-coefficient-increment}
\begin{aligned}
B\sum_{i=1}^n(a_i^{k+1}-a_i^k)
&=\frac{\tau_k B}{n}
\sum_{i=1}^n\frac{y_i}{p_i^k}
\ge\frac{\tau_k B}{2n}
\sum_{i=1}^n\frac1{p_i^k}
\ge\frac1{256}.
\end{aligned}
\end{equation}
On the other hand, $B(X_i,1)\subset K$ implies
$
\int_K\phi(X_i-\theta)\,d\theta
\ge Be^{-1/2}v_d.
$
Since $a_i\ge0$, it follows that
\[
Be^{-1/2}v_d\sum_{i=1}^n a_i
\le\int_K V_a(\theta)\,d\theta
\le(2R)^d\sup_K V_a.
\]
Summing \eqref{eq:fp-coefficient-increment}
over $0\le j<k$, and using $a^0=0$, gives
\[
\frac{k}{256}
\le B\sum_{i=1}^n a_i^k
\le\frac{e^{1/2}(2R)^d}{v_d}
\sup_K V_{a^k}.
\]
Applying Lemma~\ref{lem:fp-potential}, we obtain
\begin{equation}\label{eq:fp-budget-count}
k\le\frac{256e^{1/2}(2R)^d}{v_d}
\mathcal M(t_k).
\end{equation}
At $k=N$, this yields
$N\le C_dR^dH_\varepsilon/\varepsilon$.
Together with the preceding estimate
$N\le CH_\varepsilon(n+\varepsilon^{-1})$,
this proves \eqref{eq:fc-counts}.

We finally show that \(\tau_k>1/128\) after at most \(C_dn\log(enR)\) proximal iterations.
Choose $T_\star\ge1$ such that
\begin{equation}\label{eq:fp-good-time}
\frac{
\log((2R)^d/v_d)+1/2
+(d/2)\log(1+T_\star)}
{T_\star}
<
\frac{\log(4/3)-1/4}{8n}.
\end{equation}
Such a choice is possible with
$T_\star\le C_dn\log(enR)$.

Suppose the algorithm reaches time $T_\star$.
By \eqref{eq:fp-weighted-gap}, there is an
iterate before time $T_\star+1/64$ satisfying
\[
\widehat L_n-\mathcal L_b(\mu_{a^j})
<
\frac{\log(4/3)-1/4}{8n}
+2\eta_\varepsilon.
\]

Summing \eqref{eq:fp-monotonicity} from $j$ to $k-1$,
we obtain, for every $k\ge j$,
$
\widehat L_n-\mathcal L_b(\mu_{a^k})
\le
\widehat L_n-\mathcal L_b(\mu_{a^j})
+({\eta_\varepsilon^2}/{4})(t_k-t_j).
$

Since
$t_N<32H_\varepsilon/\varepsilon+1/64$,
every later iterate satisfies
\[
\begin{aligned}
\widehat L_n-\mathcal L_b(\mu_{a^k})
<
\frac{\log(4/3)-1/4}{8n}
+2\eta_\varepsilon
+\frac{\eta_\varepsilon^2}{4}
\left(\frac{32H_\varepsilon}{\varepsilon}
+\frac1{64}\right)
<\frac{\log(4/3)-1/4}{n}.
\end{aligned}
\]
By \eqref{eq:fp-gap-decomposition}, this implies
\[
D_b(F_i(a^k),f_i^\star)
<\log(4/3)-\frac14,
\quad 1\le i\le n.
\]
Since $f_i^\star\ge4b$, the expression for
$D_b$ in Lemma~\ref{lem:fp-equivalence} gives
\[
F_i(a^k)>\frac34f_i^\star>b,
\quad 1\le i\le n.
\]
Thus, by \eqref{eq:fp-optimality},
\[
\Psi_b[\mu_{a^k}](\theta)
=\frac1n\sum_{i=1}^n
\frac{\phi(X_i-\theta)}{F_i(a^k)}
\le\frac43,
\quad\theta\in\mathbb R^d.
\]
Using \eqref{eq:fc-evaluation-error}
and \eqref{eq:fc-step}, we obtain
\[
\frac1{64\tau_k}
\le\max\left\{
1,\,
e^{1/8}\left(
\frac43+\frac{B\delta_\varepsilon}{b^2}
\right)
+\frac{B\delta_\varepsilon}{b^2}
\right\}<2.
\]
Hence $\tau_k>1/128$ at every subsequent step.

Suppose first that the algorithm reaches time $T_\star$,
and let $k_0$ be the smallest index such that
$t_{k_0}\ge T_\star$. Since $\tau_k\le1/64$,
we have $t_{k_0}<T_\star+1/64$.
Applying the first inequality in \eqref{eq:fc-counts}
with $k=k_0$ gives $k_0\le C_dn\log(enR)$.
By the preceding argument, every step with
$k\ge k_0$ satisfies $\tau_k>1/128$.

If the algorithm terminates before $T_\star$,
the same inequality gives $N\le C_dn\log(enR)$.
Thus, in either case, we may choose an integer
$k_\star\le C_dn\log(enR)$, independent of
$\varepsilon$, such that every performed step
with $k\ge k_\star$ satisfies $\tau_k>1/128$.
Since $t_N<32H_\varepsilon/\varepsilon+1/64$,
there are at most $4096H_\varepsilon/\varepsilon+2$
additional steps.

Finally, Lemma~\ref{lem:fp-global-quadrature}
establishes \eqref{eq:fc-evaluation-error} by
induction over the proximal steps, using
\eqref{eq:fp-residual} for the preceding steps.
Hence the assumed evaluation accuracy holds
throughout Algorithm~\ref{alg:fc-global},
completing the proof.

\end{proof}

\subsection{Quadrature on the whole box}

\begin{lemma}[Size of the potential]
\label{lem:fp-potential}
Suppose that all completed proximal steps satisfy
\eqref{eq:fp-residual}. Then, for every completed
iterate $a^k$ with accumulated time $t_k$, we have
\[
\sup_{\theta\in\mathbb R^d}V_{a^k}(\theta)
\le\mathcal M(t_k).
\]
The same bound holds for every trial vector
$x(y)$ in the subsequent proximal step, with
$y\in[1/2,2]^n$.
\end{lemma}

\begin{proof}
Let $M=\sup_{\theta\in\mathbb R^d}V_a(\theta)$.
By Lemma~\ref{lem:fp-gaussian}, $V_a$ attains its
maximum at some
$z\in\operatorname{conv}\{X_1,\ldots,X_n\}$.
At this point, $\nabla V_a(z)=0$. Moreover,
\[
\begin{aligned}
\nabla^2V_a(\theta)
&=\sum_{i=1}^n a_i\phi(X_i-\theta)
\bigl\{(X_i-\theta)(X_i-\theta)^{\mathsf T}-I_d\bigr\}\\
&\succeq-V_a(\theta)I_d
\succeq-MI_d.
\end{aligned}
\]
Taylor's formula therefore gives
$V_a(z+h)\ge M-M\|h\|^2/2$.
In particular, for $\|h\|\le(1+M)^{-1/2}$,
we have $V_a(z+h)\ge M-1/2$.

Since
$\operatorname{dist}(\operatorname{conv}\{X_i\},\partial K)
\ge1$, the ball
$B(z,(1+M)^{-1/2})$ is contained in $K$.
Consequently,
\begin{equation}\label{eq:fp-normalizer}
\begin{aligned}
Z(a)
=\frac1{(2R)^d}\int_K e^{V_a(\theta)}\,d\theta
\ge\frac{e^{M-1/2}v_d}
{(2R)^d(1+M)^{d/2}}.
\end{aligned}
\end{equation}
The estimate also holds for $a=0$, since
$V_0=0$ and $Z(0)=1$.
Taking logarithms and using
$\log(1+M)\le M/(d+1)+\log(d+1)$,
we obtain
\[
\begin{aligned}
M
&\le\log Z(a)+\log\frac{(2R)^d}{v_d}
+\frac12+\frac d2\log(1+M)\\
&\le\log Z(a)+\log\frac{(2R)^d}{v_d}
+\frac12+\frac{dM}{2(d+1)}
+\frac d2\log(d+1).
\end{aligned}
\]
Rearranging gives
\[
M\le2\log Z(a)
+2\log\frac{(2R)^d}{v_d}
+1+d\log(d+1).
\]

For a completed iterate $a^k$, summing
\eqref{eq:fp-partition-step} and using $Z(0)=1$
gives $\log Z(a^k)\le(1+\eta_\varepsilon)t_k$.
Since $\eta_\varepsilon\le1/2$, it follows that
\[
\sup_{\mathbb R^d}V_{a^k}
\le3t_k+2\log\frac{(2R)^d}{v_d}
+1+d\log(d+1)
\le\mathcal M(t_k).
\]

For a trial vector $x(y)$ with
$y\in[1/2,2]^n$, \eqref{eq:fp-trial-bound}
gives $0\le V_{x(y)}-V_{a^k}\le1/32$.
Hence
$\log Z(x(y))\le\log Z(a^k)+1/32$,
and the preceding estimate yields
\[
\begin{aligned}
\sup_{\mathbb R^d}V_{x(y)}
&\le2(1+\eta_\varepsilon)t_k
+2\log\frac{(2R)^d}{v_d}
+1+d\log(d+1)+\frac1{16}\\
&\le\mathcal M(t_k).
\end{aligned}
\]
This completes the proof.
\end{proof}

\begin{lemma}[Analytic tensor quadrature]
\label{lem:fp-analytic}
Let $\mathcal E_\varrho\subset\mathbb C$ denote the
Bernstein ellipse with foci $\pm1$ and parameter
$\varrho>1$. Suppose that $g$ is analytic on a
neighborhood of $\mathcal E_\varrho^d$ and satisfies
$|g(z)|\le M$ for every $z\in\mathcal E_\varrho^d$.
Then the tensor Gauss--Legendre rule $Q_q$ of order
$q$ in each coordinate, normalized with respect to
$\operatorname{Unif}([-1,1]^d)$, satisfies
\begin{equation}\label{eq:fp-analytic-error}
\left|
Q_q(g)-\int_{[-1,1]^d}
g(x)\,d\operatorname{Unif}([-1,1]^d)(x)
\right|
\le
2^{d+1}dM(1-\varrho^{-1})^{-d}\varrho^{-2q}.
\end{equation}
\end{lemma}

\begin{proof}
By the standard Chebyshev coefficient estimate
\cite[Theorem~8.1]{trefethen2019},
the tensor Chebyshev coefficients of $g$ satisfy
$|c_\alpha|\le2^dM\varrho^{-|\alpha|_1}$.
Let $p$ be the truncation retaining the terms
with $\alpha_j<2q$ for every $j$. Then
\[
\begin{aligned}
\|g-p\|_{L^\infty([-1,1]^d)}
&\le2^dM
\sum_{\max_j\alpha_j\ge2q}\varrho^{-|\alpha|_1}\\
&\le2^ddM(1-\varrho^{-1})^{-d}\varrho^{-2q}.
\end{aligned}
\]
The tensor Gauss--Legendre rule integrates $p$
exactly. Since its weights are positive and
normalized to sum to one, we have
\[
\left|Q_q(g)-\int g\,d\operatorname{Unif}([-1,1]^d)\right|
\le2\|g-p\|_\infty,
\]
which proves \eqref{eq:fp-analytic-error}.
\end{proof}

\begin{lemma}[Construction of positive nodes and weights]
\label{lem:fp-rounded}
Let $(z_j,\omega_j)_{j=1}^q$ denote the nodes and
weights of the $q$-point Gauss--Legendre quadrature
rule on $[-1,1]$. For every $q\ge1$ and
$0<\nu\le(8q^6)^{-1}$, one can compute approximations
$(\widetilde z_j,\widetilde\omega_j)_{j=1}^q$ in
$O(q^2\log(e/\nu))$ arithmetic operations such that
\[
\begin{gathered}
|\widetilde z_j-z_j|\le\nu,\quad
\widetilde\omega_j>0,\quad
\sum_{j=1}^q\widetilde\omega_j=2,\\
\sum_{j=1}^q|\widetilde\omega_j-\omega_j|
\le64q^7\nu.
\end{gathered}
\]
Let $Q_q$ and $\widetilde Q_q$ denote the normalized
tensor quadrature rules on $[-1,1]^d$ obtained
from the exact and approximate nodes and weights,
respectively. Then, for every
$g\in C^1([-1,1]^d)$ satisfying
$\|g\|_\infty\le M$ and
$\|\nabla g\|_\infty\le G$, we have
\begin{equation}\label{eq:fp-rounded-error}
|\widetilde Q_q(g)-Q_q(g)|
\le(32dq^7M+\sqrt d\,G)\nu.
\end{equation}
\end{lemma}

\begin{proof}
Let $P_q$ denote the Legendre polynomial of degree $q$.
We use the standard recurrence, zero-counting property,
and Gauss--Legendre weight formula; see
\cite[Sections~1.2--1.5]{2004orthogonal}.
In particular,
\[
P_0(x)=1,\quad P_1(x)=x,\quad
(j+1)P_{j+1}(x)=(2j+1)xP_j(x)-jP_{j-1}(x).
\]
The zeros of $P_q$ can be counted on any interval
by the sign variations of $(P_0,\ldots,P_q)$
at its endpoints. Each count requires $O(q)$
operations.

We first approximate the nodes by bisection.
Since $\|P_q\|_\infty\le1$, Markov's inequality gives
$\|P_q'\|_\infty\le q^2$ and
$\|P_q''\|_\infty\le q^4$.
The exact Gauss--Legendre weights satisfy
\[
\omega_j=\frac{2}{(1-z_j^2)P_q'(z_j)^2}>0,
\quad
\sum_{j=1}^q\omega_j=2.
\]
In particular, $\omega_j\le2$ implies
$|P_q'(z_j)|\ge1$. For consecutive zeros
$z_j<z_{j+1}$, Rolle's theorem and the
second derivative bound give
\[
1\le |P_q'(z_j)|
\le q^4(z_{j+1}-z_j).
\]
Similarly, $P_q(\pm1)=(\pm1)^q$ gives
$1-|z_j|\ge q^{-2}$. Thus distinct zeros
are separated by at least $q^{-4}$ and lie
at distance at least $q^{-2}$ from the endpoints.

Starting from $[-1,1]$, we bisect intervals
containing zeros, discarding those containing none
and recording exact midpoint zeros separately.
Since $\nu\le(8q^6)^{-1}$, intervals of length
at most $\nu$ contain at most one zero. Their
midpoints give approximations
$|\widetilde z_j-z_j|\le\nu$ in $(-1,1)$.
At each depth there are at most $q$ intervals
containing zeros, and each zero count costs
$O(q)$ operations. Thus all nodes are computed
in $O(q^2\log(e/\nu))$ operations.

For the weights, define
$D_q(x)=(1-x^2)P_q'(x)^2$. The preceding
derivative bounds give
\[
|D_q'(x)|\le2q^4+2q^6\le4q^6,
\quad
D_q(z_j)=\frac2{\omega_j}\ge1.
\]
Consequently, for every $x$ between $z_j$
and $\widetilde z_j$,
\[
D_q(x)\ge1-4q^6\nu\ge\frac12,
\quad
\left|\left(\frac2{D_q}\right)'(x)\right|
\le32q^6.
\]
Set
\[
v_j=\frac2{D_q(\widetilde z_j)},
\quad
\widetilde\omega_j
=\frac{2v_j}{\sum_{\ell=1}^q v_\ell}.
\]
Then $\widetilde\omega_j>0$,
$\sum_j\widetilde\omega_j=2$, and
$|v_j-\omega_j|\le32q^6\nu$.
Writing $S=\sum_jv_j$, we obtain
\[
\begin{aligned}
\sum_{j=1}^q|\widetilde\omega_j-\omega_j|
\le |2-S|+\sum_{j=1}^q|v_j-\omega_j|
\le2\sum_{j=1}^q|v_j-\omega_j|
\le64q^7\nu.
\end{aligned}
\]
The values $P_q'(\widetilde z_j)$ are computed
by differentiating the three-term recurrence,
at a cost of $O(q)$ operations per node.
Thus computing and normalizing the weights
requires $O(q^2)$ additional operations.

We next establish \eqref{eq:fp-rounded-error}.
Let $p_j=\omega_j/2$ and
$\widetilde p_j=\widetilde\omega_j/2$.
Both vectors have positive entries summing to
one, and
$\|p-\widetilde p\|_{TV}\le32q^7\nu$.
Telescoping their tensor products gives
\[
\begin{aligned}
p^{\otimes d}-\widetilde p^{\otimes d}
=\sum_{j=1}^d
\widetilde p^{\otimes(j-1)}
\otimes(p-\widetilde p)
\otimes p^{\otimes(d-j)}.
\end{aligned}
\]
Hence
\[
\|p^{\otimes d}-\widetilde p^{\otimes d}\|_{TV}
\le d\|p-\widetilde p\|_{TV}
\le32dq^7\nu.
\]
For $u\in\{1,\ldots,q\}^d$, write
$z_u=(z_{u_1},\ldots,z_{u_d})$ and
$\widetilde z_u
=(\widetilde z_{u_1},\ldots,\widetilde z_{u_d})$.
Then $\|\widetilde z_u-z_u\|\le\sqrt d\,\nu$.
Therefore,
\[
\begin{aligned}
|\widetilde Q_q(g)-Q_q(g)|
\le
M\|p^{\otimes d}-\widetilde p^{\otimes d}\|_{TV}
+
G\sum_u\widetilde p_u^{\otimes d}
\|\widetilde z_u-z_u\|
\le(32dq^7M+\sqrt d\,G)\nu,
\end{aligned}
\]
which proves \eqref{eq:fp-rounded-error}.

For quadrature on $[-L,L]^d$, apply
\eqref{eq:fp-rounded-error} to $g(x)=h(Lx)$,
using $\|\nabla g\|_\infty
\le L\|\nabla h\|_\infty$, and multiply
the resulting error by $(2L)^d$.
The corresponding estimate on $K$ follows
by the change of variables $\theta=X_1+Rx$.
\end{proof}

For the global quadrature rule of order $q_K$ in
\eqref{eq:fc-global-order}, let $\overline M$ denote
the corresponding bound on the potential.
We construct the one-dimensional Gauss--Legendre nodes
and weights using Lemma~\ref{lem:fp-rounded}, with
tolerance
\begin{equation}\label{eq:fp-global-tolerance}
\nu_K=\min\left\{
\frac1{8q_K^6},
\frac{\delta_\varepsilon e^{-1/2}v_d}
{8B(2R)^d(1+\overline M)^{d/2}
[32dq_K^7+R\sqrt d(2\overline M+1)]}
\right\}.
\end{equation}
The tensor product of these nodes and weights,
mapped to $K$, defines the quadrature rule in
\eqref{eq:fc-quadrature}. Its weights are positive
and sum exactly to one.

\begin{lemma}[Accuracy of global quadrature]
\label{lem:fp-global-quadrature}
The quadrature order \eqref{eq:fc-global-order}
and tolerance \eqref{eq:fp-global-tolerance}
ensure that every fitted-value evaluation in
Algorithm~\ref{alg:fc-global} satisfies
\eqref{eq:fc-evaluation-error}.
Each evaluation requires $O(ndq_K^d)$ arithmetic
operations, and the one-dimensional nodes and
weights can be constructed in
$O(q_K^2\log(e/\nu_K))$ operations.
\end{lemma}

\begin{proof}
We first establish the quadrature error for any
coefficient vector $a$ satisfying
$M=\sup_{\mathbb R^d}V_a\le\overline M$.
Set
$\alpha=[4R\sqrt{d(1+\overline M)}]^{-1}$.
On the product of Bernstein ellipses of parameter
$e^\alpha$, mapped to $K$, the imaginary part
satisfies
\[
\|y\|\le\frac1{2\sqrt{1+\overline M}}.
\]
By \eqref{eq:fp-gaussian-max}, we have
\[
\begin{aligned}
\operatorname{Re}V_a(x+iy)-M
&\le M(e^{\|y\|^2/2}-1)\le\frac14,\\
\frac{|\phi(X_i-x-iy)|}{B}
&\le e^{\|y\|^2/2}\le e^{1/8}.
\end{aligned}
\]
Consequently, the functions
$e^{V_a-M}$ and
$B^{-1}\phi(X_i-\cdot)e^{V_a-M}$
are analytic and bounded in absolute value by
two on these ellipses.

Since $1-e^{-\alpha}\ge\alpha/2$,
Lemma~\ref{lem:fp-analytic} bounds the quadrature
error for each of these functions by
$2^{2d+2}d\alpha^{-d}e^{-2\alpha q_K}$.
We claim that
\begin{equation}\label{eq:fp-global-analytic-size}
2^{2d+2}d\alpha^{-d}e^{-2\alpha q_K}
\le
\frac{\delta_\varepsilon}{8B}
\frac{e^{-1/2}v_d}
{(2R)^d(1+\overline M)^{d/2}}.
\end{equation}
Indeed, using $v_d\ge(2/\sqrt d)^d$ and the
definition of $\alpha$, it suffices that
\[
\begin{aligned}
2\alpha q_K
&\ge(4d+5)\log2+\log d+\frac12
+d\log d+2d\log R\\
&\quad+d\log(1+\overline M)
+\log(B/\delta_\varepsilon).
\end{aligned}
\]
This follows from \eqref{eq:fc-global-order},
which gives
\[
2\alpha q_K
\ge8d\log\frac{16dR(1+\overline M)B}
{\delta_\varepsilon}.
\]

We next bound the error in the computed nodes
and weights. On the real cube, both
$e^{V_a-M}$ and
$B^{-1}\phi(X_i-\cdot)e^{V_a-M}$
are bounded by one. Moreover,
Lemma~\ref{lem:fp-gaussian} gives
$\|\nabla V_a\|_\infty\le2M$, while
$\|\nabla\phi(X_i-\cdot)\|_\infty\le B$.
Their gradients are therefore bounded by
$2\overline M+1$.

After the change of variables
$\theta=X_1+Rx$, the corresponding functions
on $[-1,1]^d$ have gradients bounded by
$R(2\overline M+1)$. Applying
\eqref{eq:fp-rounded-error} gives
\[
\begin{aligned}
|\widetilde Q_{q_K}(g)-Q_{q_K}(g)|
&\le
[32dq_K^7+R\sqrt d(2\overline M+1)]\nu_K\\
&\le
\frac{\delta_\varepsilon}{8B}
\frac{e^{-1/2}v_d}
{(2R)^d(1+\overline M)^{d/2}}.
\end{aligned}
\]
Thus, by \eqref{eq:fp-global-analytic-size},
the combined analytic and construction errors
for each shifted integrand are at most
\[
\frac{\delta_\varepsilon}{4B}
\frac{e^{-1/2}v_d}
{(2R)^d(1+\overline M)^{d/2}}.
\]

On the other hand, \eqref{eq:fp-normalizer}
gives
\[
e^{-M}Z(a)
\ge
\frac{e^{-1/2}v_d}
{(2R)^d(1+M)^{d/2}}
\ge
\frac{e^{-1/2}v_d}
{(2R)^d(1+\overline M)^{d/2}}.
\]
Consequently, the absolute errors in the shifted
denominator integral and each normalized numerator
integral are bounded by
$\delta_\varepsilon e^{-M}Z(a)/(4B)$.
Since $0\le\phi/B\le1$, taking their ratio gives
\[
\begin{aligned}
|\widehat F_i(a)-F_i(a)|
&\le
\frac{\delta_\varepsilon/2}
{1-\delta_\varepsilon/(4B)}\\
&\le\delta_\varepsilon.
\end{aligned}
\]
The shift $M$ was introduced only in the error
analysis. The common exponential factor cancels
from the ratio, so the same estimate applies to
\eqref{eq:fc-quadrature}, where the shift is the
maximum of $V_a$ over the quadrature nodes.

We now prove \eqref{eq:fc-evaluation-error}
for every evaluation by induction over the
proximal steps. At $a^0=0$, we have $V_{a^0}=0$,
so the preceding estimate applies.

Suppose that all evaluations in the first $k$
completed steps satisfy \eqref{eq:fc-evaluation-error}.
By Lemma~\ref{lem:fp-inner}, these steps satisfy
\eqref{eq:fp-residual}, and hence
Lemma~\ref{lem:fp-potential} gives
\[
\sup_{\mathbb R^d}V_{a^k}\le\mathcal M(t_k).
\]
The global rule used in the next step is
constructed with $\overline M\ge\mathcal M(t_k+1)$.
The preceding quadrature estimate therefore
establishes \eqref{eq:fc-evaluation-error}
at $a^k$, allowing $p^k$ and $\tau_k$ to
be computed by \eqref{eq:fc-step}.

For every $y\in[1/2,2]^n$, the definition
of $\tau_k$ and the Gaussian grid estimate give
\[
0\le V_{x(y)}-V_{a^k}
\le2\tau_k\|\widehat\Psi\|_\infty
\le\frac1{32}.
\]
Consequently,
\[
\sup_{\mathbb R^d}V_{x(y)}
\le\mathcal M(t_k)+\frac1{32}
\le\mathcal M(t_k+1)
\le\overline M.
\]
Thus \eqref{eq:fc-evaluation-error} holds
at every trial vector $x(y)$ with
$y\in[1/2,2]^n$.

The fixed-point iteration starts from
$y^0=\mathbf1$. By the mapping estimates
in the proof of Lemma~\ref{lem:fp-inner},
the computed map satisfies
\[
\widehat T([1/2,2]^n)
\subset[3/4,4/3]^n
\subset[1/2,2]^n.
\]
It follows by induction over the inner iterations
that every trial vector satisfies the preceding
potential bound and hence
\eqref{eq:fc-evaluation-error}.
Lemma~\ref{lem:fp-inner} then gives
\eqref{eq:fp-residual} for the accepted step.
This completes the induction over the
proximal steps.

Finally, each fitted-value evaluation involves
$q_K^d$ tensor nodes, giving
$O(ndq_K^d)$ operations. The construction cost
follows from Lemma~\ref{lem:fp-rounded}.
\end{proof}

\subsection{Proof of Theorem~\ref{thm:fc-global-cost}}
\begin{proof}[\textbf{Proof of Theorem~\ref{thm:fc-global-cost}}]
The grid $\mathcal G$ in \eqref{eq:fc-score-grid}
contains $O_d(R^d)$ points. Computing the kernel
values $\phi(X_i-\theta)$ for all
$i=1,\ldots,n$ and $\theta\in\mathcal G$
requires $O_d(ndR^d)$ operations.
Once these values are stored, each evaluation
of \eqref{eq:fc-step} requires $O_d(nR^d)$
operations.

By the stopping condition of Algorithm~
\ref{alg:fc-global},
$t_N<32H_\varepsilon/\varepsilon+1/64$.
Since $\overline M$ is the smallest power of two
greater than $\mathcal M(t_k+1)$, we have
\[
\overline M
\le2\mathcal M(32H_\varepsilon/\varepsilon+2)
\]
throughout the run. It follows from
\eqref{eq:fc-global-order} and
\eqref{eq:fp-global-tolerance} that
\[
q_K=\widetilde O_d(R\varepsilon^{-1/2}),
\quad
\log(e/\nu_K)=O_d(\log(enR/\varepsilon)).
\]
In particular, each fitted-value evaluation
requires
$\widetilde O_d(ndR^d\varepsilon^{-d/2})$
operations.

By Theorem~\ref{thm:fc-iterations}, the algorithm
performs at most
$\widetilde O_d(R^d\varepsilon^{-1})$
proximal steps, with
$O_d(\log(enR/\varepsilon))$ fitted-value
evaluations per step. The total cost of these
evaluations is therefore
\[
\widetilde O_d
\left(nR^{2d}\varepsilon^{-1-d/2}\right).
\]

It remains to account for the construction of
the quadrature rules. A new rule is constructed
only when $\overline M$ increases, in which case
it is doubled. Since
$q_K^2=\widetilde O_d(R^2(1+\overline M))$,
the total construction cost is bounded by
\[
\begin{aligned}
\sum_{\text{constructed rules}}
O(q_K^2\log(e/\nu_K))
&=\widetilde O_d
\left(R^2\sum_{\text{constructed rules}}
(1+\overline M)\right)\\
&=\widetilde O_d(R^2\varepsilon^{-1}).
\end{aligned}
\]
The costs of evaluating the score on $\mathcal G$,
constructing the kernel tables, and forming the
output measure are bounded by the preceding terms.
Thus the total arithmetic complexity is
\[
\widetilde O_d\left(
nR^{2d}\varepsilon^{-1-d/2}
+R^2\varepsilon^{-1}
\right),
\]
as claimed in \eqref{eq:fc-global-cost}.
\end{proof}

\section{Proofs of Results in Section \ref{sec:adaptive}}\label{sec:adaptive proofs}

Throughout this section, we associate each completed iterate
$a^k$ with the accumulated proximal time $t_k$. For trial
coefficient vectors used in the computation of $a^{k+1}$,
we retain the time $t_k$. All estimates below follow this
convention.

We define the following auxiliary functions to assist in our analysis.
\begin{equation}\label{eq:fa-error-functions}
 \begin{split}
 g_\varepsilon(t)&=
 \frac{\log((2R)^d/v_d)+1/2+(d/2)\log(1+t)}{t}
                           +2\eta_\varepsilon,\\
 E_\varepsilon(t)&=
 32\left[3\sqrt{ng_\varepsilon(t)}
       +\frac{4ng_\varepsilon(t)}{\log2-1/2}
       +\eta_\varepsilon+\frac1{32t}\right],\quad t\ge1.
 \end{split}
\end{equation}
Here $g_\varepsilon(t)$ bounds an average likelihood gap, while
$E_\varepsilon(t)$ controls the differences between $V_a/t$ and the
optimal score, including their first and second derivatives. Their
use in the acceptance inequalities ensures that the selected cubes
remain valid at subsequent proximal steps and inner trials.

\subsection{Convergence of the rescaled potentials}

\begin{lemma}[Approximation of the optimal score]
\label{lem:ap-potential}
Let \(a^0=0,a^1,\ldots,a^k\) be the completed proximal iterates, with \(t_k\ge1\). Suppose that every fitted-value approximation used to compute \(a^1,\ldots,a^k\) satisfies
\[
\max_{1\le i\le n}
\bigl|\widehat F_i(u)-F_i(u)\bigr|
\le\delta_\varepsilon.
\]
Let \(p_i^k=\max\{\widehat F_i(a^k),b\}\), and choose \(\tau_k\) according to \eqref{eq:fc-step}. For \(y\in[1/2,2]^n\), define
\[
x_i(y)=a_i^k+\frac{\tau_k y_i}{np_i^k},
\qquad 1\le i\le n.
\]
Set \(t=t_k\). For \(a=a^k\), and for every \(a=x(y)\) with \(y\in[1/2,2]^n\), the rescaled potential \(W_{a,t}=V_a/t\) satisfies
\[
\sup_{\theta\in\mathbb R^d}
\bigl|W_{a,t}(\theta)-\Psi^\star(\theta)\bigr|
\le E_\varepsilon(t),
\]
\[
\sup_{\theta\in\mathbb R^d}
\bigl\|\nabla W_{a,t}(\theta)-\nabla\Psi^\star(\theta)\bigr\|
\le E_\varepsilon(t),
\]
and
\[
\sup_{\substack{z\in\mathbb C^d\\ \|\operatorname{Im}z\|\le1}}
\bigl\|\nabla^2W_{a,t}(z)-\nabla^2\Psi^\star(z)\bigr\|_{\mathrm{op}}
\le E_\varepsilon(t).
\]
Moreover, \(E_\varepsilon\) is decreasing on \([1,\infty)\). 
\end{lemma}

\begin{proof}
Let $a=a^k$ be a completed iterate and write
$t=t_k\ge1$. For $0\le j<k$, set
\[
\Delta_{j+1}
=\widehat L_n-\mathcal L_b(\mu_{a^{j+1}}),
\quad
\overline f_i^{j+1}
=\max\{F_i(a^{j+1}),b\}.
\]
We first bound the weighted average of
$\Delta_1,\ldots,\Delta_k$.

Let $\mu^\star$ be an NPMLE supported on
$\operatorname{conv}\{X_1,\ldots,X_n\}$, and define
\[
\pi_t=\mu^\star*
\operatorname{Unif}\bigl(B(0,(1+t)^{-1/2})\bigr).
\]
By \eqref{eq:fp-smoothing}, this measure is
supported on $K$ and satisfies
\[
\begin{aligned}
\operatorname{KL}(\pi_t\|\mu_0)
\le\log\frac{(2R)^d}{v_d}
+\frac d2\log(1+t), \ \text{and }
\widehat L_n-\mathcal L_b(\pi_t)
\le\frac1{2(1+t)}.
\end{aligned}
\]
Applying \eqref{eq:fp-entropy} with 
$\pi_t$ and summing over $0\le j<k$, we obtain
\[
\begin{aligned}
\sum_{j<k}\tau_j
\{\mathcal L_b(\pi_t)
-\mathcal L_b(\mu_{a^{j+1}})\}
\le
\operatorname{KL}(\pi_t\|\mu_0)
-\operatorname{KL}(\pi_t\|\mu_{a^k})
-\sum_{j<k}
\operatorname{KL}(\mu_{a^{j+1}}\|\mu_{a^j})
+2\eta_\varepsilon t\\
\le
\log\frac{(2R)^d}{v_d}
+\frac d2\log(1+t)+2\eta_\varepsilon t.
\end{aligned}
\]
Since $\sum_{j<k}\tau_j=t$, it follows that

\begin{equation}
\begin{aligned}
\label{eq:weighted-gap}
\frac1t\sum_{j<k}\tau_j\Delta_{j+1}
\le\widehat L_n-\mathcal L_b(\pi_t)
+\frac{\operatorname{KL}(\pi_t\|\mu_0)}t
+2\eta_\varepsilon \\
\le
\frac{\log((2R)^d/v_d)+1/2+(d/2)\log(1+t)}t
+2\eta_\varepsilon 
&=g_\varepsilon(t).
\end{aligned}
\end{equation}

We next use this estimate to compare
$V_{a^k}/t$ with $\Psi^\star$.
By \eqref{eq:fp-residual}, the accepted coefficient
updates can be written as
\[
a_i^{j+1}-a_i^j
=\tau_j\left\{
\frac1{n\overline f_i^{j+1}}+r_i^j
\right\},
\quad
B\sum_{i=1}^n|r_i^j|\le\eta_\varepsilon.
\]
Since $a^0=0$, summing these identities gives
\[
\begin{aligned}
\frac{a_i^k}{t}-\frac1{nf_i^\star}
=\frac1t\sum_{j<k}\tau_j
\left\{
\frac1n\left(
\frac1{\overline f_i^{j+1}}
-\frac1{f_i^\star}
\right)+r_i^j
\right\}.
\end{aligned}
\]
Consequently, if
\[
P_k(\theta)=
\sum_{i=1}^n
\left|\frac{a_i^k}{t}-\frac1{nf_i^\star}\right|
\phi(X_i-\theta),
\]
then, for every $\theta\in\mathbb R^d$,
\[
\begin{aligned}
P_k(\theta)
&\le\frac1t\sum_{j<k}\tau_j
\left\{
\frac1n\sum_{i=1}^n
\left|
\frac1{\overline f_i^{j+1}}
-\frac1{f_i^\star}
\right|\phi(X_i-\theta)
+B\sum_{i=1}^n|r_i^j|
\right\}\\
&\le\frac1t\sum_{j<k}\tau_j
\left\{
3\sqrt{n\Delta_{j+1}}
+\frac{4n\Delta_{j+1}}{\log2-1/2}
+\eta_\varepsilon
\right\},
\end{aligned}
\]
where the second inequality follows from
\eqref{eq:fp-absolute-score}.

By Cauchy--Schwarz and \eqref{eq:weighted-gap},
\[
\begin{aligned}
\frac1t\sum_{j<k}\tau_j\sqrt{\Delta_{j+1}}
\le
\left(\frac1t\sum_{j<k}
\tau_j\Delta_{j+1}\right)^{1/2}
\le\sqrt{g_\varepsilon(t)}.
\end{aligned}
\]
We therefore obtain
\[
\sup_{\theta\in\mathbb R^d}P_k(\theta)
\le
3\sqrt{ng_\varepsilon(t)}
+\frac{4ng_\varepsilon(t)}{\log2-1/2}
+\eta_\varepsilon.
\]

Now let $y^r\in[1/2,2]^n$ be arbitrary and define $a=x^r$ by
$
x_i^r=a_i^k+\frac{\tau_k y_i^r}{np_i^k}$
for $1\le i\le n$, and set $t=t_k$. Since $x_i^r\ge a_i^k$, we have
\[
\begin{aligned}
&\sum_{i=1}^n
\left|\frac{x_i^r}{t}-\frac1{nf_i^\star}\right|
\phi(X_i-\theta)
\le
P_k(\theta)
+\frac{V_{x^r}(\theta)-V_{a^k}(\theta)}t.
\end{aligned}
\]
By \eqref{eq:fp-trial-bound},
$0\le V_{x^r}-V_{a^k}\le1/32$.
This bound holds for every $y^r\in[1/2,2]^n$.
Thus, for either $a=a^k$ or $a=x^r$,
\[
\begin{aligned}
\sup_{\theta\in\mathbb R^d}
\sum_{i=1}^n
\left|\frac{a_i}{t}-\frac1{nf_i^\star}\right|
\phi(X_i-\theta)
\le
3\sqrt{ng_\varepsilon(t)}
+\frac{4ng_\varepsilon(t)}{\log2-1/2}
+\eta_\varepsilon+\frac1{32t}
=\frac{E_\varepsilon(t)}{32}.
\end{aligned}
\]

Finally, recall that
\[
W_{a,t}(\theta)-\Psi^\star(\theta)
=\sum_{i=1}^n
\left(\frac{a_i}{t}-\frac1{nf_i^\star}\right)
\phi(X_i-\theta).
\]
Applying the derivative bounds for signed Gaussian
sums in Lemma~\ref{lem:fp-gaussian}, we obtain
\[
\begin{aligned}
\sup_{\theta\in\mathbb R^d}
|W_{a,t}(\theta)-\Psi^\star(\theta)|
&\le\frac{E_\varepsilon(t)}{32},\\
\sup_{\theta\in\mathbb R^d}
\|\nabla W_{a,t}(\theta)-\nabla\Psi^\star(\theta)\|
&\le\frac{E_\varepsilon(t)}{16},\\
\sup_{\substack{z\in\mathbb C^d\\
\|\operatorname{Im}z\|\le1}}
\|\nabla^2W_{a,t}(z)-\nabla^2\Psi^\star(z)\|_{\rm op}
&\le E_\varepsilon(t).
\end{aligned}
\]
This proves the three asserted estimates.

To verify monotonicity, differentiating
$g_\varepsilon$ gives
\[
g_\varepsilon'(t)
=\frac1{t^2}
\left\{
\frac{dt}{2(1+t)}
-\log\frac{(2R)^d}{v_d}
-\frac12-\frac d2\log(1+t)
\right\}<0,
\]
since $(2R)^d\ge v_d$ and
$\log(1+t)>t/(1+t)$ for $t>0$.
The definition of $E_\varepsilon$ then shows
that it is decreasing on $[1,\infty)$.
\end{proof}

\subsection{Persistence and finite selection of the cubes}

Suppose that \textsc{SelectCubes} has accepted centers
$c_1,\ldots,c_m$ and radius $r$. For each
$1\le\ell\le m$, define
\begin{equation}\label{eq:ap-complex-tube}
\mathcal T_\ell
=\left\{z\in\mathbb C^d:
\|\operatorname{Re}z-c_\ell\|_\infty\le2r,\ \ 
\|\operatorname{Im}z\|_\infty\le\frac r{\sqrt d}
\right\}.
\end{equation}

\begin{lemma}[Consequences of the acceptance inequalities]
\label{lem:ap-persistence}
Suppose that \textsc{SelectCubes} accepts centers
$c_1,\ldots,c_m$, radius $r$, and parameters
$\kappa,\gamma$ at time $t_{\rm sel}$.
Let $\mathcal U=\bigcup_{\ell=1}^m U_\ell$, and
suppose that all fitted-value evaluations up to
the iterate under consideration satisfy
\eqref{eq:fc-evaluation-error}.

For every completed iterate $a=a^k$ with
$t=t_k\ge t_{\rm sel}$, and every trial vector
$a=x^r$ in the computation of $a^{k+1}$ with
$t=t_k\ge t_{\rm sel}$, the cubes
$U_1,\ldots,U_m$ are pairwise disjoint and
\begin{align}
\kappa I_d\preceq-\nabla^2W_{a,t}(\theta)
&\preceq\Lambda I_d,
&&\theta\in\mathcal U,
\label{eq:ap-curvature}\\
\|\nabla W_{a,t}(c_\ell)\|
&\le\frac{\kappa r}{4},
&&1\le\ell\le m,
\label{eq:ap-center-gradient}\\
\sup_{\theta\in K\setminus\mathcal U}
W_{a,t}(\theta)
&\le\max_{1\le\ell\le m}W_{a,t}(c_\ell)-\gamma,
\label{eq:ap-exterior-gap}\\
\|\nabla^2W_{a,t}(z)\|_{\rm op}
&\le\Lambda,
&&z\in\mathcal T_\ell,\quad 1\le\ell\le m
\label{eq:ap-complex-hessian}
\end{align}
where $\Lambda=32\mathcal M(1)$.
\end{lemma}

\begin{proof}
Let $a^{k_0}$ be the iterate at which
\textsc{SelectCubes} accepts the cubes, and write
$W_0=V_{a^{k_0}}/t_{\rm sel}$.
By Lemmas~\ref{lem:fp-potential} and
\ref{lem:fp-gaussian}, we have
\[
\sup_{\mathbb R^d}W_0\le\mathcal M(1),
\quad
\|\nabla W_0\|_\infty\le2\mathcal M(1),
\quad
\|D^3W_0\|_{\rm op}\le48\mathcal M(1).
\]

Recall that the accepted partition in
\textsc{SelectCubes} consists of cubes
$C_u=c_u+[-s/2,s/2]^d$, where $s=2R\,2^{-h}$.
For every $\theta\in C_u$, we have
$\|\theta-c_u\|\le\sqrt d\,s/2$. Hence
\[
\begin{aligned}
\|\nabla^2W_0(\theta)-\nabla^2W_0(c_u)\|_{\rm op}
&\le24\mathcal M(1)\sqrt d\,s,\\
|W_0(\theta)-W_0(c_u)|
&\le\mathcal M(1)\sqrt d\,s.
\end{aligned}
\]
By the first inequality in \eqref{eq:fa-accept},
every cell intersecting $\mathcal U$ satisfies
\[
-\nabla^2W_0(c_u)
\succ(4\kappa+24\mathcal M(1)\sqrt d\,s)I_d.
\]
Consequently,
\[
-\nabla^2W_0(\theta)\succeq4\kappa I_d,
\quad \theta\in\mathcal U.
\]
Similarly, the cells used in \eqref{eq:fa-gap}
cover $K\setminus\mathcal U$, so the second
estimate and the definition of $\gamma$ give
\[
\sup_{K\setminus\mathcal U}W_0
\le\max_\ell W_0(c_\ell)-4\gamma.
\]
The retained centers satisfy
$\|c_\ell-c_j\|_\infty>5r$ for $\ell\ne j$.
Since $U_\ell=c_\ell+[-2r,2r]^d$, these cubes
are pairwise disjoint.

Now let $(a,t)$ correspond to a subsequent
completed iterate or trial vector. By
Lemma~\ref{lem:ap-potential} and the monotonicity
of $E_\varepsilon$, we have
\[
\begin{aligned}
\|W_{a,t}-W_0\|_\infty
&\le2E_\varepsilon(t_{\rm sel}),\\
\|\nabla W_{a,t}-\nabla W_0\|_\infty
&\le2E_\varepsilon(t_{\rm sel}),\\
\sup_{\|\operatorname{Im}z\|\le1}
\|\nabla^2W_{a,t}(z)-\nabla^2W_0(z)\|_{\rm op}
&\le2E_\varepsilon(t_{\rm sel}).
\end{aligned}
\]
Since $r\le1/8$, the last inequality in
\eqref{eq:fa-accept} implies
\[
2E_\varepsilon(t_{\rm sel})
\le\frac{\kappa r}{8}
\le\frac{\kappa}{64},
\quad
4E_\varepsilon(t_{\rm sel})\le\gamma.
\]
It follows that
\[
-\nabla^2W_{a,t}(\theta)
\succeq(4\kappa-2E_\varepsilon(t_{\rm sel}))I_d
\succeq\kappa I_d,
\quad \theta\in\mathcal U.
\]
Moreover, $\sup W_{a,t}\le\mathcal M(1)$,
so Lemma~\ref{lem:fp-gaussian} gives
$\|\nabla^2W_{a,t}\|_{\rm op}
\le4\mathcal M(1)\le\Lambda$.
This proves \eqref{eq:ap-curvature}.

For each retained center, the second inequality
in \eqref{eq:fa-accept} gives
$\|\nabla W_0(c_\ell)\|<\kappa r/8$.
Therefore,
\[
\begin{aligned}
\|\nabla W_{a,t}(c_\ell)\|
\le\|\nabla W_0(c_\ell)\|
+2E_\varepsilon(t_{\rm sel})
\le\frac{\kappa r}{4},
\end{aligned}
\]
which proves \eqref{eq:ap-center-gradient}.
Similarly,
\[
\begin{aligned}
\sup_{K\setminus\mathcal U}W_{a,t}
&\le\sup_{K\setminus\mathcal U}W_0
+2E_\varepsilon(t_{\rm sel})
\le\max_\ell W_0(c_\ell)-4\gamma
+2E_\varepsilon(t_{\rm sel})\\
&\le\max_\ell W_{a,t}(c_\ell)-4\gamma
+4E_\varepsilon(t_{\rm sel})
\le\max_\ell W_{a,t}(c_\ell)-\gamma,
\end{aligned}
\]
establishing \eqref{eq:ap-exterior-gap}.

Finally, for $z\in\mathcal T_\ell$, we have
$\|\operatorname{Im}z\|\le r\le1$.
Since $\sup_{\mathbb R^d}W_{a,t}\le\mathcal M(1)$,
Lemma~\ref{lem:fp-gaussian} gives
\[
\|\nabla^2W_{a,t}(z)\|_{\rm op}
\le32\mathcal M(1)=\Lambda.
\]
\end{proof}

\begin{lemma}[Finite selection under nondegeneracy]
\label{lem:ap-selection}
Suppose that the data satisfy \eqref{eq:fa-morse}.
Then there exist $\varepsilon_0>0$, $T_0<\infty$,
and $C_{X,d}<\infty$, depending only on
$X_1,\ldots,X_n$ and $d$, such that for every
$0<\varepsilon\le\varepsilon_0$,
\textsc{SelectCubes} accepts a tuple at some
$t_{\rm sel}\le T_0$. Moreover, the accepted
parameters satisfy
\[
m\le C_{X,d},\quad
r^{-1}+\kappa^{-1}+\gamma^{-1}+\Lambda
\le C_{X,d}.
\]
\end{lemma}

\begin{proof}
The contact set $\mathcal C^\star$ is nonempty,
compact, and contained in the data hull. Under
\eqref{eq:fa-morse}, every contact point is isolated,
so $\mathcal C^\star=\{\xi_1,\ldots,\xi_s\}$ is finite.
By continuity of $\nabla^2\Psi^\star$, there exist
$\rho_0,\lambda,\gamma_0>0$ such that the balls
$B(\xi_\ell,\rho_0)$ are disjoint, contained in
$\operatorname{int}K$, and
\begin{equation}\label{eq:ap-contact-bounds}
\begin{gathered}
8\lambda I_d\preceq-\nabla^2\Psi^\star
\preceq4I_d
\quad\text{on }\bigcup_\ell B(\xi_\ell,\rho_0),\\
4\lambda\|\theta-\xi_\ell\|^2
\le1-\Psi^\star(\theta)
\le2\|\theta-\xi_\ell\|^2
\quad(\theta\in B(\xi_\ell,\rho_0)),\\
\Psi^\star\le1-\gamma_0
\quad\text{on }K\setminus
\bigcup_\ell B(\xi_\ell,\rho_0).
\end{gathered}
\end{equation}
Here the upper Hessian bound follows from
Lemma~\ref{lem:fp-gaussian}, since
$\sup\Psi^\star=1$, and the quadratic bounds
follow from Taylor's formula and
$\nabla\Psi^\star(\xi_\ell)=0$.

Choose $p_\star$ such that
$\kappa=2^{-p_\star}\le\min\{1,\lambda\}$.
Consider the partition at level $h$, with side
length $s_h=2R2^{-h}$ and radius
$r=s_h^{1/4}$. We fix $h$ sufficiently large
in the argument below. This is possible since
$r\to0$ and $s_h=o(r^2)$ as $h\to\infty$.

Suppose that the partition is examined at a
checkpoint $(a,t)$ satisfying
\begin{equation}\label{eq:ap-selection-error}
E_\varepsilon(t)\le s_h^2.
\end{equation}
Write $W=V_a/t$. By
Lemma~\ref{lem:ap-potential},
\[
\|W-\Psi^\star\|_\infty\le s_h^2,\quad
\|\nabla W-\nabla\Psi^\star\|_\infty\le s_h^2,
\quad
\|\nabla^2W-\nabla^2\Psi^\star\|_\infty
\le s_h^2.
\]

For each $\xi_\ell$, there is a cell midpoint
$c_u$ with
$\|c_u-\xi_\ell\|\le\sqrt d\,s_h/2$.
By \eqref{eq:ap-contact-bounds},
\[
W(c_u)\ge1-(d/2+1)s_h^2.
\]
Since $\max_vW(c_v)\le1+s_h^2$, this midpoint
belongs to the high-value set in Step 2 of
\textsc{SelectCubes} for sufficiently large $h$.

Conversely, every midpoint in the high-value set
satisfies
\[
\begin{aligned}
1-\Psi^\star(c_u)
&\le1-\max_vW(c_v)
+2\mathcal M(1)\sqrt d\,s_h+s_h^2\\
&\le(d/2+2)s_h^2
+2\mathcal M(1)\sqrt d\,s_h
=O_{X,d}(s_h).
\end{aligned}
\]
Since $s_h=o(r^2)$, for sufficiently large $h$,
\[
1-\Psi^\star(c_u)
<\min\left\{\gamma_0,\,
4\lambda\left(\frac{\kappa r}{128}\right)^2\right\}.
\]
The first bound and \eqref{eq:ap-contact-bounds}
imply that $c_u\in B(\xi_\ell,\rho_0)$
for some $\ell$. Applying the quadratic lower
bound in \eqref{eq:ap-contact-bounds}, we obtain
\[
4\lambda\|c_u-\xi_\ell\|^2
\le1-\Psi^\star(c_u)
<4\lambda\left(\frac{\kappa r}{128}\right)^2.
\]
Consequently,
$\|c_u-\xi_\ell\|<\kappa r/128$,
where $\ell$ is unique since the balls are disjoint.

For sufficiently large $h$, any two midpoints
within $\kappa r/128$ of the same contact point
are separated by less than $5r$, while midpoints
near distinct contact points are separated by more
than $5r$. Recall that Step 2 of
\textsc{SelectCubes} retains a midpoint only if
its distance in the infinity norm from every
previously retained midpoint exceeds $5r$.
Since at least one midpoint satisfying the
selection inequality lies near each contact point,
the procedure retains exactly one midpoint near
each $\xi_\ell$. Relabeling these centers, we obtain
\[
m=s,\quad
\|c_\ell-\xi_\ell\|<\frac{\kappa r}{128},
\quad 1\le\ell\le s.
\]
Since the contact points are distinct and belong
to the interior of $K$, the cubes
$U_\ell=c_\ell+[-2r,2r]^d$ are disjoint, contained
in the respective balls $B(\xi_\ell,\rho_0)$,
and satisfy the interior condition in Step 3
for sufficiently large $h$.

We now verify the acceptance inequalities.
For every cell $C_u$ intersecting $\mathcal U$,
its midpoint $c_u$ belongs to one of the balls
in \eqref{eq:ap-contact-bounds}, provided $h$
is sufficiently large. Thus
\[
\begin{aligned}
-\nabla^2W(c_u)
\succeq-\nabla^2\Psi^\star(c_u)-s_h^2I_d
\succeq(8\lambda-s_h^2)I_d
\succ(4\kappa+
24\mathcal M(1)\sqrt d\,s_h)I_d.
\end{aligned}
\]
Similarly, for every retained center,
\[
\begin{aligned}
\|\nabla W(c_\ell)\|
\le\|\nabla\Psi^\star(c_\ell)\|+s_h^2
\le4\|c_\ell-\xi_\ell\|+s_h^2
<\frac{\kappa r}{32}+s_h^2
<\frac{\kappa r}{8}.
\end{aligned}
\]

For every $u\in\mathcal O$, the definition of
$\mathcal O$ gives
\[
\|c_u-\xi_\ell\|_\infty
\ge2r-\frac{s_h}{2}-\frac{\kappa r}{128}
>r,\quad 1\le\ell\le s,
\]
for sufficiently large $h$.
Using \eqref{eq:ap-contact-bounds} and
$4\lambda r^2<\gamma_0$, we obtain
\[
\max_{u\in\mathcal O}W(c_u)
\le1-4\lambda r^2+s_h^2.
\]
On the other hand,
\[
\max_\ell W(c_\ell)
\ge1-2(\kappa r/128)^2-s_h^2.
\]
Now, Since \(\kappa\le\min\{1,\lambda\}\), we have \(\kappa^2\le\lambda\). Therefore
\[
\begin{aligned}
4\gamma
\ge
\left(4\lambda-\frac{\kappa^2}{8192}\right)r^2-o(r^2)
\ge
\left(4-\frac1{8192}\right)\lambda r^2-o(r^2)
\ge2\lambda r^2
\end{aligned}
\]
for sufficiently large $h$. Hence, $
\gamma\ge\lambda r^2/2
$ for sufficiently large $h$.
Since $s_h=o(r^2)$, \eqref{eq:ap-selection-error}
implies
\[
E_\varepsilon(t)\le
\min\{\gamma/4,\kappa r/16\}.
\]
Thus all the inequalities in \eqref{eq:fa-accept}
hold for this fixed pair $(h,p_\star)$.

It remains to show that such a checkpoint occurs
within a time independent of $\varepsilon$.
Let $E_0(t)$ denote the expression obtained from
$E_\varepsilon(t)$ by setting $\eta_\varepsilon=0$.
Since $E_0(t)\to0$ as $t\to\infty$, we can choose
$j_\star\ge\max\{4h,p_\star\}$ such that
\[
E_0(2^{j_\star})<s_h^2.
\]
Since $\eta_\varepsilon\to0$ and
$32H_\varepsilon/\varepsilon\to\infty$ as
$\varepsilon\to0$, there exists
$\varepsilon_0>0$ such that the algorithm reaches
this checkpoint and \eqref{eq:ap-selection-error}
holds whenever $0<\varepsilon\le\varepsilon_0$.
As $\tau_k\le1/64$, selection has occurred by
$T_0=2^{j_\star}+1/64$.

Finally, since $t_{\rm sel}\le T_0$, the
checkpoint index satisfies $2^j\le T_0$.
The procedure examines only
\[
1\le h\le\max\{1,\lfloor j/4\rfloor\},
\quad 0\le p\le j+1.
\]
Thus only finitely many pairs $(h,p)$ can be
accepted. Since $r$ and $\kappa$ are determined
by $(h,p)$, and the number of retained centers
satisfies $m\le2^{dh}$, the quantities
$m$, $r^{-1}$, and $\kappa^{-1}$ are bounded
independently of $\varepsilon$.

Moreover, the acceptance inequalities give
\[
\gamma\ge4E_\varepsilon(t_{\rm sel})
\ge4E_0(T_0)>0,
\]
where the second inequality follows from
$t_{\rm sel}\le T_0$ and the monotonicity of
$E_0$. Finally, $\Lambda\le32\mathcal M(1)$.
These bounds prove the lemma.
\end{proof}

\subsection{Local maxima and quadrature errors}

\begin{lemma}[Projected gradient computation of the maxima]
\label{lem:ap-centers}
Suppose that \eqref{eq:ap-curvature}--\eqref{eq:ap-exterior-gap}
hold for $W_{a,t}$. Then, for each $1\le\ell\le m$,
$W_{a,t}$ has a unique maximizer $z_\ell^\star$
on $U_\ell$, satisfying
\[
\|z_\ell^\star-c_\ell\|\le\frac r4.
\]
The points $z_\ell$ computed by the projected
gradient iterations in \eqref{eq:fa-ascent} satisfy
\begin{equation}\label{eq:ap-center-error}
\|z_\ell-z_\ell^\star\|
\le\frac1{4\sqrt{\Lambda t}},
\quad
\|z_\ell-c_\ell\|\le\frac r2.
\end{equation}
Moreover,
\begin{equation}\label{eq:ap-shift-error}
0\le\max_{\theta\in K}V_a(\theta)
-\max_{1\le\ell\le m}V_a(z_\ell)
\le\frac1{32}.
\end{equation}
\end{lemma}

\begin{proof}
Fix $1\le\ell\le m$ and write $W=W_{a,t}$.
For $x\in\partial B(c_\ell,r/2)$,
\eqref{eq:ap-curvature} and
\eqref{eq:ap-center-gradient} imply
\begin{align*}
\langle\nabla W(x),x-c_\ell\rangle
=\langle\nabla W(c_\ell),x-c_\ell\rangle
+\int_0^1(x-c_\ell)^{\mathsf T}
\nabla^2W(c_\ell+s(x-c_\ell))
(x-c_\ell)\,ds\\
\le\frac r2\|\nabla W(c_\ell)\|
-\frac{\kappa r^2}{4}
\le-\frac{\kappa r^2}{8}.
\end{align*}
Thus the directional derivative of $W$ towards
$c_\ell$ is positive at every point of
$\partial B(c_\ell,r/2)$. The maximum of $W$
on the closed ball is therefore attained at an
interior point $z_\ell^\star$, where
$\nabla W(z_\ell^\star)=0$.
Since $W$ is strongly concave on $U_\ell$,
this is its unique maximizer on $U_\ell$.

Moreover, by \eqref{eq:ap-curvature},
\[
\begin{aligned}
\kappa\|z_\ell^\star-c_\ell\|^2
\le\langle\nabla W(c_\ell)
-\nabla W(z_\ell^\star),
z_\ell^\star-c_\ell\rangle
\le\|\nabla W(c_\ell)\|
\|z_\ell^\star-c_\ell\|.
\end{aligned}
\]
Together with \eqref{eq:ap-center-gradient},
this gives
$\|z_\ell^\star-c_\ell\|\le r/4$.
In particular, $z_\ell^\star\in\operatorname{int}Q_\ell$.

We next consider the projected gradient
iterations in \eqref{eq:fa-ascent}.
For $x,y\in Q_\ell$, 
$\kappa I_d\preceq-\nabla^2W\preceq\Lambda I_d$
imply
\[
\begin{aligned}
\left\|
x-y+\frac1\Lambda
(\nabla W(x)-\nabla W(y))
\right\|
=
\left\|
\left(I_d+\frac1\Lambda
\int_0^1\nabla^2W(y+s(x-y))\,ds\right)
(x-y)\right\| \\
\le
\left(1-\frac\kappa\Lambda\right)\|x-y\|.
\end{aligned}
\]
Since projection onto $Q_\ell$
is nonexpansive and
$\nabla W(z_\ell^\star)=0$, we obtain
\[
\begin{aligned}
\|z_\ell^v-z_\ell^\star\|
\le\left(1-\frac\kappa\Lambda\right)^v
\|c_\ell-z_\ell^\star\|
\le\frac r4e^{-v\kappa/\Lambda}.
\end{aligned}
\]
By the definition of $J(t)$ in
\eqref{eq:fa-ascent},
\[
\frac r4e^{-J(t)\kappa/\Lambda}
\le\frac1{4\sqrt{\Lambda t}}.
\]
Indeed, for $r\sqrt{\Lambda t}>1$,
we have
$J(t)\ge(\Lambda/\kappa)\log(r\sqrt{\Lambda t})$,
while for $r\sqrt{\Lambda t}\le1$,
we have $J(t)=0$ and
$r/4\le1/(4\sqrt{\Lambda t})$.
Consequently,
\[
\|z_\ell-z_\ell^\star\|
\le\frac1{4\sqrt{\Lambda t}},
\quad
\|z_\ell-c_\ell\|
\le\frac r4+\frac r4=\frac r2,
\]
which proves \eqref{eq:ap-center-error}.

Since $\nabla W(z_\ell^\star)=0$,
Taylor's theorem and \eqref{eq:ap-curvature}
give
\[
-\frac{\Lambda t}{2}\|\theta-z_\ell^\star\|^2
\le V_a(\theta)-V_a(z_\ell^\star)
\le-\frac{\kappa t}{2}\|\theta-z_\ell^\star\|^2,
\quad \theta\in U_\ell.
\]
In particular,
\[
0\le V_a(z_\ell^\star)-V_a(z_\ell)
\le\frac{\Lambda t}{2}
\|z_\ell-z_\ell^\star\|^2
\le\frac1{32}.
\]
By \eqref{eq:ap-exterior-gap}, every global
maximum of $V_a$ on $K$ belongs to
$\mathcal U$. Hence
\[
\max_KV_a=\max_{1\le\ell\le m}V_a(z_\ell^\star),
\]
and therefore
\[
0\le\max_KV_a-\max_\ell V_a(z_\ell)
\le\frac1{32},
\]
proving \eqref{eq:ap-shift-error}.

Finally, the upper Hessian bound also gives
\[
\begin{aligned}
\|\nabla V_a(z_\ell)\|
=\|\nabla V_a(z_\ell)
-\nabla V_a(z_\ell^\star)\|
\le\Lambda t\|z_\ell-z_\ell^\star\|
\le\frac{\sqrt{\Lambda t}}4.
\end{aligned}
\]
\end{proof}

Let $q_{\rm loc}$ be given by \eqref{eq:fa-order}.
Using Lemma~\ref{lem:fp-rounded}, construct a
one-dimensional Gauss--Legendre rule on $[-1,1]$
with order $q_{\rm loc}$ and tolerance
\begin{equation}\label{eq:ap-local-tolerance}
\nu_{\rm loc}
=\min\left\{
\frac1{8q_{\rm loc}^6},
\frac{\delta_\varepsilon e^{-1/2}v_d\Lambda^{-d/2}}
{16Bm(2L)^d[64dq_{\rm loc}^7+8d\Lambda L^2]}
\right\}.
\end{equation}
Scaling its nodes and weights by $L$ and
taking the $d$-fold tensor product gives
the rule $(y_u,\omega_u)$ on $[-L,L]^d$
used in \eqref{eq:fa-quadrature}.

\begin{lemma}[Accuracy and cost of the local rule]
\label{lem:ap-local}
Suppose that the cubes have been accepted
and that \eqref{eq:ap-curvature}--
\eqref{eq:ap-complex-hessian} hold for
$W_{a,t}$. If $t$ satisfies
\eqref{eq:fa-ready}, then the approximation
$\widehat F_i(a)$ defined in
\eqref{eq:fa-quadrature} satisfies
\[
|\widehat F_i(a)-F_i(a)|
\le\delta_\varepsilon,
\quad
0<\widehat F_i(a)\le B,
\quad 1\le i\le n.
\]
The quadrature order and evaluation cost
satisfy \eqref{eq:fa-local-cost}.
Construction of the one-dimensional rule
requires $O(q_{\rm loc}^2
\log(e/\nu_{\rm loc}))$ operations.
\end{lemma}

\begin{proof}
Let
\[
D_\ell=z_\ell+t^{-1/2}[-L,L]^d,\quad
D=\bigcup_{\ell=1}^mD_\ell,\quad
M=\max_{1\le\ell\le m}V_a(z_\ell).
\]
By \eqref{eq:ap-center-error} and
\eqref{eq:fa-ready}, for every $\theta\in D_\ell$,
\[
\begin{aligned}
\|\theta-c_\ell\|_\infty
\le\|z_\ell-c_\ell\|+\frac L{\sqrt t}
\le\frac r2+\frac r{2\sqrt d}<2r.
\end{aligned}
\]
Thus $D_\ell\subset U_\ell$, and the cubes
$D_1,\ldots,D_m$ are disjoint.

Define
\[
\begin{aligned}
D_{\rm exact}
&=t^{d/2}\int_K e^{V_a(\theta)-M}\,d\theta,\\
A_{i,\rm exact}
&=t^{d/2}\int_K
\frac{\phi(X_i-\theta)}B
e^{V_a(\theta)-M}\,d\theta.
\end{aligned}
\]
Then $0<A_{i,\rm exact}\le D_{\rm exact}$
and
$F_i(a)=BA_{i,\rm exact}/D_{\rm exact}$.

We first obtain a lower bound for $D_{\rm exact}$.
By \eqref{eq:ap-exterior-gap}, there exists
$1\le\ell_0\le m$ such that
\[
V_a(z_{\ell_0}^\star)=\max_KV_a\ge M.
\]
Since $t\ge1/(\Lambda r^2)$ and
$\|z_{\ell_0}^\star-c_{\ell_0}\|\le r/4$,
the ball
$B(z_{\ell_0}^\star,(\Lambda t)^{-1/2})$
is contained in $U_{\ell_0}$.
By \eqref{eq:ap-curvature}, for every $\theta$
in this ball,
\[
\begin{aligned}
V_a(\theta)-M
\ge V_a(z_{\ell_0}^\star)-M
-\frac{\Lambda t}{2}
\|\theta-z_{\ell_0}^\star\|^2
\ge-\frac12.
\end{aligned}
\]
Consequently,
\begin{equation}\label{eq:ap-normalizer}
D_{\rm exact}
\ge t^{d/2}
\int_{B(z_{\ell_0}^\star,(\Lambda t)^{-1/2})}
e^{V_a-M}\,d\theta
\ge e^{-1/2}v_d\Lambda^{-d/2}.
\end{equation}
For the remainder of the proof, set
\[
c_0=e^{-1/2}v_d\Lambda^{-d/2},
\quad
\eta_{\rm loc}=\frac{c_0\delta_\varepsilon}{16B}.
\]

We next bound the integrals over
$\mathcal U\setminus D$ and $K\setminus\mathcal U$.
For $\theta\in U_\ell\setminus D_\ell$, we have
$\sqrt t\|\theta-z_\ell\|\ge L$.
By \eqref{eq:ap-center-error},
\[
\begin{aligned}
\sqrt t\|\theta-z_\ell^\star\|
\ge\sqrt t\|\theta-z_\ell\|
-\sqrt t\|z_\ell-z_\ell^\star\|
\ge L-\frac1{4\sqrt\Lambda}
\ge\frac{3L}{4},
\end{aligned}
\]
where the last inequality uses $\Lambda L^2\ge8$.
Moreover, \eqref{eq:ap-curvature} and
\eqref{eq:ap-shift-error} give
\[
\begin{aligned}
V_a(\theta)-M
\le V_a(z_\ell^\star)-M
-\frac{\kappa t}{2}
\|\theta-z_\ell^\star\|^2
\le\frac1{32}
-\frac{\kappa t}{2}
\|\theta-z_\ell^\star\|^2.
\end{aligned}
\]
Changing variables
$u=\sqrt t(\theta-z_\ell^\star)$
and summing over the disjoint cubes, we obtain
\begin{equation}
t^{d/2}\int_{\mathcal U\setminus D}
e^{V_a-M}\,d\theta
\le me^{1/32}
\int_{\|u\|\ge3L/4}
e^{-\kappa\|u\|^2/2}\,du
\le me^{1/32}(4\pi/\kappa)^{d/2}
e^{-\kappa L^2/8}
\le\frac{c_0\delta_\varepsilon}{16B}
=\eta_{\rm loc}.
\label{eq:ap-inner-tail}
\end{equation}
Here the second inequality follows from the fact that
\[
\begin{aligned}
\int_{\|u\|\ge3L/4}
e^{-\kappa\|u\|^2/2}\,du
\le e^{-9\kappa L^2/64}
\int_{\mathbb R^d}
e^{-\kappa\|u\|^2/4}\,du
\le e^{-\kappa L^2/8}(4\pi/\kappa)^{d/2},
\end{aligned}
\]
while the last inequality follows from
\eqref{eq:fa-width}.

For $\theta\in K\setminus\mathcal U$,
\eqref{eq:ap-exterior-gap} and
\eqref{eq:ap-shift-error} imply
\[
\begin{aligned}
V_a(\theta)
\le\max_\ell V_a(c_\ell)-t\gamma
\le\max_KV_a-t\gamma
\le M+\frac1{32}-t\gamma.
\end{aligned}
\]
Since $K$ has volume $(2R)^d$, we obtain
\begin{equation}\label{eq:ap-outer-tail}
\begin{aligned}
t^{d/2}\int_{K\setminus\mathcal U}
e^{V_a-M}\,d\theta
\le e^{1/32}(2R)^dt^{d/2}e^{-t\gamma}
\le\frac{c_0\delta_\varepsilon}{16B}
=\eta_{\rm loc},
\end{aligned}
\end{equation}
where the final inequality is the last condition
in \eqref{eq:fa-ready}.

Combining \eqref{eq:ap-normalizer},
\eqref{eq:ap-inner-tail}, and
\eqref{eq:ap-outer-tail}, we obtain
\[
\mu_a(K\setminus D)
=\frac{t^{d/2}\int_{K\setminus D}
e^{V_a-M}\,d\theta}{D_{\rm exact}}
\le\frac{2\eta_{\rm loc}}{c_0}
=\frac{\delta_\varepsilon}{8B}.
\]
Since $0<\phi(X_i-\theta)/B\le1$,
the discarded portions of
$A_{i,\rm exact}$ are bounded by
the same estimates.

We now consider the quadrature error on $D$.
For $1\le\ell\le m$, define
\[
\begin{aligned}
H_{\ell,0}(w)
&=e^{V_a(z_\ell+t^{-1/2}w)-M},\\
H_{\ell,i}(w)
&=\frac{\phi(X_i-z_\ell-t^{-1/2}w)}B
H_{\ell,0}(w),
\quad 1\le i\le n.
\end{aligned}
\]
After the change of variables
$\theta=z_\ell+t^{-1/2}w$, the integrals over
$D$ become
\[
\begin{aligned}
t^{d/2}\int_D e^{V_a-M}\,d\theta
&=\sum_{\ell=1}^m
\int_{[-L,L]^d}H_{\ell,0}(w)\,dw,\\
t^{d/2}\int_D
\frac{\phi(X_i-\theta)}B
e^{V_a-M}\,d\theta
&=\sum_{\ell=1}^m
\int_{[-L,L]^d}H_{\ell,i}(w)\,dw.
\end{aligned}
\]

Let $\mathcal E$ be the Bernstein ellipse
of parameter two over $[-L,L]$.
For $w\in\mathcal E^d$,
\[
\|\operatorname{Re}w\|_\infty\le\frac{5L}{4},
\quad
\|\operatorname{Im}w\|_\infty\le\frac{3L}{4},
\quad
\|w\|\le2\sqrt d\,L.
\]
By \eqref{eq:ap-center-error} and
\eqref{eq:fa-ready},
\[
\begin{aligned}
\|\operatorname{Re}(z_\ell+t^{-1/2}w)
-c_\ell\|_\infty
\le\frac r2+\frac{5L}{4\sqrt t}
\le\frac r2+\frac{5r}{8\sqrt d}<2r,\\
\|\operatorname{Im}(z_\ell+t^{-1/2}w)\|_\infty
\le\frac{3L}{4\sqrt t}
\le\frac{3r}{8\sqrt d}<\frac r{\sqrt d}.
\end{aligned}
\]
Thus $z_\ell+t^{-1/2}w\in\mathcal T_\ell$,
and the segment joining this point to
$z_\ell$ also lies in $\mathcal T_\ell$.

Taylor's formula and
\eqref{eq:ap-complex-hessian} give
\[
\begin{aligned}
V_a(z_\ell+t^{-1/2}w)-V_a(z_\ell)
=t^{-1/2}\nabla V_a(z_\ell)\cdot w
+t^{-1}\int_0^1(1-s)
w^{\mathsf T}\nabla^2V_a(z_\ell+st^{-1/2}w)w\,ds.
\end{aligned}
\]
Since $V_a(z_\ell)\le M$ and
$\|\nabla V_a(z_\ell)\|
\le\sqrt{\Lambda t}/4$
by Lemma~\ref{lem:ap-centers}, we have
\[
\begin{aligned}
\operatorname{Re}\{V_a(z_\ell+t^{-1/2}w)-M\}
\le\frac{\sqrt\Lambda}{4}\|w\|
+\frac\Lambda2\|w\|^2
\le\frac12\sqrt{d\Lambda}\,L
+2d\Lambda L^2
\le3d\Lambda L^2.
\end{aligned}
\]
Furthermore,
\[
\frac{|\phi(X_i-z_\ell-t^{-1/2}w)|}{B}
\le\exp\left\{
\frac{\|\operatorname{Im}w\|^2}{2t}
\right\}
\le e^{1/2}.
\]
Consequently, using $\Lambda L^2\ge8$,
\[
\max_{0\le i\le n}
\sup_{w\in\mathcal E^d}|H_{\ell,i}(w)|
\le e^{4d\Lambda L^2}.
\]

Let $Q_q^{(L)}$ denote the exact tensor
Gauss--Legendre rule of order $q$ on
$[-L,L]^d$. The Chebyshev approximation
argument of Lemma~\ref{lem:fp-analytic}
gives
\[
\begin{aligned}
&\left|
Q_q^{(L)}(H_{\ell,i})
-\int_{[-L,L]^d}H_{\ell,i}(w)\,dw
\right|
\le
2^{2d+1}d(2L)^d
e^{4d\Lambda L^2}2^{-2q},
\quad 0\le i\le n.
\end{aligned}
\]
For $q=q_{\rm loc}$ in \eqref{eq:fa-order},
\[
e^{4d\Lambda L^2}2^{-2q}
\le
\frac{v_d\delta_\varepsilon}
{2^{3d+7}dm e^{1/2}
(\Lambda L^2)^{d/2}B}.
\]
Therefore,
\[
\begin{aligned}
\max_{0\le i\le n}\sum_{\ell=1}^m
\left|
Q_q^{(L)}(H_{\ell,i})
-\int_{[-L,L]^d}H_{\ell,i}(w)\,dw
\right|
\le
m2^{2d+1}d(2L)^d
e^{4d\Lambda L^2}2^{-2q}
\le\frac{c_0\delta_\varepsilon}{64B}
\le\eta_{\rm loc}.
\end{aligned}
\]

It remains to bound the error in the
computed nodes and weights.
For real $w\in[-L,L]^d$, we have
$z_\ell+t^{-1/2}w\in U_\ell$, so
\eqref{eq:ap-shift-error} gives
\[
|H_{\ell,i}(w)|\le e^{1/32}<2,
\quad 0\le i\le n.
\]
Moreover, by \eqref{eq:ap-curvature} and
the gradient bound from
Lemma~\ref{lem:ap-centers},
\[
\begin{aligned}
\frac{\|\nabla V_a(z_\ell+t^{-1/2}w)\|}{\sqrt t}
\le\frac{\|\nabla V_a(z_\ell)\|}{\sqrt t}
+\Lambda\|w\|
\le\frac{\sqrt\Lambda}{4}
+\Lambda\sqrt d\,L.
\end{aligned}
\]
For $g_{\ell,i}(u)=H_{\ell,i}(Lu)$,
the product rule, the bound
$\sup_\theta\|\nabla\phi(X_i-\theta)\|\le B$,
and $L/\sqrt t\le r/(2\sqrt d)$ give
\[
\begin{aligned}
\|\nabla g_{\ell,i}(u)\|
\le2\left[
\frac r{2\sqrt d}
+\frac{\sqrt\Lambda L}{4}
+\Lambda\sqrt d\,L^2
\right]
\le8\sqrt d\,\Lambda L^2,
\quad u\in[-1,1]^d.
\end{aligned}
\]
The same estimate holds for $i=0$.
Here we used $r\le1/8$ and
$\Lambda L^2\ge8$.

Let $\widetilde Q_q^{(L)}$ be the tensor
rule constructed with tolerance
$\nu_{\rm loc}$. Applying
Lemma~\ref{lem:fp-rounded} to the functions
$g_{\ell,i}$ and rescaling to $[-L,L]^d$
gives
\[
\begin{aligned}
&|\widetilde Q_q^{(L)}(H_{\ell,i})
-Q_q^{(L)}(H_{\ell,i})|\\
&\quad\le
(2L)^d[64dq^7+8d\Lambda L^2]\nu_{\rm loc}.
\end{aligned}
\]
By \eqref{eq:ap-local-tolerance},
\[
\max_{0\le i\le n}\sum_{\ell=1}^m
|\widetilde Q_q^{(L)}(H_{\ell,i})
-Q_q^{(L)}(H_{\ell,i})|
\le\eta_{\rm loc}.
\]

The approximation in \eqref{eq:fa-quadrature}
can be written as
\[
\widehat F_i(a)=\frac{B\widehat A_i}{\widehat D},
\]
where
\[
\widehat D=\sum_{\ell=1}^m
\widetilde Q_q^{(L)}(H_{\ell,0}),
\quad
\widehat A_i=\sum_{\ell=1}^m
\widetilde Q_q^{(L)}(H_{\ell,i}).
\]
By \eqref{eq:ap-inner-tail} and
\eqref{eq:ap-outer-tail}, restricting the
integrals defining $D_{\rm exact}$ and
$A_{i,\rm exact}$ from $K$ to $D$
introduces an error of at most $2\eta_{\rm loc}$.
The analytic quadrature error and the error
in computing the nodes and weights are each
at most $\eta_{\rm loc}$. Consequently,
\[
|\widehat D-D_{\rm exact}|
\le4\eta_{\rm loc},
\quad
|\widehat A_i-A_{i,\rm exact}|
\le4\eta_{\rm loc}.
\]
By \eqref{eq:ap-normalizer},
\[
\begin{aligned}
\widehat D
\ge D_{\rm exact}-4\eta_{\rm loc}
\ge\left(1-\frac{\delta_\varepsilon}{4B}\right)
D_{\rm exact}>0.
\end{aligned}
\]
Since $A_{i,\rm exact}\le D_{\rm exact}$,
we conclude that
\[
\begin{aligned}
|\widehat F_i(a)-F_i(a)|
=\left|
\frac{B\widehat A_i}{\widehat D}
-\frac{BA_{i,\rm exact}}{D_{\rm exact}}
\right|
\le
\frac{B|\widehat A_i-A_{i,\rm exact}|}{\widehat D}
+\frac{BA_{i,\rm exact}
|\widehat D-D_{\rm exact}|}
{D_{\rm exact}\widehat D}\\
\le\frac{8B\eta_{\rm loc}}{\widehat D}
\le\frac{\delta_\varepsilon}
{2(1-\delta_\varepsilon/(4B))}
\le\delta_\varepsilon.
\end{aligned}
\]
This proves \eqref{eq:fc-evaluation-error}.
Since the quadrature weights are positive and
$0<\phi(X_i-\theta)/B\le1$ for all
$\theta\in\mathbb R^d$, we have
$0<\widehat A_i\le\widehat D$. Consequently,
$
0<\widehat F_i(a)
=\frac{B\widehat A_i}{\widehat D}\le B.
$. 

Finally, \eqref{eq:fa-width} gives
\[
\begin{aligned}
\kappa L^2
=8\max\left\{1,
\log\frac{16me^{17/32}
(4\pi\Lambda/\kappa)^{d/2}B}
{v_d\delta_\varepsilon}
\right\}
\le C_d
\log\frac{em(1+\Lambda/\kappa)B}
{\delta_\varepsilon}.
\end{aligned}
\]
Consequently, substitution in
\eqref{eq:fa-order} yields
\[
q_{\rm loc}
\le C_d\frac{\Lambda}{\kappa}
\log\frac{em(1+\Lambda/\kappa)B}
{\delta_\varepsilon}.
\]

For each quadrature node
$\theta_{\ell u}=z_\ell+t^{-1/2}y_u$,
computing the $n$ values
$\phi(X_i-\theta_{\ell u})$ and
$V_a(\theta_{\ell u})
=\sum_{i=1}^n a_i\phi(X_i-\theta_{\ell u})$
requires $O(nd)$ operations.
The resulting values are used to compute
\[
\widehat D
=\sum_{\ell,u}\omega_u
e^{V_a(\theta_{\ell u})-M},
\quad
\widehat A_i
=\frac1B\sum_{\ell,u}\omega_u
\phi(X_i-\theta_{\ell u})
e^{V_a(\theta_{\ell u})-M}.
\]
Thus computing $\widehat F_i(a)
=B\widehat A_i/\widehat D$ for all
$1\le i\le n$ requires
$O(ndm q_{\rm loc}^d)$ operations.

Similarly, each projected gradient step in
\eqref{eq:fa-ascent} requires evaluating
\[
\nabla W_{a,t}(\theta)
=\frac1t\sum_{i=1}^n
a_i(X_i-\theta)\phi(X_i-\theta),
\]
at a cost of $O(nd)$ operations.
Computing all $m$ centers therefore costs
$O(ndmJ(t))$ operations. Consequently,
the total cost of one fitted-value evaluation is
\[
O\bigl(ndm\{q_{\rm loc}^d+J(t)\}\bigr),
\]
as asserted in \eqref{eq:fa-local-cost}.
Construction of the one-dimensional quadrature
rule costs
$O(q_{\rm loc}^2\log(e/\nu_{\rm loc}))$
operations by Lemma~\ref{lem:fp-rounded}.
\end{proof}

\subsection{Proof of Theorem~\ref{thm:fa-validity}}
\begin{proof}[\textbf{Proof of Theorem~\ref{thm:fa-validity}}]

Before the restriction of the integration domain,
all fitted-value evaluations satisfy
\eqref{eq:fc-evaluation-error} by
Lemma~\ref{lem:fp-global-quadrature}.
Suppose that cubes have been accepted, that
\eqref{eq:fa-ready} holds, and that all evaluations
preceding the computation of $a^{k+1}$ are accurate.
By Lemma~\ref{lem:ap-potential}, the estimates
\eqref{eq:ap-curvature}--\eqref{eq:ap-complex-hessian}
hold for $W_{a^k,t_k}$.
Moreover, every trial vector $x^r$ used to
compute $a^{k+1}$ satisfies
\eqref{eq:fp-trial-bound}, independently of the
accuracy of its fitted-value evaluation.
Lemma~\ref{lem:ap-persistence} therefore gives
the same estimates for $W_{x^r,t_k}$, and
Lemma~\ref{lem:ap-local} implies
\[
\max_{1\le i\le n}
|\widehat F_i(x^r)-F_i(x^r)|
\le\delta_\varepsilon.
\]

Since the evaluations of fitted-values are accurate
at every inner trial, Lemma~\ref{lem:fp-inner}
implies that the accepted iterate $a^{k+1}$
satisfies
\[
B\sum_{i=1}^n
\left|
\frac{a_i^{k+1}-a_i^k}{\tau_k}
-\frac1{n\max\{F_i(a^{k+1}),b\}}
\right|
\le\eta_\varepsilon.
\]
Repeating this argument proves that all subsequent
fitted-value evaluations satisfy
\eqref{eq:fc-evaluation-error}, and that every
accepted proximal step satisfies the inexact
update estimate of Lemma~\ref{lem:fp-inner}.
Consequently, the accuracy and iteration bounds
of Theorem~\ref{thm:fc-iterations} also hold for
Algorithm~\ref{alg:fa-adaptive}.

We next bound the time at which local quadrature
begins. For fixed data satisfying
\eqref{eq:fa-morse}, Lemma~\ref{lem:ap-selection}
gives
\[
t_{\rm sel}\le T_0,\quad
m+r^{-1}+\kappa^{-1}+\gamma^{-1}+\Lambda
\le C_{X,d},
\]
for all sufficiently small $\varepsilon$.
Since
$\log(B/\delta_\varepsilon)
=O_{X,d}(\log(e/\varepsilon))$,
the definition \eqref{eq:fa-width} gives
$
L^2=O_{X,d}(\log(e/\varepsilon)).
$

For every $t>0$, we have
$\log t\le\gamma t/d-1+\log(d/\gamma)$,
and consequently
\[
t\gamma-\frac d2\log t
\ge\frac{\gamma t}{2}
+\frac d2-\frac d2\log(d/\gamma).
\]
Thus all the inequalities in
\eqref{eq:fa-ready} hold whenever
\begin{equation}\label{eq:ap-ready-time}
\begin{gathered}
t\ge\max\left\{
t_{\rm sel},1,\frac1{\Lambda r^2},
\frac{4dL^2}{r^2},\frac d{2\gamma}
\right\},\\
t\ge\frac2\gamma\left[
\left(
\log\frac{16e^{17/32}(2R)^d
\Lambda^{d/2}B}{v_d\delta_\varepsilon}
\right)_+
+\frac d2\max\{0,\log(d/\gamma)\}
\right].
\end{gathered}
\end{equation}
Both right-hand sides are
$O_{X,d}(\log(e/\varepsilon))$.
Let $T$ be the maximum of the lower bounds in
\eqref{eq:ap-ready-time}, and choose
$j\ge1$ such that $T\le2^j\le2T$.
At the first completed iterate $a^k$ satisfying
$t_k\ge2^j$, we have
\[
2^j\le t_k<2^j+\frac1{64}
\le2T+\frac1{64},
\]
since $t_k-t_{k-1}=\tau_{k-1}\le1/64$.
Thus \eqref{eq:fa-ready} holds at this checkpoint,
and the algorithm switches to local quadrature.
Moreover,
\[
\frac{32H_\varepsilon}{\varepsilon}
\bigg/\log(e/\varepsilon)\longrightarrow\infty
\quad(\varepsilon\downarrow0).
\]
Therefore, for all sufficiently small
$\varepsilon$, the algorithm switches to local
quadrature at a time
$O_{X,d}(\log(e/\varepsilon))$.

Finally, the function
$t\mapsto t\gamma-(d/2)\log t$ is increasing
for $t\ge d/(2\gamma)$, and the remaining
inequalities in \eqref{eq:fa-ready} are
preserved as $t$ increases. Thus
\eqref{eq:fa-ready} holds at every subsequent
evaluation, completing the proof.
\end{proof}

\subsection{Proof of Theorem~\ref{thm:fa-total-cost}}

\begin{proof}[\textbf{Proof of Theorem~\ref{thm:fa-total-cost}}]
Let $\mathcal V$ denote the set of pairs $(j,h)$
examined by \textsc{SelectCubes}, and let $m_{jh}$
be the number of centers retained at each such pair.
The partition at level $h$ contains $2^{hd}$ cells.
Computing $W$, $\nabla W$, and $\nabla^2W$ at their
midpoints requires $O(nd^2 2^{hd})$ operations.
Sorting these values, comparing the distances
between candidate centers, and checking the cell
intersections require at most
$O((hd+m_{jh}d)2^{hd})$ operations.
Finally, for each $0\le p\le j+1$, checking
the inequalities in \eqref{eq:fa-accept} requires
$O((d^3+m_{jh}d)2^{hd})$ operations.
Thus the total cost of selecting the cubes is
bounded by
\begin{equation}\label{eq:ap-search-cost}
C\sum_{(j,h)\in\mathcal V}
\left[nd^2+(j+2)(d^3+m_{jh}d)+hd\right]2^{hd}.
\end{equation}

Since $h\le\max\{1,\lfloor j/4\rfloor\}$,
we have
\[
2^{hd}\le2^d(2^j)^{d/4},
\quad
m_{jh}\le2^{hd}.
\]
Moreover, the checkpoints examined before
local quadrature begins satisfy
$2^j\le t_{\rm loc}$.
Summing \eqref{eq:ap-search-cost} over these
checkpoints and partition levels therefore gives
\begin{equation}\label{eq:ap-search-bound}
\widetilde O_d\left(
n(1+t_{\rm loc})^{d/4}
+(1+t_{\rm loc})^{d/2}
\right).
\end{equation}

We next consider the proximal steps performed
using quadrature on $K$.
The fixed grid $\mathcal G$ contains
$O_d(R^d)$ points. Its construction requires
$O_d(ndR^d)$ operations, and evaluating the
score on this grid at each proximal step costs
$O_d(nR^d)$ operations. The total cost of
these evaluations is therefore
$\widetilde O_d(nR^dN)$.

By \eqref{eq:fc-global-order}, the order of
the global quadrature rule at any time
$t\le t_{\rm loc}$ satisfies
\[
q_K\le
C_dR\sqrt{1+\log R+t_{\rm loc}}
\log\frac{eR(1+\log R+t_{\rm loc})B}
{\delta_\varepsilon}.
\]
By Lemma~\ref{lem:fp-inner}, each proximal
step requires at most
\[
2+\left\lceil
\frac{\log(8B/(b\eta_\varepsilon))}{\log8}
\right\rceil
\]
fitted-value evaluations. Since $G$ steps
are computed using global quadrature,
including the evaluations at their initial
and terminal coefficient vectors gives at most
\[
(G+1)\left[
2+\left\lceil
\frac{\log(8B/(b\eta_\varepsilon))}{\log8}
\right\rceil\right]
\]
evaluations using these rules.

Each evaluation requires computing
$\phi(X_i-\theta)$ for all $1\le i\le n$
at the $q_K^d$ tensor-product nodes,
followed by the evaluation of $V_a(\theta)$
and the fitted-value ratios. This costs
$O(ndq_K^d)$ operations. Consequently,
the total cost of global quadrature evaluations is
\[
\widetilde O_d\left(
nR^d(G+1)(1+\log R+t_{\rm loc})^{d/2}
\right).
\]
The construction of the one-dimensional
Gauss--Legendre rules at the dyadic times
requires an additional
\[
\widetilde O_d\left(
R^2(1+\log R+t_{\rm loc})\right)
\]
operations, by Lemma~\ref{lem:fp-rounded}.
Since $R\ge1$, the bound in
\eqref{eq:ap-search-bound} is absorbed by
the cost of the global quadrature evaluations.

We now consider the remaining $N-G$ proximal
steps, which use the local quadrature rule.
By \eqref{eq:fa-local-cost} and
Lemma~\ref{lem:ap-local}, each fitted-value
evaluation costs
\[
O\left(ndm\{q_{\rm loc}^d+J(t)\}\right),
\]
where
\[
q_{\rm loc}
=\widetilde O_d(\Lambda/\kappa),
\quad
J(t)=\widetilde O_d(\Lambda/\kappa).
\]
Since each proximal step requires
$O(\log(8B/(b\eta_\varepsilon)))$
fitted-value evaluations, and
$d\ge1$, $\Lambda/\kappa\ge1$, the
total cost of these steps is
\[
\widetilde O_d\left(
nm(N-G)(\Lambda/\kappa)^d\right).
\]

The one-dimensional local quadrature rule
is constructed once, at a cost of
$O(q_{\rm loc}^2\log(e/\nu_{\rm loc}))$
operations. By \eqref{eq:ap-local-tolerance},
the logarithm depends on $L$ and $\Lambda$
through the product
\[
(2L)^d\Lambda^{d/2}
=2^d(\Lambda L^2)^{d/2}.
\]
Since \eqref{eq:fa-width} gives
\[
\Lambda L^2
=O_d\left(
\frac{\Lambda}{\kappa}
\log\frac{em(1+\Lambda/\kappa)B}
{\delta_\varepsilon}
\right),
\]
we obtain a construction cost of
$\widetilde O_d((\Lambda/\kappa)^2)$.

The output correction in Algorithm~2
requires $O(nd)$ operations. Combining
these estimates gives
\[
\begin{aligned}
\widetilde O_d\bigg(
nR^dN
+nR^dG(1+\log R+t_{\rm loc})^{d/2}
+R^2(1+\log R+t_{\rm loc})\\
+nm(N-G)(\Lambda/\kappa)^d
+(\Lambda/\kappa)^2
\bigg),
\end{aligned}
\]
which proves \eqref{eq:fa-total-cost}.

Finally, suppose that the data satisfy
\eqref{eq:fa-morse}. By
Theorems~\ref{thm:fa-validity} and
\ref{thm:fc-iterations}, for sufficiently
small $\varepsilon$,
\[
t_{\rm loc},G
=O_{X,d}(\log(e/\varepsilon)),
\quad
N=O_{X,d}(\varepsilon^{-1}\log(e/\varepsilon)).
\]
Moreover, $m$, $\kappa^{-1}$, and $\Lambda$
are bounded independently of $\varepsilon$.
It follows from the definitions of the
global and local quadrature orders that
\[
q_K=O_{X,d}(\log^{3/2}(e/\varepsilon)),
\quad
q_{\rm loc},J(t_N)
=O_{X,d}(\log(e/\varepsilon)).
\]
The number of fitted-value evaluations per
proximal step is also
$O_{X,d}(\log(e/\varepsilon))$.
Consequently, the global and local
quadrature evaluations require, respectively,
\[
\begin{aligned}
&O_{X,d}\left(
n\log^{2+3d/2}(e/\varepsilon)\right),\\
&O_{X,d}\left(
nm_\varepsilon\varepsilon^{-1}
\log^{d+2}(e/\varepsilon)\right)
\end{aligned}
\]
operations.

The construction of the global quadrature
rules costs
$O_{X,d}(\log^4(e/\varepsilon))$, while
the local rule is constructed in
$O_{X,d}(\log^3(e/\varepsilon))$
operations. By Lemma~\ref{lem:ap-selection},
the selection procedure examines only
finitely many partitions, independently of
$\varepsilon$. These costs are absorbed by
the two preceding bounds, proving
\eqref{eq:fa-fixed-data}.
\end{proof}

\section{Proofs of Results in Section \ref{sec:generic}}\label{sec:generic proofs}
\subsection{Auxiliary Lemmas}
\begin{lemma}
\label{lem:independence}
Let $z_1,\ldots,z_r\in\R^d$ be distinct, and let $P_1,\ldots,P_r$
be real polynomials. If
\begin{equation}
\label{eq:independence}
    \sum_{j=1}^r P_j(x)e^{\langle x,z_j\rangle}=0
                   \quad(x\in\R^d),
\end{equation}
then every $P_j$ is identically zero.
\end{lemma}

\begin{proof}
Choose $u\in\mathbb S^{d-1}$ such that 
$\lambda_j:=\langle u,z_j\rangle$ are pairwise distinct for all
$1\le j\le r$.
Next, fix $x_0\in\mathbb R^d$ and set
$
Q_j(t):=e^{\langle x_0,z_j\rangle}P_j(x_0+tu).
$
Then, restricting \eqref{eq:independence} to $x=x_0+tu$ yields
\[
\sum_{j=1}^r Q_j(t)e^{\lambda_jt}=0,
\quad t\in\mathbb R.
\]
After relabeling, we may assume $\lambda_1<\cdots<\lambda_r$.
Then
\[
Q_r(t)
=-\sum_{j=1}^{r-1}
Q_j(t)e^{-(\lambda_r-\lambda_j)t}
\longrightarrow0
\quad\text{as }t\to\infty.
\]
Since $Q_r$ is a polynomial, it follows that
$Q_r\equiv0$. Repeating the argument gives
$Q_j\equiv0$ for every $1\le j\le r$. In particular,
$
0=Q_j(0)
=e^{\langle x_0,z_j\rangle}P_j(x_0).
$
Since $x_0\in\mathbb R^d$ was arbitrary,
$P_j\equiv0$ for every $j$.
\end{proof}

\begin{lemma}
\label{lem:jet}
Let $f,g\in C^\infty(\mathbb R^d)$, with $g>0$.
Suppose that there exist $r \geq 1$, $z_1,\ldots,z_r\in\mathbb R^d$
and nonnegative integers $q_1,\ldots,q_r$ such that,
for every $c\in\mathbb R$,
\begin{equation}\label{eq:numerator}
g(x)(f(x)-c)
=\sum_{j=1}^r P_{j,c}(x)e^{\langle x,z_j\rangle},
\quad \deg P_{j,c}\le q_j,
\end{equation}
for some real polynomials $P_{1,c},\ldots,P_{r,c}$,
where the right-hand side is not identically zero.
Set $ J=\sum_{j=1}^r(q_j+1) \geq 2$.
Then, for every $x_0\in\mathbb R^d$, there exist
$u\in\mathbb S^{d-1}$ and $1\le\ell\le J-1$ such that
\[
D_u^\ell f(x_0)\ne0.
\]
\end{lemma}

\begin{proof}
Fix $x_0\in\mathbb R^d$ and set $c=f(x_0)$. By assumption,
$H=g(f-c)$ is not identically zero. Choose $x_1$ such that
$H(x_1)\ne0$, and set
\[
u=\frac{x_1-x_0}{\|x_1-x_0\|},
\quad h(t)=H(x_0+tu).
\]
Then $h(0)=0$ and $h(\|x_1-x_0\|)\ne0$. Also we have,
$
h(t)=\sum_{j=1}^r
P_{j,c}(x_0+tu)e^{\langle x_0,z_j\rangle}
e^{t\langle u,z_j\rangle}.
$
Since $t\mapsto P_{j,c}(x_0+tu)$ has degree at most $q_j$,
\[
\prod_{j=1}^r
\left(\frac{d}{dt}-\langle u,z_j\rangle\right)^{q_j+1}
h(t)=0.
\]
This is a homogeneous linear ODE of order $J$. By uniqueness
for its initial-value problem, and since $h$ is not identically
zero, there exists $1\le\ell\le J-1$ such that
$
h^{(\ell)}(0)\ne0.
$

Now, suppose instead that
\[
D_u^m f(x_0)=0,\quad 1\le m\le J-1.
\]
Since $f(x_0)-c=0$, Leibniz's rule applied to
\[
h(t)=g(x_0+tu)(f(x_0+tu)-c)
\]
gives
\[
h^{(m)}(0)=0,\quad 0\le m\le J-1,
\]
a contradiction. Hence $D_u^\ell f(x_0)\ne0$ for some
$1\le\ell\le J-1$.
\end{proof}

We next establish a version of the van der Corput
estimate for oscillatory integrals
\cite{stein1993harmonic}.

\begin{lemma}
\label{lem:vdc}
Let $I\subset\mathbb R$ be a bounded interval, $q\ge2$,
and $f\in C^q(I)$ satisfy
\[
|f^{(q)}(t)|\ge c>0,\quad t\in I.
\]
Then, for every $A\in C^1(I)$ and $\lambda>0$,
\begin{equation}\label{eq:vdc}
\left|\int_I e^{i\lambda f(t)}A(t)\,dt\right|
\le C_q(\lambda c)^{-1/q}
\left(\|A\|_\infty+\int_I|A'(t)|\,dt\right).
\end{equation}
If $f\in C^2(I)$ and satisfies
\[
|f'(t)|\ge c>0,\quad
|f''(t)|\le M,\quad t\in I,
\]
then
\begin{equation}\label{eq:vdc-one}
\left|\int_I e^{i\lambda f(t)}A(t)\,dt\right|
\le\frac1\lambda
\left[
\frac{2\|A\|_\infty+\int_I|A'(t)|\,dt}{c}
+\frac{M|I|\|A\|_\infty}{c^2}
\right].
\end{equation}
\end{lemma}

\begin{proof}
Let $p\ge1$,
$g\in C^p(I)$, and suppose that $|g^{(p)}|\ge c>0$ on $I$.
We claim that, for every $\eta>0$,
\begin{equation}\label{eq:sublevel}
|\{t\in I:|g(t)|\le\eta\}|
\le 2(p+1)\left(\frac{\eta}{c}\right)^{1/p}.
\end{equation}
Set $E=\{t\in I:|g(t)|\le\eta\}$,
and $L=|E|$ We may assume $L>0$. Define
\[
F(t)=|(E\cap(-\infty,t])|,\quad t\in\mathbb R.
\]
For $s<t$,
\[
0\le F(t)-F(s)
=|E\cap(s,t]|\le t-s.
\]
Thus $F$ is continuous. For $0\le j\le p$, let
\[
t_j=\inf\left\{t\in\mathbb R:
F(t)\ge\frac{j+1/2}{p+1}L\right\}.
\]
Since $F$ is continuous, $F(t_j)=(j+1/2)L/(p+1)$.
Moreover, $t_j\in E$. Indeed, if $t_j\notin E$,
continuity of $g$ implies that $E$ is empty in a
neighborhood of $t_j$, so $F$ is constant there,
contradicting the definition of $t_j$. Now, for $0\le i<j\le p$, we have
\[
\begin{aligned}
t_j-t_i
\ge|E\cap(t_i,t_j]|
=F(t_j)-F(t_i)
=\frac{j-i}{p+1}L.
\end{aligned}
\]
Applying the mean-value theorem for divided differences
to $g$ at $t_0,\ldots,t_p$, there exists
$\xi\in(t_0,t_p)$ such that
\[
\frac{g^{(p)}(\xi)}{p!}
=
\sum_{j=0}^p
\frac{g(t_j)}{\prod_{h\ne j}(t_j-t_h)}.
\]
Since $|g(t_j)|\le\eta$ and
\[
\prod_{h\ne j}|t_j-t_h|
\ge\left(\frac{L}{p+1}\right)^p j!(p-j)!,
\]
we obtain
\[
\begin{aligned}
\frac{c}{p!}
\le
\left|\sum_{j=0}^p
\frac{g(t_j)}{\prod_{h\ne j}(t_j-t_h)}\right|
\le
\eta\left(\frac{p+1}{L}\right)^p
\sum_{j=0}^p\frac1{j!(p-j)!}
=\frac{2^p\eta}{p!}
\left(\frac{p+1}{L}\right)^p.
\end{aligned}
\]
We can rearrange this to prove \eqref{eq:sublevel}.

Now suppose $q\ge2$. Applying \eqref{eq:sublevel}
to $g=f'$ with $p=q-1$ gives
\begin{equation}\label{eq:vdc-small}
|\{t\in I:|f'(t)|\le\eta\}|
\le 2q\left(\frac{\eta}{c}\right)^{1/(q-1)}.
\end{equation}
Since $f^{(q)}$ does not vanish, repeated application
of Rolle's theorem shows that $f''$ has at most $q-2$
distinct zeros. Consequently, $I$ can be partitioned
into at most $q-1$ subintervals on which $f'$ is monotone.
Intersecting these intervals with $\{|f'|>\eta\}$
gives at most $2(q-1)$ intervals on which $f'$ is
monotone and $|f'|\ge\eta$.

Let $J$ be one such interval.
Integrating by parts, we get
\[
\begin{aligned}
\int_J e^{i\lambda f(t)}\,dt
=
\int_J\frac1{i\lambda f'(t)}
\frac{d}{dt}e^{i\lambda f(t)}\,dt
=
\left[\frac{e^{i\lambda f(t)}}{i\lambda f'(t)}
\right]_{\partial J}
+\frac1{i\lambda}\int_J
e^{i\lambda f(t)}
\frac{f''(t)}{f'(t)^2}\,dt.
\end{aligned}
\]
As $f'$ is monotone and does not vanish on $J$,
\[
\begin{aligned}
\int_J\frac{|f''(t)|}{|f'(t)|^2}\,dt
=\left|\left[\frac1{f'(t)}\right]_{\partial J}\right|
\le\frac1\eta.
\end{aligned}
\]
Therefore,
$
\left|\int_J e^{i\lambda f(t)}\,dt\right|
\le\frac3{\lambda\eta}.
$
Splitting $I$ into the sets where $|f'|\le\eta$
and $|f'|>\eta$, and using \eqref{eq:vdc-small},
we obtain
\[
\left|\int_I e^{i\lambda f(t)}\,dt\right|
\le
2q\left(\frac{\eta}{c}\right)^{1/(q-1)}
+\frac{6(q-1)}{\lambda\eta}.
\]
Taking
$
\eta=c^{1/q}\lambda^{-(q-1)/q}
$
gives
\begin{equation}\label{eq:vdc-unweighted}
\left|\int_I e^{i\lambda f(t)}\,dt\right|
\le C_q(\lambda c)^{-1/q}.
\end{equation}
The same argument applies to every subinterval of $I$.

Let $a=\inf I$ and define
$
H(t)=\int_a^t e^{i\lambda f(s)}\,ds.
$
By \eqref{eq:vdc-unweighted},
$
\|H\|_\infty\le C_q(\lambda c)^{-1/q}.
$
Integration by parts now gives us
\[
\begin{aligned}
\int_I e^{i\lambda f(t)}A(t)\,dt
=\int_I H'(t)A(t)\,dt
=[H(t)A(t)]_{\partial I}
-\int_I H(t)A'(t)\,dt.
\end{aligned}
\]
Consequently,
\[
\left|\int_I e^{i\lambda f(t)}A(t)\,dt\right|
\le C_q(\lambda c)^{-1/q}
\left(\|A\|_\infty+\int_I|A'(t)|\,dt\right),
\]
which proves \eqref{eq:vdc}.

Finally, suppose $f\in C^2(I)$,
$|f'|\ge c>0$, and $|f''|\le M$. We integrate by parts again to get
\[
\begin{aligned}
\int_I e^{i\lambda f(t)}A(t)\,dt
=
\left[\frac{e^{i\lambda f(t)}A(t)}
{i\lambda f'(t)}\right]_{\partial I}
-\frac1{i\lambda}
\int_I e^{i\lambda f(t)}
\left(
\frac{A'(t)}{f'(t)}
-\frac{A(t)f''(t)}{f'(t)^2}
\right)\,dt.
\end{aligned}
\]
Hence,
\[
\left|\int_I e^{i\lambda f(t)}A(t)\,dt\right|
\le\frac1\lambda\left[
\frac{2\|A\|_\infty+\int_I|A'(t)|\,dt}{c}
+\frac{M|I|\|A\|_\infty}{c^2}
\right],
\]
proving \eqref{eq:vdc-one}.
\end{proof}

\begin{lemma}[Uniform Fourier decay]\label{lem:fourier}
Let $C\subset\mathbb R^p$ be compact, and suppose that
$\gamma_\eta:\mathbb R^d\to\mathbb R^m$ is jointly smooth
in $(x,\eta)$ on a neighborhood of $\mathbb R^d\times C$.
Suppose that, for every $\eta\in C$, $a\in\mathbb S^{m-1}$,
and $c\in\mathbb R$, the function $a\cdot\gamma_\eta-c$
admits a representation of the form \eqref{eq:numerator}
with $g>0$,
\[
\sum_{j=1}^r(q_j+1)\le J,
\]
for some fixed integer $J\ge2$, and
$a\cdot\gamma_\eta-c$ is not identically zero.
Then, for every $\rho\in C_c^\infty(\mathbb R^d)$,
\begin{equation}\label{eq:fourier}
\sup_{\eta\in C}
\left|\int_{\mathbb R^d}
e^{i\langle\xi,\gamma_\eta(x)\rangle}
\rho(x)\,dx\right|
\le C_{C,\rho}(1+\|\xi\|)^{-1/(J-1)},
\quad \xi\in\mathbb R^m.
\end{equation}
\end{lemma}

\begin{proof}
Let
$f_{\eta,a}(x)=\langle a,\gamma_\eta(x)\rangle$
and choose a compact cube $K$ such that
$\operatorname{supp}\rho\subset\operatorname{int}K$.
Fix $(x_0,\eta_0,a_0)\in K\times C\times\mathbb S^{m-1}$.
Applying Lemma~\ref{lem:jet} to $f_{\eta_0,a_0}$ with
$c=f_{\eta_0,a_0}(x_0)$, there exist
$u\in\mathbb S^{d-1}$ and $1\le\ell\le J-1$ such that
\[
D_u^\ell f_{\eta_0,a_0}(x_0)\ne0.
\]
By joint smoothness, there exist an open rectangular box
$Q$ containing $x_0$, with a coordinate direction $u$,
and an open neighborhood $V$ of $(\eta_0,a_0)$ such that
\[
|D_u^\ell f_{\eta,a}(x)|\ge c>0,
\quad (x,\eta,a)\in Q\times V.
\]
We may choose $Q$ and $V$ such that their (compact) closures contained
in the domain of smoothness.

By compactness, finitely many such neighborhoods
$Q_h\times V_h$, $1\le h\le H$, cover
$K\times C\times\mathbb S^{m-1}$. Denote their corresponding
directions, derivative orders, and lower bounds by
$u_h,\ell_h,c_h$, respectively. Thus
\begin{equation}\label{eq:fourier-local}
|D_{u_h}^{\ell_h}f_{\eta,a}(x)|\ge c_h,
\quad (x,\eta,a)\in Q_h\times V_h,
\end{equation}
where $1\le\ell_h\le J-1$.

Choose nonnegative smooth functions
$b_h\in C_c^\infty(Q_h\times V_h)$ such that
$\sum_{h=1}^H b_h>0$ on a neighborhood of
$K\times C\times\mathbb S^{m-1}$. Let $\zeta$ be a smooth
compactly supported function, equal to one on this compact
set and supported where $\sum_h b_h>0$. Define
\[
\chi_h=\frac{\zeta b_h}{\sum_{j=1}^H b_j},
\]
extended by zero outside the support of $\zeta$.
Then
\begin{equation}\label{eq:partition-bounds}
\sum_{h=1}^H\chi_h=1
\quad\text{on }K\times C\times\mathbb S^{m-1},
\quad
\sup_{h,x,\eta,a}
\left(|\chi_h|+\|\nabla_x\chi_h\|\right)<\infty.
\end{equation}
Moreover, $\operatorname{supp}\chi_h$ is a compact subset
of $Q_h\times V_h$.

Let $\xi=\lambda a$, where $\lambda=\|\xi\|\ge1$ and
$a\in\mathbb S^{m-1}$. By \eqref{eq:partition-bounds},
\[
\int_{\mathbb R^d}
e^{i\langle\xi,\gamma_\eta(x)\rangle}\rho(x)\,dx
=
\sum_{h=1}^H
\int_{Q_h}e^{i\lambda f_{\eta,a}(x)}
\rho(x)\chi_h(x,\eta,a)\,dx.
\]

Next, we fix $h$ and choose an orthonormal basis of $\mathbb R^d$
containing $u_h$. Since $Q_h$ was chosen to be a
rectangular box in these coordinates, we may write
$
Q_h=\{y+tu_h:y\in P_h,\ t\in I_h\},
$
where $P_h\subset u_h^\perp$ is a bounded rectangular box
and $I_h\subset\mathbb R$ is a bounded interval.
For fixed $y,\eta,a$, set
\[
A(t)=\rho(y+tu_h)\chi_h(y+tu_h,\eta,a),
\quad t\in I_h.
\]
Since $\operatorname{supp}\chi_h$ is a compact subset of
$Q_h\times V_h$, the function $A$ is supported in the
interior of $I_h$.
Furthermore,
\[
\begin{aligned}
|A(t)|
&\le \|\rho\|_\infty|\chi_h(y+tu_h,\eta,a)|, \text{ and }\\
|A'(t)|
&\le \|\nabla\rho\|_\infty|\chi_h(y+tu_h,\eta,a)| +\|\rho\|_\infty
\|\nabla_x\chi_h(y+tu_h,\eta,a)\|.
\end{aligned}
\]
Consequently, \eqref{eq:partition-bounds} gives
\begin{equation}\label{eq:fourier-amplitude}
\sup_{y,\eta,a}
\left(\|A\|_\infty+\int_{I_h}|A'(t)|\,dt\right)
\le M_h<\infty.
\end{equation}

Suppose first that $\ell_h\ge2$. For every $(\eta,a)\in V_h$,
\eqref{eq:fourier-local} gives
\[
\left|\frac{d^{\ell_h}}{dt^{\ell_h}}
f_{\eta,a}(y+tu_h)\right|
\ge c_h,\quad t\in I_h.
\]
Applying Lemma~\ref{lem:vdc} and
\eqref{eq:fourier-amplitude}, we obtain
\[
\left|
\int_{I_h}e^{i\lambda f_{\eta,a}(y+tu_h)}
A(t)\,dt
\right|
\le C_{\ell_h}M_h(\lambda c_h)^{-1/\ell_h}.
\]

If $\ell_h=1$, joint smoothness and compactness of
$\overline{Q_h}\times\overline{V_h}$ give
$
\sup_{(x,\eta,a)\in Q_h\times V_h}
|D_{u_h}^2f_{\eta,a}(x)|\le M_h'<\infty.
$

Applying \eqref{eq:vdc-one} with
$|D_{u_h}f_{\eta,a}|\ge c_h$, we obtain
\[
\left|
\int_{I_h}e^{i\lambda f_{\eta,a}(y+tu_h)}
A(t)\,dt
\right|
\le C_h\lambda^{-1}.
\]
If $(\eta,a)\notin V_h$, then
$\chi_h(x,\eta,a)=0$ for every $x\in\mathbb R^d$,
and hence $A\equiv0$. Thus, both estimates hold
uniformly for all $\eta\in C$ and $a\in\mathbb S^{m-1}$.

Since $\ell_h\le J-1$ and $\lambda\ge1$, it follows
in either case that
\[
\left|
\int_{I_h}e^{i\lambda f_{\eta,a}(y+tu_h)}
A(t)\,dt
\right|
\le C_h\lambda^{-1/(J-1)}.
\]
Integrating over $P_h$ and applying Fubini's theorem gives
\[
\left|
\int_{Q_h}e^{i\lambda f_{\eta,a}(x)}
\rho(x)\chi_h(x,\eta,a)\,dx
\right|
\le C_h|P_h|\lambda^{-1/(J-1)}.
\]
Here $|P_h|=1$ when $d=1$.

Summing over $h=1,\ldots,H$, we conclude that
\[
\sup_{\eta\in C}
\left|
\int_{\mathbb R^d}
e^{i\langle\xi,\gamma_\eta(x)\rangle}
\rho(x)\,dx
\right|
\le C_{C,\rho}\|\xi\|^{-1/(J-1)},
\quad \|\xi\|\ge1.
\]
For $\|\xi\|<1$, the integral is bounded by
$\|\rho\|_{L^1}$. Enlarging $C_{C,\rho}$ proves
\eqref{eq:fourier}.
\end{proof}

\begin{lemma}[Uniform small-ball bound]\label{lem:smallball}
Under the hypotheses of Lemma~\ref{lem:fourier}, let
$Y_1,\ldots,Y_n$ be independent random vectors in
$\mathbb R^d$ with smooth, compactly supported densities.
If $n>m(J-1)$, then there exists $C'>0$ such that
\begin{equation}\label{eq:smallball}
\sup_{\eta\in C}\sup_{b\in\mathbb R^m}
\Pr\left(
\left\|\frac1n\sum_{i=1}^n\gamma_\eta(Y_i)-b\right\|
\le\delta
\right)
\le C'\delta^m,
\quad 0<\delta\le1.
\end{equation}
The same conclusion holds if the densities of
$Y_1,\ldots,Y_n$ are merely bounded and compactly supported.
\end{lemma}

\begin{proof}
Let
\[
S_\eta=\frac1n\sum_{i=1}^n\gamma_\eta(Y_i),
\]
and denote its characteristic function by
$\widehat\nu_\eta$. By independence and
Lemma~\ref{lem:fourier}, there exists $C>0$,
independent of $\eta$, such that
\[
\begin{aligned}
|\widehat\nu_\eta(\xi)|
=\prod_{i=1}^n
\left|\mathbb E
e^{i\langle\xi/n,\gamma_\eta(Y_i)\rangle}\right|
\le C\left(1+\frac{\|\xi\|}{n}\right)^{-n/(J-1)}.
\end{aligned}
\]
Since $n/(J-1)>m$, the right-hand side is
integrable over $\mathbb R^m$. Consequently,
$
\sup_{\eta\in C}
\|\widehat\nu_\eta\|_{L^1(\mathbb R^m)}
<\infty.
$
By Fourier inversion, $S_\eta$ has a bounded
continuous density $h_\eta$, given by
\[
h_\eta(z)
=\frac1{(2\pi)^m}
\int_{\mathbb R^m}
e^{-i\langle\xi,z\rangle}
\widehat\nu_\eta(\xi)\,d\xi.
\]
In particular,
$
\sup_{\eta\in C}\|h_\eta\|_\infty
\le (2\pi)^{-m}
\sup_{\eta\in C}
\|\widehat\nu_\eta\|_{L^1}
<\infty.
$
Thus, for every $b\in\mathbb R^m$ and $\delta>0$,
\[
\begin{aligned}
\Pr\{\|S_\eta-b\|\le\delta\}
=\int_{B(b,\delta)}h_\eta(z)\,dz
\le v_m\delta^m\|h_\eta\|_\infty,
\end{aligned}
\]
where $v_m$ denotes the volume of the unit ball
in $\mathbb R^m$. Taking the supremum over $\eta$
and $b$ proves \eqref{eq:smallball}.

Finally, suppose that $Y_i$ has a bounded density
$p_i$ supported on a compact set $K_i$. Choose
$\rho_i\in C_c^\infty(\mathbb R^d)$ such that
$\rho_i\ge0$, $\int\rho_i=1$, and
$\inf_{x\in K_i}\rho_i(x)>0$. Then there exists
$M_i<\infty$ such that
\[
p_i(x)\le M_i\rho_i(x),
\quad x\in\mathbb R^d.
\]
Let $Z_1,\ldots,Z_n$ be independent random vectors
with respective densities $\rho_1,\ldots,\rho_n$.
For every Borel set $A\subset\mathbb R^m$,
independence gives
\[
\begin{aligned}
\Pr\left\{
\frac1n\sum_{i=1}^n\gamma_\eta(Y_i)\in A
\right\}
=
\int_{(\mathbb R^d)^n}
\mathbf1_{\{
n^{-1}\sum_i\gamma_\eta(x_i)\in A
\}}
\prod_{i=1}^n p_i(x_i)\,dx_1\cdots dx_n
\\
\le 
\left(\prod_{i=1}^nM_i\right)
\Pr\left\{
\frac1n\sum_{i=1}^n\gamma_\eta(Z_i)\in A
\right\}.
\end{aligned}
\]
Applying the first part with $A=B(b,\delta)$
proves the result.
\end{proof}

\begin{lemma}\label{lem:avoidance}
Let $C$ be a compact parameter set, and suppose that,
for every $0<\delta\le1$, there exists a $\delta$-net
of $C$ with cardinality at most $C_0\delta^{-p}$,
where $p<m$. Let
$\gamma_\eta:\mathbb R^d\to\mathbb R^m$ satisfy the
hypotheses of Lemma~\ref{lem:smallball}.
Suppose also that $\gamma_\eta(x)$ is uniformly locally
Lipschitz in $\eta$ for $x$ in compact subsets of
$\mathbb R^d$.
Then, for every fixed $b\in\mathbb R^m$, the set
\begin{equation}\label{eq:bad-datasets}
\left\{
X\in(\mathbb R^d)^n:
\exists\eta\in C,\quad
\frac1n\sum_{i=1}^n\gamma_\eta(X_i)=b
\right\}
\end{equation}
has Lebesgue measure zero.
\end{lemma}

\begin{proof}
Fix $M\ge1$ and let $K_M=[-M,M]^d$. Choose a
smooth compactly supported probability density $\rho_M$
such that
\[
\inf_{x\in K_M}\rho_M(x)>0,
\]
and let $Y_1,\ldots,Y_n$ be independent random vectors
with density $\rho_M$.

By the assumption of uniform local Lipschitzness and the compactness
of $K_M\times C$, there exist $L_M\ge1$ and $\delta_M>0$
such that
\[
\|\gamma_\eta(x)-\gamma_{\eta'}(x)\|
\le L_M\|\eta-\eta'\|,
\quad x\in K_M,
\]
whenever $\eta,\eta'\in C$ and
$\|\eta-\eta'\|\le\delta_M$.

For $0<\delta\le\min\{\delta_M,L_M^{-1}\}$,
let $C_\delta$ be a $\delta$-net of $C$ satisfying
$
|C_\delta|\le C_0\delta^{-p}.
$
If $Y_i\in K_M$ for every $i$ and
$
n^{-1}\sum_{i=1}^n\gamma_\eta(Y_i)=b
$
for some $\eta\in C$, choose $\eta'\in C_\delta$
such that $\|\eta-\eta'\|\le\delta$. Then
\[
\begin{aligned}
\left\|\frac1n\sum_{i=1}^n\gamma_{\eta'}(Y_i)-b\right\|
\le\frac1n\sum_{i=1}^n
\|\gamma_{\eta'}(Y_i)-\gamma_\eta(Y_i)\|
\le L_M\delta.
\end{aligned}
\]
Consequently, by the union bound and
Lemma~\ref{lem:smallball},
\[
\begin{aligned}
\Pr\left\{Y\in K_M^n,\
\exists\eta\in C:
\frac1n\sum_{i=1}^n\gamma_\eta(Y_i)=b\right\}
\le
\sum_{\eta'\in C_\delta}
\Pr\left\{
\left\|\frac1n\sum_{i=1}^n
\gamma_{\eta'}(Y_i)-b\right\|
\le L_M\delta
\right\}\\
\le C'C_0L_M^m\delta^{m-p}.
\end{aligned}
\]
Since $m>p$, letting $\delta\downarrow0$ gives
\begin{equation}\label{eq:avoidance-zero}
\Pr\left\{Y\in K_M^n,\
\exists\eta\in C:
\frac1n\sum_{i=1}^n\gamma_\eta(Y_i)=b\right\}=0.
\end{equation}

Let
\[
E_M=
\left\{X\in K_M^n:
\exists\eta\in C,\quad
\frac1n\sum_{i=1}^n\gamma_\eta(X_i)=b
\right\}.
\]
The set
$
\{(X,\eta)\in K_M^n\times C:
n^{-1}\sum_{i=1}^n\gamma_\eta(X_i)=b
\}
$
is closed in the compact set $K_M^n\times C$.
Its projection onto $K_M^n$ is therefore compact,
so $E_M$ is measurable.

Since
$
c_M:=\inf_{x\in K_M}\rho_M(x)>0,
$
independence and \eqref{eq:avoidance-zero} give
\[
\begin{aligned}
0
=\int_{E_M}\prod_{i=1}^n\rho_M(x_i)\,
dx_1\cdots dx_n
\ge c_M^n|E_M|.
\end{aligned}
\]
Thus $|E_M|=0$. Finally,
\[
\left\{X\in(\mathbb R^d)^n:
\exists\eta\in C,\quad
\frac1n\sum_{i=1}^n\gamma_\eta(X_i)=b
\right\}
=\bigcup_{M=1}^\infty E_M,
\]
which has Lebesgue measure zero.
\end{proof}

We now introduce the functions needed to establish the generic nondegeneracy of the optimal score. Consider the $k$-atomic mixing measure 
\begin{equation}\label{eq:mixture}
\mu=\sum_{j=1}^k w_j\delta_{\theta_j},
\end{equation}
where 
\[
w_j>0,\text{ and } \sum_{j=1}^k w_j=1
\]
.Taking $w_1,\ldots,w_{k-1}$
as free weight parameters, with
$w_k=1-\sum_{j<k}w_j$, the parameter dimension is
\begin{equation}\label{eq:s}
s=(k-1)+kd=k(d+1)-1.
\end{equation}
For $x,z\in\mathbb R^d$, define
\begin{equation}\label{eq:T}
T_\mu(x,z)=\frac{\phi_d(x-z)}{p_\mu(x)},
\quad
E_z(x)=e^{\langle x,z\rangle-\|z\|^2/2},
\quad
D_\mu(x)=\sum_{j=1}^k w_jE_{\theta_j}(x).
\end{equation}
Since $D_\mu>0$, we have
$T_\mu(x,z)=E_z(x)/D_\mu(x)$. Differentiating with respect
to $z$, while keeping $\mu$ fixed, gives
\begin{equation}\label{eq:derivatives}
\begin{aligned}
\nabla_zT_\mu(x,z)
&=(x-z)T_\mu(x,z),\\
\nabla_z^2T_\mu(x,z)
&=\bigl[(x-z)(x-z)^\top-I_d\bigr]T_\mu(x,z).
\end{aligned}
\end{equation}

The first-order optimality conditions for a $k$-atomic
NPMLE involve the values and first derivatives of
$T_\mu(x,\theta_j)$. We therefore define
\begin{equation}\label{eq:base}
G_\mu(x)=
\left(
(T_\mu(x,\theta_j))_{j<k},
(\nabla_zT_\mu(x,\theta_j))_{j\le k}
\right)\in\mathbb R^s.
\end{equation}
Here the $k$th value coordinate is omitted since
$
\sum_{j=1}^k w_jT_\mu(x,\theta_j)=1.
$
To exclude additional contact points and degenerate
Hessians, we consider the following two maps:
\begin{align}
\gamma^{\mathrm{ext}}_{\mu,z}(x)
&=\left(G_\mu(x),T_\mu(x,z),
\nabla_zT_\mu(x,z)\right),
&&z\notin\{\theta_1,\ldots,\theta_k\},
\label{eq:extra-map}\\
\gamma^{\mathrm{hes}}_{\mu,v,j}(x)
&=\left(G_\mu(x),
\nabla_z^2T_\mu(x,\theta_j)v\right),
&&v\in\mathbb S^{d-1},\quad 1\le j\le k.
\label{eq:hessian-map}
\end{align}
Their respective dimensions are
\begin{equation}\label{eq:m-dimensions}
m_1=s+d+1=(k+1)(d+1)-1,
\quad
m_2=s+d=(k+1)(d+1)-2.
\end{equation}

\begin{lemma}[Independence and polynomial degrees]
\label{lem:features}
For either map $\gamma$ in
\eqref{eq:extra-map}--\eqref{eq:hessian-map},
with target dimension $m$, we have
\[
a\cdot\gamma(x)-c\equiv0
\quad\Longrightarrow\quad a=0,\quad c=0.
\]
Moreover, for every $a\in\mathbb S^{m-1}$ and
$c\in\mathbb R$, the function $a\cdot\gamma-c$
satisfies \eqref{eq:numerator} with $g=D_\mu>0$,
a nonzero right-hand side, and
\begin{equation}\label{eq:orders}
\begin{aligned}
&r=k+1,\quad q_\ell=1\ (1\le\ell\le k+1),
&&\gamma=\gamma^{\mathrm{ext}}_{\mu,z},\\
&r=k,\quad q_j=2,\quad q_\ell=1\ (\ell\ne j),
&&\gamma=\gamma^{\mathrm{hes}}_{\mu,v,j}.
\end{aligned}
\end{equation}
In particular, $J_1=2(k+1)$ and $J_2=2k+1$.
On compact subsets of the respective admissible parameter
spaces, both maps are jointly smooth and uniformly locally
Lipschitz in their parameters. Each compact parameter set
admits $\delta$-nets of cardinality
$O(\delta^{-(m-1)})$.
\end{lemma}

\begin{proof}
We first show that the coordinates of $G_\mu$, together
with the constant function, are linearly independent.
Suppose that
\[
c+\sum_{j<k}\alpha_jT_\mu(x,\theta_j)
+\sum_{j=1}^k\langle b_j,x-\theta_j\rangle
T_\mu(x,\theta_j)=0
\]
for every $x\in\mathbb R^d$. Multiplying by $D_\mu$ gives
\[
\sum_{j=1}^k
\left(cw_j+\mathbf1_{\{j<k\}}\alpha_j
+\langle b_j,x-\theta_j\rangle\right)
E_{\theta_j}(x)=0.
\]
Since the $\theta_j$ are distinct,
Lemma~\ref{lem:independence} gives
\[
cw_j+\mathbf1_{\{j<k\}}\alpha_j
+\langle b_j,x-\theta_j\rangle=0,
\quad 1\le j\le k.
\]
Comparing the linear coefficients gives $b_j=0$
for every $j$. Taking $j=k$ then gives $cw_k=0$.
Since $w_k>0$, we obtain $c=0$, and consequently
$\alpha_j=0$ for every $j<k$.

For \eqref{eq:extra-map}, suppose that
\[
\begin{aligned}
0={}&c+\sum_{j<k}\alpha_jT_\mu(x,\theta_j)
+\sum_{j=1}^k\langle b_j,\nabla_zT_\mu(x,\theta_j)\rangle
+\alpha T_\mu(x,z)
+\langle b,\nabla_zT_\mu(x,z)\rangle
\end{aligned}
\]
for every $x\in\mathbb R^d$. Multiplying by $D_\mu$ gives
\[
\begin{aligned}
0={}&\sum_{j=1}^k
\left(cw_j+\mathbf1_{\{j<k\}}\alpha_j
+\langle b_j,x-\theta_j\rangle\right)E_{\theta_j}(x)
+\left(\alpha+\langle b,x-z\rangle\right)E_z(x).
\end{aligned}
\]
Since $z\notin\{\theta_1,\ldots,\theta_k\}$,
Lemma~\ref{lem:independence} gives
$
\alpha+\langle b,x-z\rangle=0
$
and
$
cw_j+\mathbf1_{\{j<k\}}\alpha_j
+\langle b_j,x-\theta_j\rangle=0,
\ \forall 1\le j\le k.
$
Thus $\alpha=0$, $b=0$, and the first part of
the proof gives $c=0$, $\alpha_j=0$, and $b_j=0$
for every $j$.
Hence the coordinates of \eqref{eq:extra-map},
together with the constant function, are linearly
independent.

For \eqref{eq:hessian-map}, suppose that
\[
\begin{aligned}
0={}&c+\sum_{\ell<k}\alpha_\ell T_\mu(x,\theta_\ell)
+\sum_{\ell=1}^k
\langle b_\ell,\nabla_zT_\mu(x,\theta_\ell)\rangle
+\langle h,\nabla_z^2T_\mu(x,\theta_j)v\rangle
\end{aligned}
\]
for every $x\in\mathbb R^d$, where $h\in\mathbb R^d$.
By \eqref{eq:derivatives}, multiplication by $D_\mu$ gives
\[
\begin{aligned}
0=\sum_{\ell=1}^k\Bigl(
cw_\ell+\mathbf1_{\{\ell<k\}}\alpha_\ell
+\langle b_\ell,x-\theta_\ell\rangle
+\mathbf1_{\{\ell=j\}}
\bigl[
\langle h,x-\theta_j\rangle
\langle v,x-\theta_j\rangle-\langle h,v\rangle
\bigr]\Bigr)E_{\theta_\ell}(x).
\end{aligned}
\]
Since the $\theta_\ell$ are distinct,
Lemma~\ref{lem:independence} implies that
\[
\begin{aligned}
0={}&cw_j+\mathbf1_{\{j<k\}}\alpha_j
+\langle b_j,x-\theta_j\rangle
+\langle h,x-\theta_j\rangle
\langle v,x-\theta_j\rangle-\langle h,v\rangle
\end{aligned}
\]
for every $x\in\mathbb R^d$.
Comparing the homogeneous quadratic terms gives
$
\langle h,x\rangle\langle v,x\rangle=0,
\quad x\in\mathbb R^d.
$
Since $v\ne0$, we must have $h=0$.
Hence, we obtain
\[
c+\sum_{\ell<k}\alpha_\ell T_\mu(x,\theta_\ell)
+\sum_{\ell=1}^k
\langle b_\ell,\nabla_zT_\mu(x,\theta_\ell)\rangle=0
\]
for every $x\in\mathbb R^d$. By the first part of
the proof, $c=0$, $\alpha_\ell=0$ for $\ell<k$,
and $b_\ell=0$ for $1\le\ell\le k$.

Consequently, for either map $\gamma$, we have
\[
a\cdot\gamma(x)-c\equiv0
\quad\Longrightarrow\quad a=0,\quad c=0.
\]
In particular, for every $a\in\mathbb S^{m-1}$
and $c\in\mathbb R$,
\[
D_\mu(x)(a\cdot\gamma(x)-c)
\]
is not identically zero. The preceding expressions
also show that the polynomial coefficients have
degrees at most one at each of the $k+1$ distinct
exponents for \eqref{eq:extra-map}, and at most two
at $\theta_j$ and one at the remaining $k-1$
exponents for \eqref{eq:hessian-map}. Thus,
\[
J_1=2(k+1),\quad
J_2=3+2(k-1)=2k+1.
\]
Here the factors $e^{-\|z\|^2/2}$ in $E_z$
are absorbed into the polynomial coefficients.

Finally, let $K\subset\mathbb R^d$ and $C$ be
compact, with $C$ contained in the admissible
parameter space of either map. Since
$D_\mu(x)>0$ and is continuous,
\[
\inf_{(x,\eta)\in K\times C}D_\mu(x)>0.
\]
The expressions \eqref{eq:T}--\eqref{eq:hessian-map}
therefore extend smoothly to a neighborhood of
$K\times C$ in the free weight, location, and
vector coordinates. Their first derivatives with
respect to these coordinates are bounded on a
smaller neighborhood of $K\times C$. Hence,
for some $L_{K,C}<\infty$ and $\delta_0>0$,
\[
\|\gamma_\eta(x)-\gamma_{\eta'}(x)\|
\le L_{K,C}\|\eta-\eta'\|,
\]
whenever $x\in K$, $\eta,\eta'\in C$, and
$\|\eta-\eta'\|\le\delta_0$.

For \eqref{eq:extra-map}, the parameters belong
to a subset of $\mathbb R^{s+d}$, so every compact
parameter set admits $O(\delta^{-(s+d)})$
net points. For \eqref{eq:hessian-map}, the
parameters belong to a subset of
$\mathbb R^s\times\mathbb S^{d-1}$, which admits
$O(\delta^{-(s+d-1)})$ net points, since
$\mathbb S^{d-1}$ has a finite covering by
Lipschitz coordinate charts of dimension $d-1$.
For $d=1$, the sphere consists of two points.

By \eqref{eq:m-dimensions}, these exponents are
$m_1-1$ and $m_2-1$, respectively. This proves
the remaining assertions.
\end{proof}

\subsection{Proof of Theorem \ref{thm:generic-nondegeneracy}}
\begin{proof}[\textbf{Proof of Theorem \ref{thm:generic-nondegeneracy}}]
Let $\mu$ be a $k$-atomic global NPMLE. By the first-order optimality conditions, we have
\begin{equation}\label{eq:optimality}
\Psi_{\mu,X}(z)\le1,\quad
\Psi_{\mu,X}(\theta_j)=1,\quad
\nabla\Psi_{\mu,X}(\theta_j)=0,
\end{equation}
for every $z\in\mathbb R^d$ and $1\le j\le k$.
Next, for
$
b_0=(\mathbf1_{k-1},\mathbf0_{kd})\in\mathbb R^s,
$
equations \eqref{eq:optimality} give
\begin{equation}\label{eq:base-equations}
\frac1n\sum_{i=1}^nG_\mu(X_i)=b_0.
\end{equation}

For every integer $M\ge\max\{2,k\}$, let $C_M$
be the set of parameters in \eqref{eq:mixture} satisfying
\begin{equation}\label{eq:compact-mixtures}
w_j\ge M^{-1},\quad
\|\theta_j\|\le M,\quad
\|\theta_j-\theta_h\|\ge M^{-1}
\quad(j\ne h).
\end{equation}
These sets are compact in the $s$-dimensional parameter
space, and every $k$-atomic probability measure belongs
to some $C_M$. Moreover, each $C_M$ admits
$O_M(\delta^{-s})$ net points.

Suppose first that $\Psi_{\mu,X}$ has a global
maximizer $z\notin\{\theta_1,\ldots,\theta_k\}$.
For sufficiently large $M$, we have
\[
\mu\in C_M,\ \ 
\|z\|\le M,\ \ 
\min_{1\le j\le k}\|z-\theta_j\|\ge M^{-1}.
\]
Let $C_M^{\mathrm{ext}}$ denote the set of parameters
$(\mu,z)$ satisfying these inequalities.
It is compact and admits $O_M(\delta^{-p_1})$
net points, where $p_1=s+d=m_1-1.$

Since $z$ is a global maximizer,
$\Psi_{\mu,X}(z)=1, \text{and }\nabla\Psi_{\mu,X}(z)=0.
$
Together with \eqref{eq:base-equations}, this gives
\begin{equation}\label{eq:extra-equations}
\frac1n\sum_{i=1}^n
\gamma^{\mathrm{ext}}_{\mu,z}(X_i)
=(b_0,1,\mathbf0_d)\in\mathbb R^{m_1}.
\end{equation}
By Lemma~\ref{lem:features}, the map
$\gamma^{\mathrm{ext}}_{\mu,z}$ satisfies
Lemma~\ref{lem:fourier} with $J_1=2(k+1)$.
The assumption \eqref{eq:generic-sample-size} gives
\[
n>(2k+1)m_1=m_1(J_1-1).
\]
Since $p_1<m_1$, Lemmas~\ref{lem:smallball}
and \ref{lem:avoidance} imply that the set of
$X\in(\mathbb R^d)^n$ for which
\eqref{eq:extra-equations} holds for some
$(\mu,z)\in C_M^{\mathrm{ext}}$ has Lebesgue
measure zero. Taking the countable union over
$M$ shows that, outside a Lebesgue-null set,
every $k$-atomic global NPMLE satisfies
\[
\text{argmax}_{z\in\mathbb R^d}\Psi_{\mu,X}(z)
=\{\theta_1,\ldots,\theta_k\}.
\]

Next, fix $1\le j\le k$. Since $\theta_j$
is a global maximizer,
$
\nabla^2\Psi_{\mu,X}(\theta_j)\preceq0.
$
Suppose that this matrix is not negative definite.
It is then singular, so there exists
$v\in\mathbb S^{d-1}$ such that
$
\nabla^2\Psi_{\mu,X}(\theta_j)v=0.
$
For sufficiently large $M$, the parameters
$(\mu,v)$ belong to
\[
C_M^{\mathrm{hes}}=C_M\times\mathbb S^{d-1}.
\]
This set is compact and admits
$O_M(\delta^{-p_2})$ net points, where
$
p_2=s+d-1=m_2-1.
$
For $d=1$, the sphere consists of two points,
and the same bound holds with $d-1=0$.

By \eqref{eq:base-equations} and the definition
of \eqref{eq:hessian-map},
\begin{equation}\label{eq:hessian-equations}
\frac1n\sum_{i=1}^n
\gamma^{\mathrm{hes}}_{\mu,v,j}(X_i)
=(b_0,\mathbf0_d)\in\mathbb R^{m_2}.
\end{equation}
By Lemma~\ref{lem:features}, the map
$\gamma^{\mathrm{hes}}_{\mu,v,j}$ satisfies
Lemma~\ref{lem:fourier} with $J_2=2k+1$.
Since $m_2=m_1-1$, we have
\[
m_2(J_2-1)
=2k(m_1-1)
<(2k+1)m_1<n.
\]
Moreover, $p_2<m_2$. Thus, by
Lemmas~\ref{lem:smallball} and
\ref{lem:avoidance}, the set of datasets
satisfying \eqref{eq:hessian-equations} for
some $(\mu,v)\in C_M^{\mathrm{hes}}$
has Lebesgue measure zero.

Taking the union over $1\le j\le k$ and
all integers $M\ge\max\{2,k\}$, we conclude
that, outside a Lebesgue-null set,
\[
\nabla^2\Psi_{\mu,X}(\theta_j)\prec0,
\quad 1\le j\le k,
\]
for every $k$-atomic global NPMLE $\mu$.

The union of the two exceptional sets has
Lebesgue measure zero. Hence both conclusions
in \eqref{eq:generic-conclusion} hold simultaneously
for every $k$-atomic global NPMLE outside
this set. If the joint distribution of
$(X_1,\ldots,X_n)$ is absolutely continuous,
the exceptional set has probability zero.
\end{proof}


\end{document}